\documentclass[12pt]{amsart}
\usepackage{amssymb,amscd}  
\usepackage{pb-diagram}   
\usepackage{verbatim}

\usepackage[dvipdfmx]{graphicx,color}  

\usepackage{graphicx}
\usepackage{mathrsfs}

\usepackage{here}   

\usepackage[english]{babel}
\usepackage{amsfonts,amssymb,amsmath,amstext,amsbsy,amsopn,amscd,amsthm,graphicx,euscript,graphicx}
\usepackage{graphics}
\usepackage[all]{xy}

\def\X{{\mathcal X}}

\DeclareMathOperator{\sign}{sign}

\DeclareMathOperator{\im}{im}

\newtheorem*{Whitney towers}{Theorem~\ref{Whitney towers}}
\newtheorem*{h-towers}{Theorems ~\ref{half} \& \ref{$(n)$-solvable}}

\newtheorem*{surgery curves}{Theorem~\ref{surgery curves}}
\newtheorem*{cg=0}{Theorem~\ref{vanish}}

\newtheorem{thm}{Theorem}[section]
\newtheorem{mth}[thm]{Main Theorem}
\newtheorem{proposition}[thm]{Proposition} 
\newtheorem{lem}[thm]{Lemma}

\newtheorem{fact}[thm]{Fact}

\newtheorem{cla}[thm]{Claim}

\theoremstyle{definition}
\newtheorem{definition}[thm]{Definition}
\newtheorem{example}[thm]{Example}
\newtheorem{rev}[thm]{Review}
\newtheorem{construction}[thm]{Construction}
\newtheorem{property}[thm]{Property}

\newtheorem{que}[thm]{Question}
\newtheorem{remark}[thm]{Remark}

\newtheorem{assumption}[thm]{Assumption}

\numberwithin{equation}{section}
\numberwithin{figure}{section}
\numberwithin{table}{section}

\newcommand{\vs}{\vskip10mm}
\newcommand{\bb}{\bigbreak}
\newcommand{\bs}{\smallbreak}
\newcommand{\h}{\noindent}
\newcommand{\x}{\times}
\newcommand{\np}{\newpage}
\newcommand{\Z}{\mathbb{Z}}
\newcommand{\N}{\mathbb{N}}
\newcommand{\C}{\mathbb{C}}

\newcommand{\R}{\mathbb{R}}

\def\yen{{\setbox0=\hbox{Y}Y\kern-.97\wd0\vbox{hrule height.lex width.98%
\wd0\kern.33ex\hrule height.lex width.98\wd0\kern.45ex}}}

\def\np{\newpage}

\begin{document}
{

\vspace*{-50pt}

\pagestyle{plain}




\h{\bf 
Steenrod square for virtual links toward Khovanov-Lipshitz-Sarkar stable homotopy type for virtual links} 

\bb
\bb

Louis H. Kauffman and  Eiji Ogasa 

\bb
\bb

\h {\bf Abstract.}
We define a second Steenrod square  for virtual links, 
 which is 
stronger than Khovanov homology for virtual links, 
 toward  constructing Khovanov-Lipshitz-Sarkar
stable homotopy type for virtual links.
This induces the first meaningful nontrivial example of 
the second Steenrod square operator on the Khovanov homology for 
links in a 3-manifold other than $S^3$. 

 \tableofcontents

\part{Introduction}\label{partintro}
\bigbreak
\section{Main result}\label{vk}
\bb

\noindent 
In this paper we define a second Steenrod square operator acting on Khovanov homology for virtual links. 
This induces the first meaningful nontrivial example of 
the second Steenrod square operator on the Khovanov homology for 
links in a 3-manifold other than $S^3$ (see \S\ref{tuke}). 
Our work extends the work of Lipshitz and Sarkar \cite{LSk, LSs} and Seed \cite{Seed} 
(See
the part 
above Theorem \ref{LSSe} 
and Theorem \ref{LSSe}
in \S\ref{seclo}).
It also raises many new questions about extending more fully their work on 
Khovanov-Lipshitz-Sarkar stable homotopy type. 
The list of contents for the paper given above shows the structure of this paper. 
In \S\ref{secv1k} 
we review the basics of virtual knot theory.
In this paper we consider Khovanov homology for vrtual links defined by Manturov \cite{Man}. 
We review its definition in  \S\ref{seclo}, which is  self-contained.  
See the part below Theorem \ref{LSSe} for other versions of  Khovanov homology for vrtual links. 
In \S\ref{tuke}  we explain 
why virtual knot theory is important for research on the Jones polynomial, Khovanov homology, and Khovanov homotopy, and 
describe key relationships of virtual knot theory with problems 
in low dimensional topology.
\\

We state our main theorem,  which follows from 
 ``Theorem \ref{main} with Definition \ref{ofuro}'',   
Theorem \ref{YY}, and Theorem \ref{oreikhalkos}.

\begin{mth}\label{ogon}
$(1)$We define the second Steenrod square operator for virtual links.  

\smallbreak\noindent$(2)$
If $L$ is a classical link diagram, 
our second Steenrod square  is the same 
as the second Steenrod square  in the case of classical links which is defined in \cite{LSs}.

\smallbreak\noindent$(3)$  
In the case of non-classical, virtual links, the second Steenrod square  is stronger than 
Khovanov homology. 
%
That is, there is a pair of non-classical, virtual links 
that  have different second Steenrod square operators while they have the same Khovanov homology. 
\end{mth}

The paper consists of two main parts. 
Part \ref{partKhh} is a review and technical definition of Khovanov homology of links and virtual links. 
If the reader is already familiar with these theories, this part will be useful for reference. 
It is important for the reader to understand that single-cycle resmoothings can occur in virtual link theory and that these lead to the differences between classical and virtual Khovanov homology. 
Part \ref{partsht} consists in the work needed to prove the main results of the paper and uses the constructions in Part \ref{partKhh}.

\section{Making a CW complex from a given chain complex
}\label{secffc}
We will make a CW complex from a Khovanov chain complex for a  given virtual link diagram.
However we begin by discussing chain complexes 
that are not necessarily a Khovanov chain complex for classical links or virtual links.

Suppose that we are gven a chain complex whose basis is a set $\{g_i\}$, 
where $\{i\}$ is a finite set. 
Note that the given chain complex is not necessarily a Khovanov chain complex for classical links or virtual links.  
Is there a CW complex whose chain complex is the given chain complex?
The answer is positive. The proof is easy 
(See \cite[Theorem in Exercise 4,  section 39, page 231]{Munkres}.). 
However, there may be more than one CW complex whose chain complex is the given one.
We want to give a specific stable homotopy type to a given chain complex.
It is an important step in our procedure. We summarize our idea below.

\begin{example}\label{exam0}
Both 
the one point union $S^2\vee S^4$ and $\C P^2$ 
have a natural CW decomposition 
 (the base point)$\cup e^{2}\cup e^{4}$. 
Both have the same chain complex and 
 different stable homotopy types. 
\end{example}

If we specify maps to $\partial e^4$ to the 2-skeleton   (the base point)$\cup e^{2}$ 
as follows, we determine which the resulting 4-skeleton is 
$S^2\vee S^4$ or $\C P^2$. 
Take the trivial knot 
$S^1$ in the 3-sphere $\partial e^4$. \\
\bs\h(1)  We give a framing on the normal bundle of $S^1$ in  $\partial e^4$
so that the framing extends to an embedded disc that bounds $S^1$. 
\\\bs\h(2)  We give it as follows;  
 	A framing consists locally of two orthonormal vectors (arrows). 
Let the normal frame to the circle so that the original circle and either circle traced by one of the arrows makes the Hopf link.\\
\bs\h
Use the Pontrjagin-Thom construction associated with (1) and (2).
Then the condition (1) defines $S^2\vee S^4$, 
and the condition (2) induces $\C P^2$. 
Recall a relation between the Hopf fibration 
and the Pontrjagin-Thom construction associated with (2).

The circle $S^1$ in $\partial e^4$ is an example of modulis. 
The term moduli is defined in LS and framed moduli are the basis for constructing CW complexes by Lipshitz and Sarkar
 (See \cite[\S3.3 and \S4]{LSk} and \cite[\S3.3]{LSs}.). 
We direct the reader to LS for the full definition of moduli and point out that the term is related to equivalence classes of points in the flow related to a Morse singularity. The flow concept is generalized in \cite{LSk} to a flow category. We give the definition of flow category below and again direct the reader to \cite{LSk} for the details. The examples we have cited are useful for understanding the general principle of the framed moduli.
Those normal framings on $S^1$ in the conditions (1) and (2) are examples of framings on modulis. 
The method of modulis and framings is obtained by generalizing the Pontrjagin-Thom construction. 

Note that, for any given CW complex, in each cell, we can define modulis and framings. 

We construct as follows.

\begin{construction}\label{constkore}
We are given a chain complex and a basis. The basis is a poset by the differential.
We attach the lowest dimensional cell to the base point. 
Note that if it is  0-dimensional cell, the attaching map is the empty. 
From lower dimensional cells to higher dimensional cells, step by step, 
we attach each cell to the lower dimensional skeleton as follows. 
When we attach each cell $e$, 
 we define a moduli $\mathcal M$ in $\partial e$ and a framing on $\mathcal M$. 
 So $\partial\mathcal M\cap$(the lower dimensional skeleton)  has been framed. 
 We must define the framing on $\mathcal M$ 
 so that the given framing on $\partial\mathcal M$ extends to $\mathcal M$. 

In other words, 
we define modulis and framings in lower dimensional cells 
to those in higher dimensional ones.
\end{construction}

Construction \ref{constkore} does not complete for any basis of any chain complex. 
We show an example below. 
We use 
\cite[Definition 3.12.(M-2)]{LSk} crucially, 
so we cite it here.


\bb

{\bf Definition 3.12 of \cite{LSk}.}\label{pageda}
A {\it flow category} is a pair 
$(\mathscr C,$ gr$)$ where 
$\mathscr C$ is a category with finitely many
objects Ob = Ob$(\mathscr C)$ and 
gr : Ob $\to\Z$ is a function, called the grading, satisfying the
following additional conditions:

\begin{enumerate}
\item[(M-1)]  
Hom$(x, x)=\{$Id$\}$ for all 
$x\in$ Ob, and for distinct 
$x, y\in$ Ob, Hom$(x, y)$ is a compact
(gr$(x)-$gr$(y)-1)$-dimensional 
$<$gr$(x)-$gr$(y)-1>$-manifold (with the understanding
that negative dimensional manifolds are empty).

\item[(M-2)]  
For distinct $x, y, z\in$Ob with 
gr$(z)-$gr$(y)=m$, the composition map

$$\circ: \text{Hom}(z, y)\x \text{Hom}(x, z)\to\text{Hom}(x, y)$$ 

\noindent 
is an embedding into $\partial_m \text{Hom}(x, y)$. 
Furthermore,

$$\circ^{-1}(\partial_i \text{Hom}(x, y)) =
\begin{cases}
\partial_i \text{Hom}(z, y)\x Hom(x, z) & \text{for $i<m$} \\
\text{Hom}(z, y)\x\partial_{i-m} \text{Hom}(x, z) & \text{for $i>m$} 
\end{cases}$$

\item[(M-3)] 
For distinct $x, y\in$Ob, 
$\circ$  induces a diffeomorphism
\vskip2mm

\noindent
(3.1) $\partial_i \text{Hom}(x, y)\cong
\displaystyle\cup_{
z\hskip1mm 
gr(z)=gr(y)+i\hskip1mm } 
\text{Hom}(z, y)\x\text{Hom}(x, z).$
\end{enumerate}

\bigbreak
\bigbreak

Therefore, if $\mathcal D$ is the diagram whose vertices are the spaces 
$$\text{Hom}(z_m, y)\x\text{Hom}(z_{m-1}, z_m)\x\cdot\cdot\cdot\x 
\text{Hom}(x, z_1)$$ 
for $m\geqq1$ and distinct 
$z_1, . . . , z_m\in$ Ob $-\{x, y\}$, 
and whose arrows correspond to composing a single adjacent pair of Hom's, then
$$\partial\text{Hom}(x, y)\cong
\cup_i
\partial_i\text{Hom}(x, y)\cong
\text{colim}\mathcal D.$$ 
Given objects $x, y$ in a flow category $\mathscr C$, 
define the {\it moduli space from $x$ to $y$} to be
$$\mathcal M(x, y) =
\begin{cases}
\phi& \text{if x=y}\\
\text{Hom$(x, y)$} &\text{otherwise.}
\end{cases}$$ 

Given a flow category $\mathscr C$ and an integer $n$, 
let $\mathscr C [n]$ be the flow category obtained from
$\mathscr C$ by increasing the grading of each object by $n$.

\bb





\begin{example}\label{exae5}
Take a CW decomposition (the base point)$\cup e^{2}\cup e^{4}$ of $\C P^2$. 
Consider a question whether we can attach a 5-cell $e^5$ 
to $\C P^2$ so that $H_4(\C P^2\cup e^5)=0$. 
The answer is negative because $\pi_4\C P^2=0$. 
We show an alternative proof by using modulis and framings and Definition 3.12.(M-2).
We use reductio ad absurdum. 
We assume that we can attach $e^5$ to $\C P^2$ so that $H_4(\C P^2\cup e^5)=0$. 
Then there is a moduli $\mathcal M(e^5, e^2)$ in $e^5$. 
By Definition 3.12.(M-2), 
  $\partial\mathcal M(e^5, e^2)=\mathcal M(e^5, e^4)\x\mathcal M(e^4, e^2)$.
Since $H_4(\C P^2\cup e^5)=0$, 
$\mathcal M(e^5, e^4)$ is one point and the framing on the point is determined by the differential. 
$\mathcal M(e^4, e^2)$ is a circle and 
the framing on the circle defined as in  the above condition (2).  
The framing on $\partial\mathcal M(e^5, e^2)$ 
is induced from that on $\mathcal M(e^5, e^4)$ and that on $\mathcal M(e^4, e^2)$ 
(See \cite[Definition 3.18]{LSk}).    
This framing on $\partial\mathcal M(e^5, e^2)$ 
cannot extend to $\mathcal M(e^5, e^2)$. 
We arrived at a contradiction. 
\end{example}



We can realize Construction \ref{constkore} 
 if a given chain complex has the following prpperty.  

\begin{property}\label{propertyFM}
Suppose that we are given a basis of a chain complex.  
Each basis element is represented by a cell. 
Suppose that we can attach a cell $e$ to the llower dimensional skeleton. 
Then we have a moduli $\mathcal M$ in $\partial e$. 

We have the following important property: Not all sub-modulis in $\partial \mathcal M$ 
have been framed in the lower dimensional skeleton $X$.  
Furthermore we can  extend the framing on $\partial\mathcal M$
that has been framed in $X$, to $\mathcal M$.  
\end{property}

 Khovanov chain complexes for classical links 
 always   have Property \ref{propertyFM}  (See \cite[Proposition 4.12 and Definition 5.5]{LSk}). 
In this paper we will prove that a partial case of  the virtual case has Property \ref{propertyFM}.  
So Property \ref{propertyFM} is a key of our construction of Khovanov partial CW complex.
Note that   there are many chain complexes that satisfy Property \ref{propertyFM} 
and that are not a Khovanov chain complex for classical or virtual links.  
Note that Example \ref{exae5} does not satisfy Property \ref{propertyFM}.

\section{Strategy of the construction of CW complexes for virtual links}\label{secstr}

 Lipshitz and Sarkar 
 \cite{LSk} 
 defined 
 modulis for all $P(e,g)$, where $e$ and $g$ run over all basis elements, 
 and introduced framings on them, 
 and constructed a stable homotopy type of a CW complex  for a given classical link.

In the virtual case, it is very complicated to construct a  moduli $P(e,g)$ 
for a pair $e$ and $g$ such that 
the difference of the homological gradings is greater than four. 
The reason is that the property of coefficients in the virtual case (Definition \ref{korekoso}) 
is different from the classical case (see also \S\ref{secs1}). 
In this paper  we show an explicit way to assign to the moduli space 
when the difference  {\rm gr}$_h \bold x-${\rm gr}$_h \bold y$ of the homological gradings is no greater than four. 
We construct a CW complex which consists of 
only 
$(m-1)$-cells, 
 $m$-cells, 
$(m+1)$-cells, 
$(m+2)$-cells, and 
$(m+3)$-cells, 
where $m$ is any integer,    
for the dual Khovanov chain complex in this case by using these moduli spaces. 
We prove that the second Steenrod square of the CW complex is invariant under any Reidemeister move 
although we do not prove 
whether the stable homotopy type of the CW complex is invariant under any Reidemeister move.

\begin{remark}\label{remiron}
	In Lipshitz and Sarkar's paper \cite[\S5.4]{LSk} it is important how one assigns a moduli to the ladybug configuration.
(We review the ladybug configuration in \S\ref{lady} and \S\ref{qlady}.).  
Lipshitz and Sarkar consider all ladybug configurations in a given classical link diagram 
and show that they can make the specific choice either the  the right or the left pair. That is, all ladybug configurations are chosen to have the same pair.
In our paper, we consider the following situation: 
Each ladybug configuration in a given virtual link diagram 
may have different pairs. 
We consider all cases. Therefore we may consider more than one set of modulis for one Khovanov chain complex. 
Note that classical link diagrams are virtual link diagrams. 
%
We discuss 
both 
the condition that Lipshitz and Sarkar\cite{LSk} discussed 
and 
the condition that they 
did not discuss.
Thus, in the case of classical links, 
we also handle modulis that they 
did not handle, 
but our second Steenrod square on a classical link is the same as theirs.
\end{remark}

\bigbreak
\part{Review of Khovanov homology for virtual links}\label{partKhh}

\section{Virtual knots and virtual links}\label{secv1k}

We work in the smooth category. 
Let $L=\{L_1,...,L_m\}$ be a 1-dimensional submanifold of a connected 3-manifold $M$, 
where $L_{k}$ denotes a connected component. 
Then $L$ is called an {\it $m$-component link} in $M$. 
In this paper we only discuss 1-dimensional links, 
so we say link, not 1-dimensional link nor 1-link.
Any 1-component link $L$ in $M$ is called a {\it knot} in $M$. 
If $M$ is the 3-sphere $S^3$, any link (respectively, knot) $L$ in $M$ is called a {\it classical link} (respectively, {\it classical knot}).

The theory of {\it virtual knot}s  was introduced in \cite{Kauffman1,Kauffman, Kauffmani} 
as a generalization of classical knot theory, 
and  
studies the embeddings of circles in thickened 
oriented closed surfaces
modulo isotopies and orientation preserving diffeomorphisms
plus one-handle stabilization of the surfaces.   \\

By a one-handle stabilization, 
we mean a surgery on the surface that is performed on a curve 
in the complement of the link embedding and 
that either increases or decreases the genus of the surface. 
The reader should note that knots and links in thickened surfaces can be represented by diagrams on the surface in the same sense as link diagrams drawn in the plane or on the two-sphere. 
From this point of view, a one handle stabilization is obtained by cutting the surface along a curve in the complement of the link diagram and 
capping the two new boundary curves with disks, or taking two points on the surface in the link diagram complement and cutting out two disks, and then adding a tube between them. 
The main point about handle stabilization is that it allows the virtual knot to be eventually placed in a least genus surface in which it can be represented. 
A theorem of Kuperberg \cite{Kuperberg} asserts that such minimal representations are topologically unique. \\

\begin{figure}
\centering
\includegraphics[width=30mm]{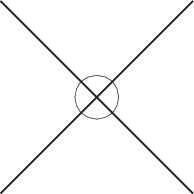}
\caption{{\bf Virtual crossing point}\label{vcro0}}   
\end{figure}

Virtual knot theory has a diagrammatic formulation.
A {\it virtual knot} can be represented by a {\it virtual knot diagram} 
in $\R^2$ (respectively, $S^2$) 
containing a finite number of real crossings, and {\it virtual crossings} 
indicated by a small circle placed around the crossing point as shown in Figure \ref{vcro0}.   
A virtual crossing is a combinatorial structure that shows how the diagram for the virtual knot is connected, and allows the reconstruction of a surface in which the knot has a diagram, up to handle stabilization.

The moves on virtual knot diagramsin $\R^2$ are generated by the usual Reidemeister moves plus the {\it detour move}.
The detour move allows a segment with a consecutive sequence of virtual crossings
to be excised and replaced by any other such a segment with a consecutive virtual
crossings, as shown in Figure \ref{detour}.    

Virtual 1-knot diagrams $\alpha$ and $\beta$ are changed into each other 
by a sequence of the usual Reidemeister moves  and detour moves 
if and only if 
$\alpha$ and $\beta$ are changed into each other 
by a sequence of all Reidemeister moves   
drawn in Figure \ref{all-1}. 
\\

\begin{figure}
\includegraphics[width=110mm]{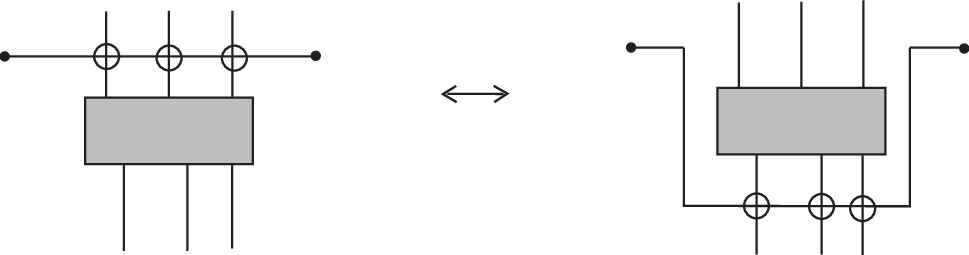}
\caption{{\bf An example of detour moves}\label{detour}}   
\end{figure}

\begin{figure}
\includegraphics[width=140mm]{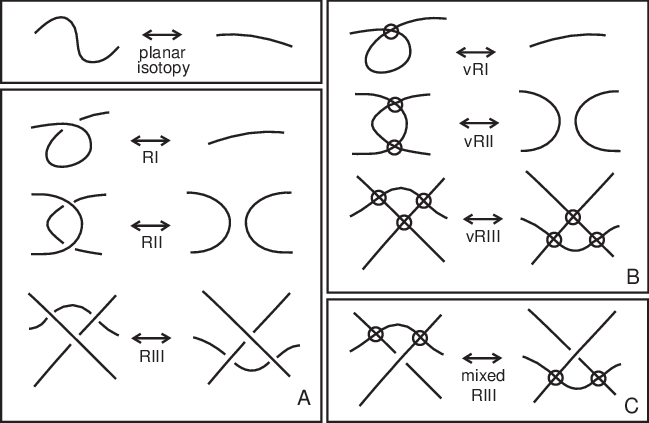}
\caption{{\bf All Reidemeister moves}   
\label{all-1}}   
\end{figure}

Virtual knot and link diagrams that can be related to each other by a finite sequence of the
Reidemeister and detour moves are said to be
{\it virtually equivalent} or {\it virtually isotopic}. \\
The virtual isotopy class of a virtual knot diagram is called a {\it virtual knot}.\\

There is a one-to-one correspondence between the topological and the diagrammatic approach
to virtual knot theory. The following theorem providing the transition between the
two approaches is proved by abstract knot diagrams, see  \cite{Kauffman1,Kauffman, Kauffmani}.\\


\begin{thm}\label{kihon} {\bf (\cite{Kauffman1,Kauffman, Kauffmani})} 
Two virtual link diagrams are virtually isotopic if and only if their surface embeddings are equivalent up to isotopy in the thickened surfaces, orientation preserving diffeomorphisms of the surfaces, and the  addition/removal of empty handles. 
\end{thm}

\begin{figure}
     \includegraphics[width=120mm]{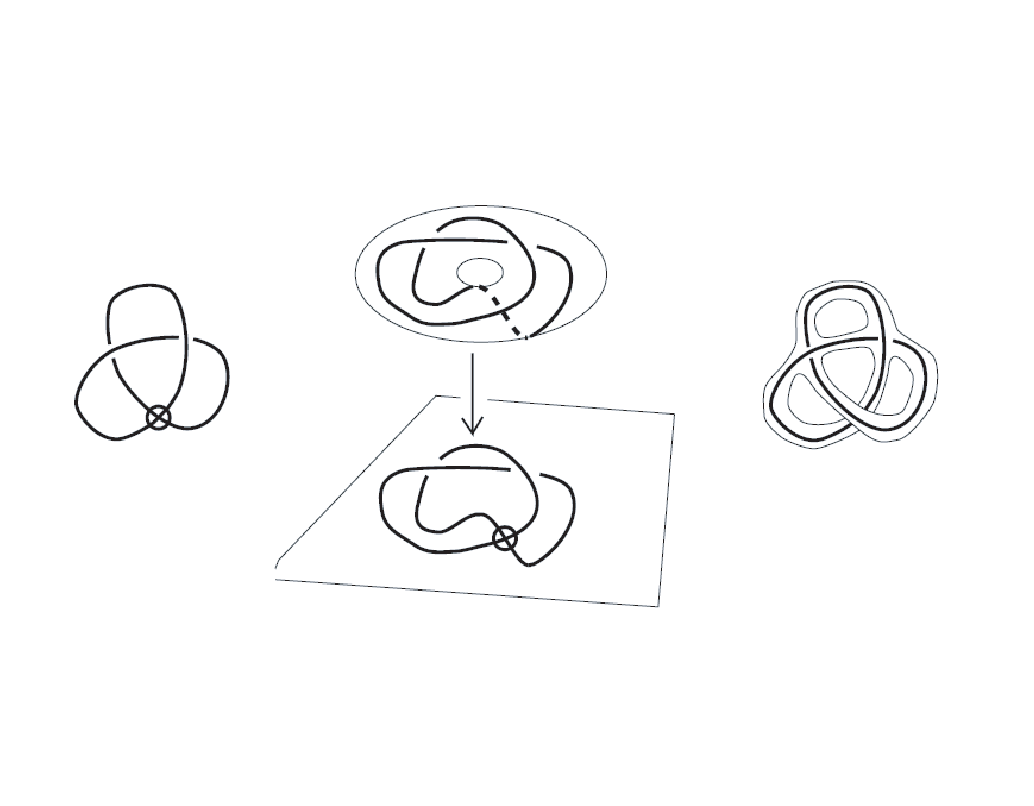}
\caption{{\bf 
How to make a representing surface from 
the tubular neighborhood of  a virtual knot diagram in $\R^2$
}\label{vtube}}   
\bigbreak  
\end{figure}

\noindent{\bf Remark.}  
A handle is said to be {\it empty} if the knot diagram does not thread through the handle.  One way to say this more precisely is to model the addition of and removal of handles via the location of surgery curves in the surface that do not intersect the knot diagram.
Here, an oriented surface with a  link diagram using only classical crossings 
appears. 
This surface is called a {\it representing surface}.   
In Figure \ref{vtube}
we show an example of a way to make a representing surface from a virtual knot diagram. 
Take the tubular neighborhood of a virtual knot diagram in $\R^2$. 
Near a virtual crossing point, double the tubular neighborhood. 
Near a classical crossing point, keep    the tubular neighborhood and the classical crossing point. Thus we obtain a compact representing surface with non-vacuous boundary.  
%
 We may start with  a representing surface that is oriented and not closed, 
 and then  embed the surface in a closed oriented surface to obtain a new representing surface. 
 Taking representations of virtual knots up to such cutting (removal of exterior of neighborhood of the diagram in a given surface) and re-embedding, plus isotopy in the given surfaces, corresponds to a unique diagrammatic virtual knot type.
\\

%
%
%

\bb

The presence of the single cycle zero map in the virtual Khovanov chain complex 
and 
the properties of the coefficients in the virtual Khovanov chain complex  
take us out of the cube complex method for defining a framed flow category. For this reason we use a truncated homotopy type for this paper.



\bigbreak
\section{Khovanov homology for virtual links}\label{seclo}

In this section we give a detailed version of  the Manturov method for defining Khovanov homology for virtual links. The reader interested in this definition should examine previous papers by 
Manturov \cite{Man}  and by   
Dye, Kaestner and Kauffman \cite{DKK} 
and also by  
Nikonov \cite{Igor}. 
The definition has a necessary complexity due to the presence of single cycle resmoothings in the Kauffman states of virtual link diagrams. In the papers mentioned there are different motivating points of view given for the definitions. Here we use these ideas and give a strict defintion that will require some work on the part of the reader to absorb. The advantage of this definition is that it permits specific formulas for the coefficients in the boundary formulas for the chain complex and it has some other advantages that are of a technical nature. The reader already familiar with Khovanov homology for virtual links should review what he knows before reading this section.

\bb


We begin by reviewing the case of classical links.   
Let $\mathcal L$ be a classical link. 
Let $L$ be a classical link diagram which represents $\mathcal L$. 
In \cite{K} Khovanov defined a chain complex for $L$, and 
proved that the homology defined by the chain complex is an invariant of the link type of $\mathcal L$. 
We call  this chain complex the {\it Khovanov chain complex}, and  
this homology {\it Khovanov homology}.  
Khovanov proved that 
%
%
%
%
the Jones polynomial of any classical link $L$ is a graded Euler characteristic of the Khovanov homology of $L$. 
In \cite{B} Bar-Natan reformulated the definition of Khovanov homology, and  
proved, by direct calculation,  that Khovanov homology is stronger than the Jones polynomial.\\

The differential on each element in the Khovanov chain complex
increases the degree of the element.  
This is not an ordinary convention but this terminology has been used for a long time. 
So we adopt this notation in this paper. Furthermore we define the following notations.  
We use the dual of the Khovanov chain complex in this paper. 
The differential on it decreases the degree. 
We call this chain complex the {\it dual Khovanov chain complex}.   
We call the homology defined by the dual Khovanov chain complex, 
the {\it homology of the dual Khovanov chain complex}. 
We use often the terminologies, ``chain complex'', ``chain homotopy'', and ``homological'',  
in both the case of   
Khovanov chain complexes and 
the case of the dual Khovanov chain complexes. 

Lipshitz and Sarkar \cite{LSk} introduced a stable homotopy type of CW complexes for 
each classical link, 
and proved that the stable homotopy type  is a topological invariant of classical links. 
We call it 
{\it Khovanov-Lipshitz-Sakar stable homotopy type} for classical links. 
Stable homotopy types of CW complexes have Steenrod squares. 
  Seed \cite{Seed} 
 proved the following fact by making a computer program 
according to Lipshitz and Sarkar's method in \cite{LSs}.

\begin{thm}\label{LSSe}  {\bf (\cite{Seed}.)}
The second Steenrod square for classical links is stronger than 
Khovanov homology for classical links.  
That is, there are classical links $K_1$ and $K_2$ such that 
 the Khovanov homology of $K_1$ is the same as that of $K_2$, 
but the second Steenrod square  of $K_1$ is different from that of $K_2$. 
Therefore 
Khovanov-Lipshitz-Sakar stable homotopy type for classical links 
is stronger than Khovanov homology for classical links. 
\end{thm}

We start the discussion on the case of virtual links now.  
Manturov defined Khovanov homology for virtual links in \cite{Man} (arXiv 2006). 
Rushworth \cite{Ru} and Tubbenhauer \cite{Tub} defined it in different methods. 
See also Viro \cite{Viro}.
Dye, Kaestner, and  Kauffman \cite{DKK} made an alternative definition of \cite{Man}. 
Nikonov \cite{Igor} described 
an alternative definition of \cite{DKK,Man}. 

We review the definition of Khovanov chain complexes and 
that of Khovanov homology for virtual links defined by Manturov \cite{Man}. 
In this paper, Khovanov homology for virtual links 
means Manturov's Khovanov homology for virtual links. 
Our exposition below of virtual Khovanov homology is self-contained. \\


In this paper we will generalize the result about the Steenrod square operator on 
Khovanov homology for classical links in \cite{LSr}  
to the virtual link case, 
by extending many results and methods in \cite{LSk, LSr}. 
So we first 
explain the mode of
definition for Khovanov chain complexes and 
Khovanov homology for virtual links 
in the fashion of \cite[section two]{LSk}. 

\begin{definition}\label{2.1} 
A {\it resolution configuration} $D$ is a pair $(Z(D),A(D))$, 
where $Z(D)$ is a set of pairwise-disjoint immersed circles in $S^2$, 
and $A(D)$ is a totally ordered collection of disjoint arcs embedded in $S^2$, 
with $A(D)\cap Z(D)=\partial A(D)$.  
We call the number of arcs in $A(D)$ the {\it index} of the resolution configuration $D$, 
and denote it by ind$(D)$. 
We sometimes abuse notation and write $Z(D)$ 
to mean $\cup_{Z\in Z(D)} Z$ and $A(D)$ to mean $\cup_{A\in A(D)} A$.
Occasionally, we will describe the total order on $A(D)$ 
by numbering the arcs: a lower numbered arc precedes a higher numbered one. 

We sometimes call an element of $Z(D)$, an immersed circle, a circle,  a loop or a component of $D$. 
\end{definition}

Note that 
resolution configurations are the same 
as what many people call {\it Kauffman state}s,   
which Kauffman first introduced in \cite{Kauffmanstate}. 
(The ways to draw arcs 
in both papers are different.)
\\

\begin{definition}\label{2.2} 
Given a virtual link diagram $L$ with $n$ classical crossings, 
an ordering of the crossings in $L$, 
and a {\it vector} $v\in\{0,1\}^n$ 
there is an {\it associated resolution configuration} $D_L(v)$ 
obtained by taking the resolution of $L$ corresponding to $v$ 
(that is, taking the 0-resolution at the $i$-th crossing
if $v_i =0$, and the 1-resolution otherwise) 
and then placing arcs corresponding to each of the crossings labeled 
by 0’s in $v$ 
(that is, at the $i$-th crossing if $v_i=0$). See Figure \ref{r}.\\

\begin{figure}
\includegraphics[width=140mm]{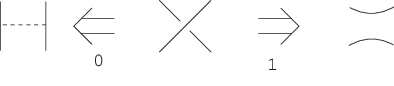}
\caption{{\bf The 0- and 1-resolutions
}\label{r}}   
\end{figure}

Therefore, $n-$ind$(D_L(v)) = |v| =\sum v_i ^2=\sum v_i$, 
the (Manhattan) norm of $v$.
	We use Ind as defined in Definition \ref{2.1}.
\end{definition}

Note that the 0-(respectively, 1-)resolution is 
the same as the $A$-(respectively, $B$-)type split in \cite{Kauffmanstate}.

Compare Definitions \ref{2.1} and \ref{2.2} with \cite[Definitions 2.1 and 2.2]{LSk}. 
In the case of classical link diagrams, 
$Z(D)$ consists in only embedded circles, but 
in the case of virtual link diagrams,  $Z(D)$ consists in immersed circles whose singular points are virtual crossing points.
The dotted arc in Figure \ref{r} corresponds to the 0-smoothing, while the 1- smoothing is undotted. This corresponds to the red  arc in \cite[Figure 2.1.a]{LSk}. 

It is important that we do not carry out a resolution on any virtual crossing.

\begin{figure}
\includegraphics[width=100mm]{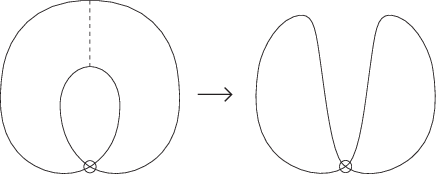}
\caption{{\bf A single cycle surgery
}\label{single}}   
\end{figure}

\begin{definition}\label{2.3} 
{\bf (\cite[Definition 2.3]{LSk}.)}  
Given resolution configurations $D$ and $E$ 
there is a new resolution configuration
$D-E$ defined by
$$Z(D-E)=Z(D)-Z(E)\hskip3mm A(D-E) = \{A\in A(D) | \forall Z\in Z(E): \partial A\cap Z=\phi\}.$$
Let $D\cap E=D-(D-E)$. 
\end{definition}

Note that $Z(D\cap E)=Z(E\cap D)$ and 
$A(D\cap E)=A(E\cap D)$; however, the total orders
on $A(D\cap E)$ and $A(E\cap D)$ could be different.

\begin{definition}\label{2.4} 
{\bf (\cite[Definition 2.4]{LSk}.)}  
The {\it core} $c(D)$ of a resolution configuration $D$ is the resolution configuration
obtained from $D$ by deleting all the circles in $Z(D)$ 
that are disjoint from all the arcs in $A(D)$.  
A resolution configuration $D$ is called {\it basic} if $D = c(D)$, that is, 
if every circle in $Z(D)$ intersects an arc in $A(D)$.
\end{definition}

\begin{definition}\label{scs}
Let $D$ be a resolution configuration.
Suppose that, when we carry out a surgery along 
one arc of $A(D)$ on an immersed circle of $Z(D)$, 
the number of the elements of $Z(D)$ is not changed. 
Then we call this surgery {\it single cycle surgery}. 
See Figure \ref{single} for an example. 
\end{definition}

\vs

Khovanov homology of virtual links is independent of orientations of virtual links 
but we use orientations of virtual links when we define it.

\begin{definition}
{\rm (\cite{Igor})}   
Let $D$ be a (unoriented) virtual link diagram and $\X(D)$ be the set of classical crossings of $D$. For any crossing $c\in\X(D)$ we choose one of the two {\em source-sink orientations}~\cite{IM13} 
(see Fig.~\ref{fig:source_sink_orientations}). 
Denote the choice of source-link orientation at each classical crossing by $\lambda$ 
and call it a {\em local source-sink structure $($LSSS$)$} of the diagram $D$. 
Denote the set of all local source-sink structures by $\Lambda(D)$.

\begin{figure}[h]
\centering\includegraphics[width=0.6\textwidth]{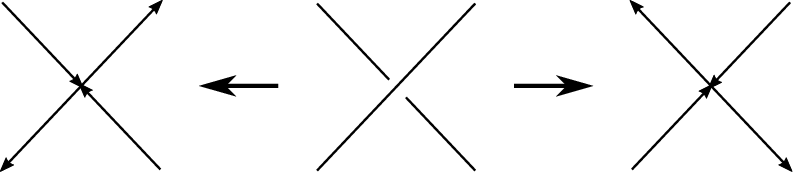}
\caption{Source-sink orientations}\label{fig:source_sink_orientations}
\end{figure}

If the diagram $D$ is oriented, one can define the {\em canonical source-sink structure}  
 \cite{DKK} 
as shown in Fig.~\ref{fig:canonical_source_sink_orientation}.
\begin{figure}[h]
\centering\includegraphics[width=0.4\textwidth]{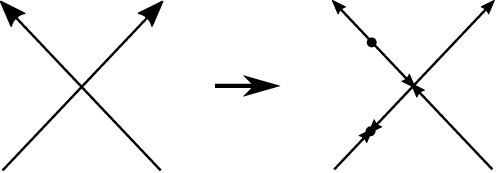}
\caption{Canonical source-sink orientation. The cut loci mark the places where the local source-sink orientation changes to the orientation of the link}\label{fig:canonical_source_sink_orientation}
\end{figure}
\end{definition}

\begin{remark}
{\rm (\cite{Igor})}  
Source-sink structure can be used to define (local) orientation of Kauffman states of the diagram, see Fig.~\ref{fig:source_sink_smoothings}.
\begin{figure}[h]
\centering\includegraphics[width=0.6\textwidth]{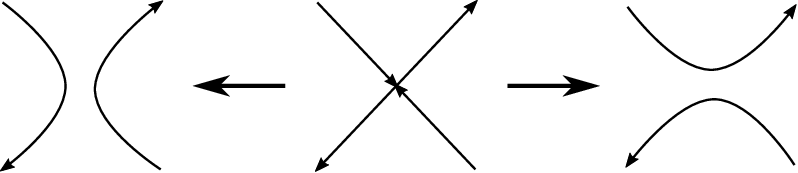}
\caption{Local orientation of smoothed components}
\label{fig:source_sink_smoothings}
\end{figure}
\end{remark}

Any LSSS $\lambda$ has the {\em opposite LSSS} $-\lambda$ that is obtained by the {\em global change of orientation}, i.e. when one switches the source-sink structure at every classical crossing of $D$.

\begin{figure}
\bigbreak
\includegraphics[width=90mm]{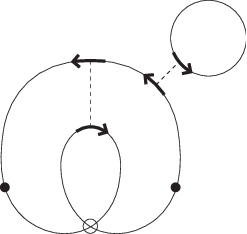}
\caption{{\bf Cut loci}\label{cl}}   
\bigbreak
\end{figure}

\begin{definition}\label{cutlocus}
 {\bf The cut locus}. 
Let $C$ be an immersed circle in $Z(D)$. 
Take all arrows drawn in Figure \ref{fig:source_sink_smoothings}, in $C$. 
See an example in Figure \ref{cl}.  
$C$ is a union of arrows and  immersed curved segments as follows: 
Any immersed curved segment is the image of an immersion $f$ of the closed interval $[0,1]$.  
$f(0)\neq f(1)$.   
Immersed curved segments may intersect other ones at virtual crossings. 
Put a point in $f((0,1))$ 
only if both of $f(0)$ and $f(1)$ touch arrowtails (respectively, arrowheads).     
Don't put the point at any virtual crossing point.   
We call this point the {\it cut locus}. Note that the cut locus point is always at the boundary between opposite local orientations induced by the arrows.
\end{definition}

\noindent{\bf Remark.} 
It is trivial that,
after a surgery along 
one arc of a resolution configuration, 
the placements of cut loci do not change  in $S^2$. 
\\

By Figure \ref{fig:canonical_source_sink_orientation}, we have the following. 

\begin{proposition}\label{soba}
We can put all cut loci near classical crossing points.  
\end{proposition}

By the method to give sink source orientations, 
we have the following. 

\begin{proposition}\label{cc}
There is no cut locus in an arbitrary classical link diagram. 
\end{proposition}

\begin{thm}\label{eve}
The number of cut loci in any immersed loop $C$ in $Z(D)$ is even. 
\end{thm}

\noindent{\bf Proof of Theorem \ref{eve}.}
Let $n$ be a nonnegative integer. 
There are $n$ arrows in $C$ if and only if 
there are $n$ immersed curved segments in $C$. 
If there is a set of two arrows and one immersed curved segment in $C$ 
as drawn in Figure \ref{pico}, 
there is no cut locus in this immersed curved segment. 
 If there is a set  as drawn in Figure \ref{pico}, 
change it into one arrow 
whose orientation of the arrow 
is the same as the two ones in Figure \ref{pico}. 
Repeat this procedure. 
If each immersed curved segment in the result includes  a cut locus, 
then stop the procedure. \\

{\it Claim.} The number of immersed curved segments in 
the result of the process in Figure \ref{pico}.
is even. 
\quad {\it Reason.} 
If the number of arrows in the result is odd, 
there is at least one set of the type in Figure \ref{pico}. 
The immersed curved segment in the set does not include a cut locus by the definition. 
We arrived at a contradiction. 
Therefore the number of immersed curved segments in $C$ is even, 
and the number of cut loci in $C$ is even. 
This completes the proof of Theorem \ref{eve}. 
\qed\\

\begin{figure}
\includegraphics[width=50mm]{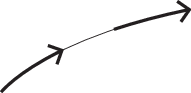}
\caption{{\bf An example that the orientations of two arrowheads are the same. 
The immersed curved segment touches both arrows. 
In this figure, the immersed curved segment is drawn as an embedding, but 
it is not always embedding. 
}\label{pico}}   
\end{figure}

\begin{definition}\label{starting}
Let $D$ be a resolution configuration. 
We notate a star in each immersed circle in $Z(D)$ 
so that it does not  touch an arrow, a cut locus,  nor  a virtual crossing, as indicated below, 
and call it a {\it starting star}. \\

Let $E$ be a resolution configuration that we obtain from $D$ by a surgery along 
one arc of $E$. 
Let $S$ be an immersed circle in $Z(D)\cap Z(E)$. 
Note that the place of $S$ in $S^2$ before this surgery is the  same as that after it. 
We put the starting star in the same place in $S$ of $E$ as  in that of $D$.  
See Figure \ref{star}. 
\end{definition}

\begin{figure}
\vskip-20mm
\includegraphics[width=120mm]{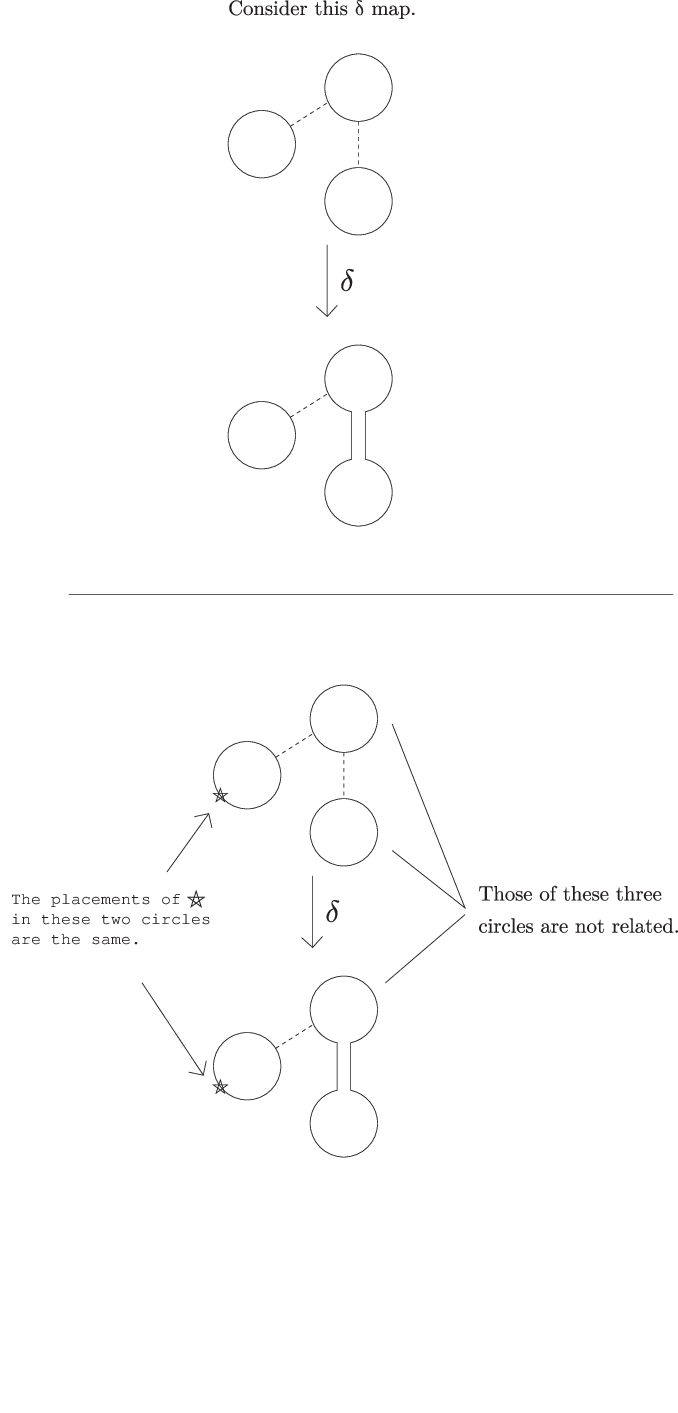}
\vskip-30mm
\caption{{\bf In this figure, circles may not be embedded circles, but 
we draw (not necessarily embedded) immersed circles as embedded ones abstractly. 
}\label{star}}   
\end{figure}

\bigbreak
\noindent{\bf Remark.}  
In 
Figures \ref{star}, 
circles may not be embedded circles, but 
we draw 
them 
as embedded ones abstractly. 
\\

\begin{definition}\label{2.9}  {\bf (\cite[Definition 2.9]{LSk}).} 
A {\it labeled resolution configuration} is a pair $(D, x)$ of a resolution configuration 
$D$ and 
a labeling $x$ 
of each element of $Z(D)$ by either $x_+$ or $x_-$.
\end{definition}

Note that 
 labeled resolution configurations are the same 
as what many people call  
{\it enhanced Kauffman state}s 
or 
{\it enhanced state}s. 
Some people use 
$v_+$ (respectively, $v_-$) 
for 
$x_+$ (respectively, $x_-$).

\begin{definition}\label{2.10}  {\bf (\cite[Definition 2.10]{LSk}).} 
There is a partial order $\prec$ on labeled resolution configurations defined as
follows. 
We declare that $(E, y)\prec(D, x)$ if:

\begin{enumerate}
\item[(1)]
 The labelings $x$ and $y$ induce the same labeling on 
$D\cap E = E\cap D.$ 

\item[(2)] 
$D$ is obtained from $E$ 
by surgering along a single arc of $A(E)$. In particular, either:

\begin{enumerate}
\item[(a)]
$Z(E-D)$ contains exactly one circle, say $Z_i$, and 
$Z(D-E)$ contains exactly two
circles, say $Z_j$ and $Z_k$, or

\item[(b)]
$Z(E-D)$ contains exactly two circles, say $Z_i$ and $Z_j$, and 
$Z(D-E)$ contains exactly one circle, say $Z_k$.

\end{enumerate}

\item[(3)]
 In Case (2a), either $y(Z_i) = x(Z_j) = x(Z_k) = x_-$ or 
$y(Z_i) = x_+$ and $\{x(Z_j), x(Z_k)\} =\{x_+, x_-\}$.  

In Case (2b), either $y(Z_i) = y(Z_j) = x(Z_k) = x_+$ 
or $\{y(Z_i), y(Z_j)\} = \{x_-, x_+\}$ and $x(Z_k) = x_-$.
\end{enumerate}
\end{definition}

Note that no surgery in Definition \ref{2.10} is a single cycle surgery although we now consider virtual links. 

See \cite[\S2]{LSk} for $s(D)$.

\begin{definition}\label{2.11}   {\bf (\cite[Definition 2.11]{LSk}.)}
A {\it decorated resolution configuration} is a triple $(D, x, y)$ 
where $D$ is a resolution configuration and 
$x$ (respectively,  $y$) is a labeling of each component of $Z(s(D))$ 
(respectively, $Z(D)$) by an element of ${x_+, x_-}$. 
If we do not have $(D, y) \underline{\prec} (s(D), x)$, 
we say that $(D, x, y)$ is empty.

Associated to a decorated resolution configuration $(D, x, y)$ is the {\it poset}  
 $P(D, x, y)$ consisting
of all labeled resolution configurations 
$(E, z)$ with 
$(D, y)\underline{\prec}(E, z)\underline{\prec}(s(D), x)$.
We say that  $P(D, x, y)$ is the poset of a decorated resolution configuration $(D, x, y)$.   
\end{definition}

\h{\bf Remark.} In \cite[Definition 2.11]{LSk}, 
if we write $(D, x, y)$, then $(D, x, y)$ is not empty. 
We define a term, ``empty decorated resolution configuration''  for convenience.
If $(D, x, y)$ is empty, $P(D, x, y)$ is the empty set.

\begin{definition}\label{2.15gr}   {\bf(A part of Definition 2.15 in \cite{LSk}.)} 
For labeled resolution configurations, 
{\it homological grading $\text{gr}_h$} and 
{\it a quantum grading $\text{gr}_q$}, defined as follows:

\vskip2mm
\hskip7mm
$\text{gr}_h((DL(u), x)) = -n_- + |u|,$ 

\vskip2mm
\hskip7mm
$\text{gr}_q((DL(u), x)) = n_+ - 2n_- + |u|$

\hskip40mm
$+ 
\sharp\{Z\in Z(DL(u)) | x(Z) = x_+\}
-\sharp\{Z\in Z(DL(u)) | x(Z) = x_-\}.$ 
\vskip2mm

Here $n_+$ denotes the number of positive crossings in $L$; and 
$n_- = n-n_+$ denotes the number of negative crossings.
\end{definition}

\begin{remark}\label{muki}
If $L$ has only one component, 
when we define $n$, $n_+$, $n_-$, $\text{gr}_h$, and $\text{gr}_q$, 
we do not need to use the orientation of $L$. 
If $L$ has greater than one component, 
we use that of $L$ when we define $n_+$ and $n_-$. 
Note that, if we change the orientation of $L$ into 
the opposite  one, 
then 
neither $n_+$ nor $n_-$ changes.

Each of $\text{gr}_h(D,x)$ and $\text{gr}_q(D,x)$ 
is independent of which orientation we give $Z(D)$, 
and is independent of which we choose $L$ or $-L$. 
Here, let $-L$ be a link made from $L$ by reversing the orientation of $L$. 
\end{remark}

%

\begin{fact}\label{qde}
Assume that a single cycle surgery changes  
a $($non-labeled$)$ resolution configuration $D_L(u)$ into $D_L(v)$. 
Let $A_i$ $($respectively, $A_j$$)$ be a labeled resolution configuration defined on  
 $D_L(u)$ $($respectively, $D_L(v)$$)$.  
 Then $A_i$ and $A_j$ have different quantum gradings. 
\end{fact}

\h{\bf Proof of Fact \ref{qde}.}
Recall the definition of quantum gradings, $\text{gr}_q((DL(u), x))$, above. 
Since $n_+$ and $n_-$ are deteremined by a given virtual link diagram, 
a single cycle surgery does not change $n_+ - 2n_-$. 
By the definition of a (single cycle) surgery and that of $|\hskip2mm|$, 
 a single cycle surgery changes the parity of $|u|$. 
Since a single cycle surgery does not change the number of immersed circles in 
a labelled resolution configuration, 
 a single cycle surgery does not change the parity of 
$ \sharp\{Z\in Z(DL(u)) | x(Z) = x_+\}
-\sharp\{Z\in Z(DL(u)) | x(Z) = x_-\}.$ 
Therefore a single cycle surgery always changes 
a quantum grading $\text{gr}_q((DL(u), x))$. 
\qed\bb

We define 
the integral ($\Z$-coefficient) Khovanov chain complex for $L$ 
in Definition \ref{korekoso} 
after we state the conditions with which it should satisfy.\\

We define the Khovanov chain complex to be generated by 
all labeled resolution configurations made from a virtual link diagram $L$.  
Let $\{A_i\}_{i\in\Lambda}$ be the set of all labeled resolution configurations made from $L$. 
Note that $\Lambda$ is a finite set. 
We will define $$\delta A_i = \displaystyle\sum_{j\in\Lambda} c[A_i;A_j]\cdot A_j.$$
Here, $c[A_i;A_j]$ is an integer coefficient.  
We only have to define  $c[A_i;A_j]$,  
which should have the following properties. 

\bigbreak\noindent(1)
If $A_i$ and $A_j$ have different quantum gradings, then $c[A_i;A_j]=0$. 
Note: This condition holds in the case of Khovanov homology for classical links.
If $A_i$ and $A_j$ are given as in 
Fact \ref{qde} above, then 
we want to define $c[A_i;A_j]=0$. 
We will explain why we want this condition in Remark \ref{point}.

\bigbreak\noindent(2) 
If (the homological grading $A_i)+1\neq$ (that of $A_j$),  then $c[A_i;A_j]=0$. 
\bigbreak

\bigbreak\noindent(3) 
Suppose that  (the homological grading $A_i)+1=$ (that of $A_j$), \\
 and that  (the quantum grading $A_i)=$ (that of $A_j$). 
Then $c[A_i;A_j]$ may not be zero. 
\vskip9mm


\vs

We change notations, and continue to explain how we define virtual Khovanov homology. 
Let $\alpha$ and $\beta$ be enhanced Kaffman states.
Let $\alpha^\text{no}$ (respectively, $\beta^\text{no}$) 
be a (non-enhanced) Kaffman state under 
$\alpha$ (respectively, $\beta$). 

We define 
the coefficient $[\alpha:\beta]$ of 
$\beta$ in 
$\delta\alpha$ to be nonzero only if 
$\alpha$ and $\beta$ satisfy the following condition $(\ast)$. 
\\

\h$(\ast)$
$\beta$ is obtained from $\alpha$ 
by a single multiplication or a single co-multiplication 
drawn in Figure \ref{hitsuyo}.  
\\

\h Furthermore we define  $[\alpha:\beta]$ in the case of $(\ast)$ to be $+1$ or $-1$.  
We explain how we define it.  
\\

We define $[\alpha:\beta]$ to be a product,  
$\mathcal O(\alpha^\text{no}, \beta^\text{no})\mathcal P(\alpha, \beta)$, 
of 
$\mathcal O(\alpha^\text{no}, \beta^\text{no})$ and $\mathcal P(\alpha, \beta)$. 
 
$\mathcal P(\alpha, \beta)$ depends on labelings  in general.
  $\mathcal P(\alpha, \beta)$ is defined by using cut loci and starting stars as below.

$\mathcal O(\alpha^\text{no}, \beta^\text{no})$ is independent of labelings. 
 $\mathcal O(\alpha^\text{no}, \beta^\text{no})$ is defined 
 by using local and global orders of loops in Kauffman states below. \\

\begin{definition}\label{gloc}
{\bf Global and local orders of loops in Kauffman states, 
and the sign $\mathcal O(\quad, \quad).$}  

\h(1)
Let $D$ be a diagram of a virtual link and $\X(D)$ be the set of the classical crossings of $D$. Denote the set of Kauffman states of the diagram $D$ by $S(D)=\{s\colon\X(D)\to\{0,1\}\}$.
Let  $\Gamma(s)$ be the sets of components, circles,  or loops, 
of Kauffman states. 
Fix an arbitrary order $\sigma_s\colon \{1,2,\dots, |\Gamma(s)|\}\to \Gamma(s)$ of the components for all Kauffman states $D_s$, $s\in S(D)$.
This order is called a {\it global order} of  $D_s$. 

\h (2)
At each crossing we choose one of the incoming edges in the source-sink orientation as shown in Fig.~\ref{fig:source_sink_orders}. We mark the distinguished edge with a box. The choice of the edges for all classical crossings will be denoted by $o=\{o_c\}_{c\in\X(D)}$ and called an {\em oriented direction system} on the LSSS $\lambda$.

\begin{figure}[h]
\centering\includegraphics[width=0.4\textwidth]{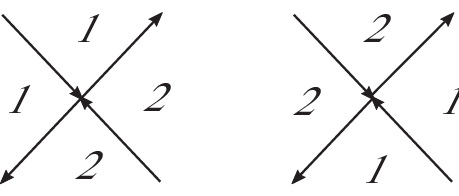}
\caption{Directions of a source-sink orientation. It is associated with Figure 
\ref{fig:source_sink_orientations}. 
}\label{fig:source_sink_orders}
\end{figure}

If $\lambda$ is a canonical LSSS of the oriented diagram $D$, we can consider the {\em canonical oriented direction system} choosing the incoming edges as shown in Fig.~\ref{fig:canonical_order}.

\begin{figure}[h]
\centering\includegraphics[width=0.4\textwidth]{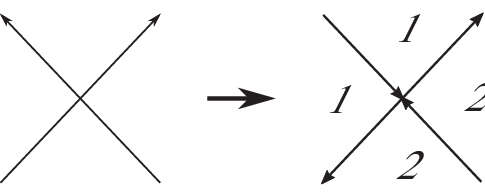}
\caption{Canonical oriented direction system. 
It is related to Figure \ref{fig:canonical_source_sink_orientation}.}
\label{fig:canonical_order}
\end{figure}

The {\it oriented direction} or {local order} $o_c^s$ in a crossing $c$ defines a local ordering of the components in Kauffman states which pass by $c$: Here, $s$ defines a enhanced Kauffman state.  
See Fig.~\ref{fig:source_sink_smoothing_order}. We can consider the local order as a family of maps $o^s_c\colon\{1,2\}\to\Gamma(s)$, $c\in\X(D)$, $x\in S(D)$ (the images of $1$ and $2$ may coincide).

\begin{figure}[h]
\centering\includegraphics[width=0.6\textwidth]{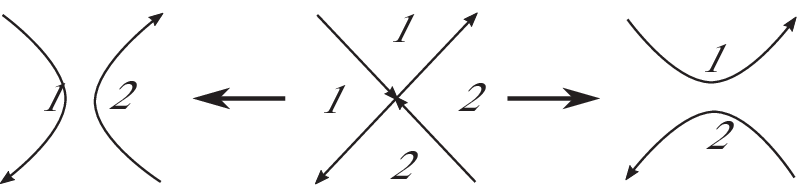}
\caption{Local ordering of component induced by an oriented direction}\label{fig:source_sink_smoothing_order}
\end{figure}

For any crossing $c\in\X(D)$, state $s\in S(D)$, global order $\sigma_s$ and oriented direction system $o_c^s$ we define two new orders:
\begin{itemize}
\item an order $\sigma_s\setminus c$ on the set $\Gamma(s)\setminus c=\Gamma(s)\setminus\im o^s_c$ consisting of the components in $D_s$ which don't pass by the crossing $c$. The order $\sigma_s\setminus c$ is the restriction of $\sigma_s$ to this set. Formally, let
    $j_1=\min\sigma_s^{-1}(\im o^s_c)$ and $j_2=\max\sigma_s^{-1}(\im o^s_c)$ be the indices of the component which pass by $c$ in the order $\sigma_s$. Then $j_1\le j_2$, they may coincide. We define
\begin{equation}\label{eq:sigma-c}
(\sigma_s\setminus c)(i)=\left\{\begin{array}{cl}
\sigma_s(i), & i<j_1,\\
\sigma_s(i+1), & j_1\le i<j_2-1,\\
\sigma_s(i+2), & i>=j_2-1.
\end{array}\right.
\end{equation}

\item an order $\sigma_s\triangleleft(o,c)$ on the set $\Gamma(s)$. Informally, we start the numbering of the components of $D_s$ with those that pass the crossing $c$, and enumerate them in the order $o_c^s$. The other components are ordered according $\sigma_s$. Let us give the explicit formulas. Let $n_s(c)$ be the number of component of $D_s$ which pass by the crossing $c$. Then
\begin{equation}\label{eq:sigma*oc}
(\sigma_s\triangleleft(o,c))(i)=\left\{\begin{array}{cl}
o_c^s(i), & i\le n_s(c),\\
(\sigma_s\setminus c)(i-n_s(c)), & i>n_s(c).
\end{array}\right.
\end{equation}
\end{itemize}

Given two orders $\sigma_1,\sigma_2\colon\{1,\dots,|Z|\}\to Z$ on a finite set $Z$, we define the sign $\epsilon(\sigma_1,\sigma_2)$ as the sign of the permutation $\sigma_1^{-1}\circ\sigma_2$ on the set $\{1,\dots,|Z|\}$.

Now we can define $\sign(s,s')$ of the differential $\partial_{s\to s'}$ as follows

\begin{equation}\label{eq:differential_sign}
\sign(s,s')=\epsilon(\sigma_s,\sigma_s\triangleleft(o,c))\cdot\epsilon(\sigma_s\setminus c,\sigma_{s'}\setminus c)\cdot\epsilon(\sigma_{s'},\sigma_{s'}\triangleleft(o,c)).
\end{equation}

If $s$ (respectively, $s'$) defines an enhanced Kauffman state 
$\alpha$ (respectively, $\beta$), 
we write  $\sign(s,s')$ by $\mathcal O(\alpha, \beta)$. 
\end{definition}

\begin{definition}\label{korekoso} 
(1) {\bf The differential $\delta$.}
Given an oriented virtual link diagram $L$ with $n$ crossings, 
an ordering of the crossings in $L$ and global and local orders of
loops, or circles, in each Kauffman state, 
the {\it Khovanov chain} complex is defined as follows.
The Khovanov chain group $KC(L)$ is the $\Z$-module freely generated 
by labeled resolution configurations of the form $(DL(u), x)$ for $u\in\{0, 1\}^n$.
The differential preserves the quantum grading, increases the homological grading by 1,
and is defined as follows. 
\\

$\displaystyle\delta(D_L(v),y)=
\displaystyle\sum_\text{all $(D_L(u),x)$ as below},  
(\mathcal O((D_L(v)), (D_L(u)))(-1)^{\zeta}. 
$

\h Here, $(D_L(u),x)$ satisfies the condition  
$|u|=|v|+1, (D_L(v),y)\prec (D_L(u),x)$
The number 
$\zeta=\zeta((D_L(u),x), (D_L(v),y))$ 
is defined for a pair $(D_L(u),x)$ and $(D_L(v),y)$
in  Definition \ref{korekoso}.(3) 
by using the number $\xi$ defined in Definition \ref{korekoso}.(2). 

Here is an out line of how to define $\xi$ and $\zeta$.  
We count the parity of cut loci traversed  between a surgered point and a star, on a loop. 
If the loop is labeled by $x_-$ and if the parity on it is odd, 
we multiply $-1$ (once) in the coefficient. 
If the loop is labelled by $x_+$, we do not multiply $-1$. 

The number $\mathcal O((D_L(v)), (D_L(u)))$ is defined above.

The numbers $\mathcal O((D_L(v)), (D_L(u)))$ and 
$(-1)^{\zeta((D_L(u),x), (D_L(v),y))}$ are $+1$ or $-1$. 
The product 
$\mathcal O((D_L(v)), (D_L(u)))
(-1)^{\zeta((D_L(u),x), (D_L(v),y))}$ 
is a coefficient.

Recall 
the above notation 
$\mathcal P(\quad,\quad)$: 
We define 
$\mathcal P((D_L(u),x), (D_L(v),y))$ 
to be 
$(-1)^{\zeta((D_L(u),x), (D_L(v),y))}$.  
\end{definition}

\begin{remark}\label{kazoe}
The reader should note the following facts, before reading the definitions below. 
\\

Let 
$\alpha$ and $\beta$ 
be  two enhanced Kauffman states of a classical link diagram. 
Let 
$\alpha^\text{no}$ (respectively, $\beta^\text{no}$)  
be a (non-enhanced) Kauffman state under 
$\alpha$ (respectively, $\beta$)  
 Let 
$\mathcal A(\alpha^\text{no}, \beta^\text{no})$ be 
an integer defined for  
$\alpha^\text{no}$ and $\beta^\text{no}$, which is 
the same one as $(-1)^{s_0(\mathcal C_{u,v})}$ 
in \cite[Definition 2.15]{LSk}: 
Each (non-enhanced) Kauffman state is characterized 
by each vector $w=(w_1,...w_{n})$, where $w_*\in\{0,1\}$ (See Definition {2.2}.). 
Let $\alpha$ (respectively, $\beta$)  be characterized by 
$u$ (respectively, $v$). 
We define $s_0(C_{u,v})\in\Z_2$  as follows: 
if $u = (\epsilon_1,..., \epsilon_{i-1}, 1, \epsilon_{i+1}, . . . , \epsilon_n)$ and 
$v =(\epsilon_1,..., \epsilon_{i-1}, 0, 
\epsilon_{i+1}, . . . , \epsilon_n)$, 
then $s_0(C_{u,v}) = (\epsilon_1+\cdot\cdot\cdot+ \epsilon_{i-1})$; 
see also \cite[Definition 2.15]{LSk}.\\

In Khovanov's original case, the case of classical links in $S^3$, 
we define the coefficient of $\beta$ in $\delta\alpha$
to be $\mathcal A(\alpha^\text{no}, \beta^\text{no})$.  

The issue is as follows. 
See Figure \ref{eK1}. 
See also \cite[Figure 13 and the explanation about it]{DKK}. 
Call the left upper Kauffman state with a labeling $x_+$, $A$ as in Figure \ref{eK2}.
If we 
define the coefficient of $\beta$ in $\delta\alpha$
by using $\mathcal A(\alpha^\text{no}, \beta^\text{no})$ as in the case of classical links, 
we have $\delta^2A\neq0$.
So we introduce 
$\mathcal O(\alpha^\text{no}, \beta^\text{no})$ and $\mathcal P(\alpha, \beta)$, 
and settle this issue. 

See Theorem \ref{ittoku} on a relation among 
 $\mathcal A(\alpha^\text{no}, \beta^\text{no})$,  
$\mathcal O(\alpha^\text{no}, \beta^\text{no})$ and $\mathcal P(\alpha, \beta)$.
\end{remark}

\begin{figure}
\bigbreak
\includegraphics[width=130mm]{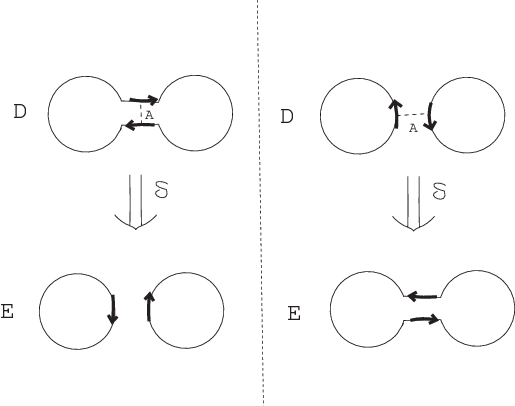}
\caption{{\bf 
A multiplication surgery and a co-multiplication surgery
}\label{atokara}}   
\end{figure}

\vskip3mm
\noindent{\bf Definition \ref{korekoso}.(2)  
The number $\xi$.}  
We define the coefficient in 
Definition \ref{korekoso}.(1) to be nonzero 
only if      $(D_L(u),x)$ is obtained from $(D_L(v),y)$ by 
one multiplication surgery or one co-multiplication surgery. 
We define a number $\xi$ 
before we introduce a number $\zeta$.  
These numbers are obtained by appropriate counts of the parity of cut points on a path or paths between the site of the algebraic operation (multiplication or co-multiplication) and the starting stars on the loops of the state configuration.

See Figure \ref{atokara}. 
Kauffman states $D$ and $E$ are the same in the part other than Figure \ref{atokara}: 
We carry out a surgery along an arc $A$. 
  $(D_L(u),x)$ and $(D_L(v),y)$ are different 
 only in the part like 
Figure \ref{atokara}.  
 
In Figure \ref{atokara}, 
we put arrows at two foots of an arc where we carry out surgery 
according to Figures 
\ref{fig:canonical_source_sink_orientation}
and 
\ref{fig:source_sink_smoothings}. 
We also put arrows the points in Kauffman states after a surgery 
according to Figures 
\ref{fig:canonical_source_sink_orientation}
and 
\ref{fig:source_sink_smoothings}.

Since we do not consider a single cycle surgery here, 
there are three loops related to this surgery.
We define a number $\xi$ 
for each loop in $D$ (respectively, $E$) and the arc $A$ 
in Figure \ref{atokara}. 
There are two cases.
\\

We count cut loci in each the waved curve of  Figure \ref{atokara2}. 
More precisely, we do as follows.

\begin{figure}
\bigbreak
\includegraphics[width=130mm]{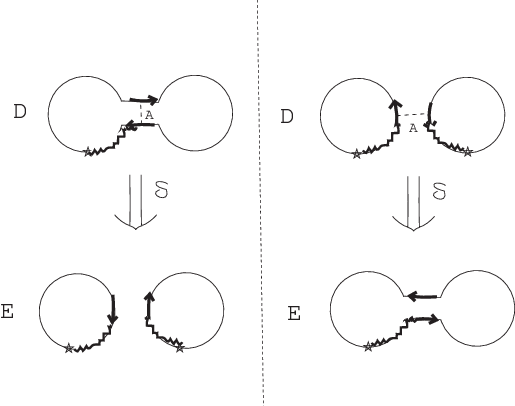}
\caption{{\bf 
A multiplication surgery and a co-multiplication surgery
}\label{atokara2}}   
\end{figure}

\noindent 
{\bf Case 1.} 
The case where the loop includes only one arrow: 
Each of two loops in the left lower $E$ of Figure \ref{atokara2}, 
and 
Each of two loops in the right upper  $D$   of Figure \ref{atokara2}. 

See an example in Figure \ref{shan}.  
See $X$ and $Y$ in Figure \ref{shan}. 

%
\begin{figure}
\includegraphics[width=60mm]{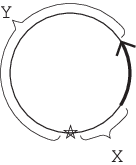}
\caption{{\bf A figure associated with the definition of $\xi$ 
in Definition \ref{korekoso}.(2).   
}\label{shan}}   
\end{figure}

Let $\xi$ be 0 (respectively, 1) 
if the number of the cut loci in $X$ is even (respectively, odd). 
Note the following: 
By Theorem \ref{eve}, the parity of the number of cut loci in 
both immersed curved segments, $X$ and $Y$,  
are the same. 
If we change $X$ into $Y$ in the definition of $\xi$, the value of $\xi$ is the same. \\
%
\\
%
%
%

\noindent 
{\bf Case 2.} 
The case where the loop includes the two arrows: 
A unique loop in the left upper $D$ of Figure \ref{atokara2}, 
and a unique loop in the right lower  $E$   of Figure \ref{atokara2}.

See four cases in Figure \ref{shan2}. 

\begin{figure}
\includegraphics[width=100mm]{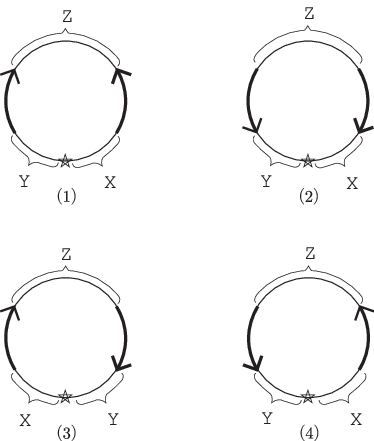}
\caption{{\bf Figures associated with the definition of $\xi$ 
in Definition \ref{korekoso}.(2). 
}\label{shan2}}   
\end{figure}

The case of  Figure \ref{shan2}.(1) and that of Figure \ref{shan2}.(2) 
occur if and only if 
the surgery is a single cycle surgery. 
We do not need these cases now. 
In the case of Figure \ref{shan2}.(3) and that of Figure \ref{shan2}.(4),  
let $\xi$ be 0 (respectively, 1) if the number of the cut loci in $X$ is even (respectively, odd).
Note the following:
Each of $Z$ in Figure \ref{shan2}.(3) and $Z$ in Figure \ref{shan2}.(4) is an immersed circle 
before (respectively, after) this surgery. 
Hence the number of cut loci in $Z$ is even by Theorem \ref{eve}.  
Furthermore, 
by Theorem \ref{eve},  
the sum of the number of cut loci in $X$, that in $Y$, and that in $Z$ is even. 
Therefore 
the parity of the number of cut loci in $X$ and that in $Y$ 
are the same. 
If we change $X$ into $Y$ in the definition of $\xi$, the value of $\xi$ is the same. \\

Note that, when we consider classical links, then $\xi=0$ in all cases. \\

\begin{remark}\label{point}
In the case of  Figure \ref{shan2}.(1) and (2), 
$Z$ does not make an immersed circle before (respectively, after) this surgery, 
because it is a single cycle surgery. 
See Figure \ref{pipi} for an example. 
Hence we cannot claim that the parity of the number of   cut loci  in $Z$ is zero, 
unlike the case of  Figure \ref{shan2}.(3) and (4). 
Therefore we can not conclude that  
the parity of the number of cut loci  in $X$ is equivalent to  that in $Y$, 
unlike the case of  Figure \ref{shan2}.(3) and (4).  
We cannot determine which we choose that of $X$ or that of $Y$.
\begin{figure}
\includegraphics[width=80mm]{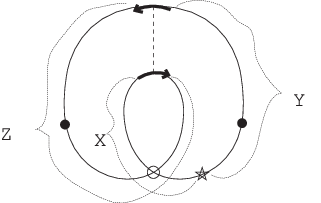}
\caption{{\bf 
A resolution configuration with cut loci and a starting star
}\label{pipi}}   
\end{figure}
One way of settling this situation 
is defining  
  the coefficient in the differential associated with single cycle surgeries to be zero. 
More precisely:  
As written in (1) right above Definition \ref{korekoso}.(1),  
 if $A_i$ and $A_j$ are given as 
 in Fact \ref{qde}, 
 then we define $c[A_i;A_j]=0$. 
Recall Remark \ref{kazoe}. 
We can summarize this by saying that 
the coefficient 
that corresponds to a single cycle surgery is taken to be zero. 
\end{remark}

\noindent{\bf Definition \ref{korekoso}.(3)  
The definition of $\zeta((D_L(u),x), (D_L(v),y))$.}        
In order to define $\zeta((D_L(u),x), (D_L(v),y))$, 
we need to check only four cases of surgeries along each arc, 
which are drawn in Figure \ref{hitsuyo}. 
Here, 
$\Delta_*, \nabla_*, M_*$, and $W_*$ denote a labeled resolution configuration, 
and $\alpha_*, \alpha_{\sharp\natural},  \beta_*, \gamma_*, \sigma_*$ an immersed circle. 
We carry out 
one surgery in each of Figures \ref{hitsuyo}.(1)-(4).
We define the number 
$\xi(\alpha_*)$ (respectively,  
$\xi(\alpha_{\sharp\natural}),  
\xi(\beta_*), 
\xi(\gamma_*),$ and 
$\xi(\sigma_*)$)  
 for this surgery and an immersed circle
$\alpha_*$ (respectively, $\alpha_{\sharp\natural},  \beta_*, \gamma_*, \sigma_*$). 
We define 
\\

\begin{figure}
\includegraphics[width=120mm]{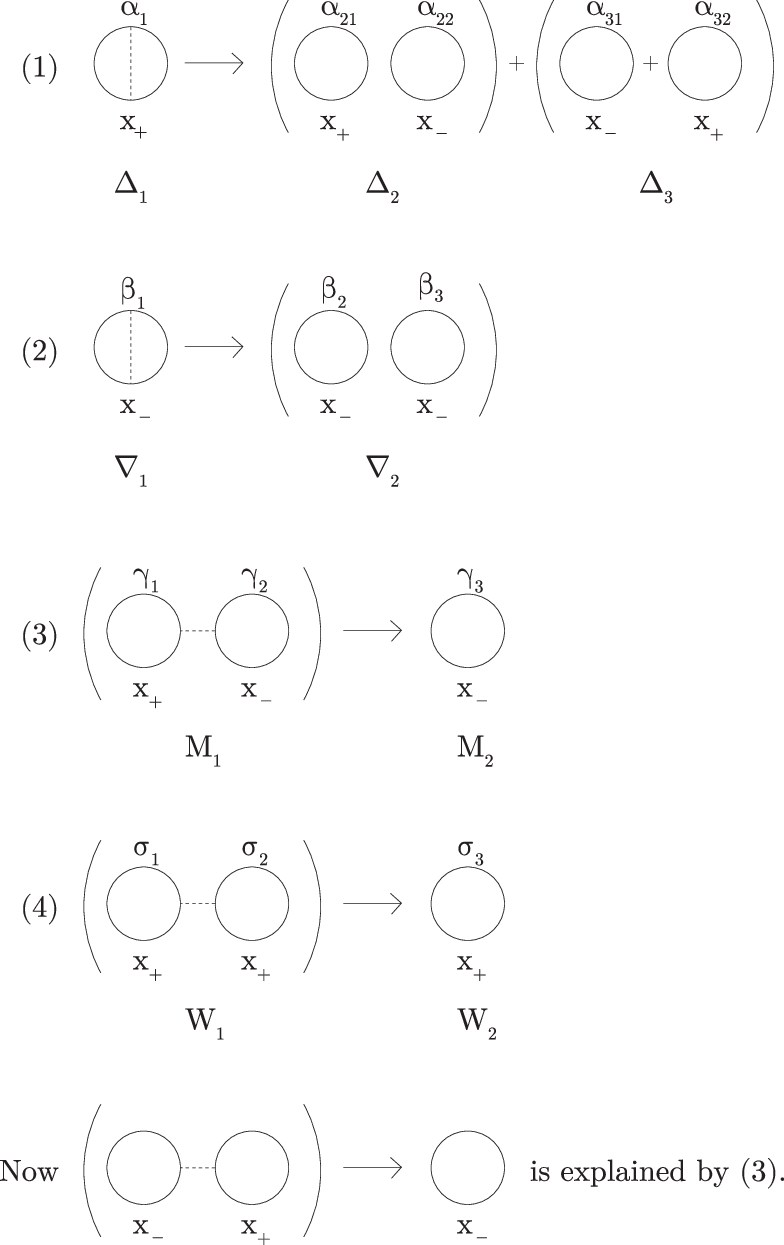}
\caption{{\bf We define $\xi$ by using these figures.
}\label{hitsuyo}}   
\end{figure}

$\zeta(\Delta_1,\Delta_2)=\xi(\alpha_{22})$

$\zeta(\Delta_1,\Delta_3)=\xi(\alpha_{31})$

$\zeta(\nabla_1,\nabla_2)=\xi(\beta_1)+\xi(\beta_2)+\xi(\beta_3)$

$\zeta(M_1, M_2)=\xi(\gamma_2)+\xi(\gamma_3)$

$\zeta(W_1,W_2)=0$. 
\\

\h{\bf Remark.} 
In the case of Figures \ref{hitsuyo}.(1), 
we have two cases   
$\zeta(\Delta_1,\Delta_2)$
and 
$\zeta(\Delta_1,\Delta_3)$. 
The coefficient is defined for each pair of enhanced Kauffman states. 
\\


By an explicit calculation, we have the following.

\begin{thm}\label{wellde}{\bf (\cite{Man})}
For $\delta$ in Definition \ref{korekoso},   
we have $\delta\cdot\delta=0$. 
That is, 
the $\Z$ coefficient Khovanov homology for virtual links is defined. 
\end{thm}

Thus 
Definition \ref{korekoso} 
is well-defined. 

This homology of $L$ is the same as that of $-L$. 

\begin{definition}\label{Khteigi}
Let $L$ be a virtual link diagram. 
Take the Khovanov chain complex for $L$.  
The Khovanov chain complex with the differential yields the 
{\it Khovanov homology for the virtual link diagram $L$}.  
Let $A_i$ be each resolution configuration made from $L$
($i\in\Lambda$.).  
Note that $\Lambda$ is a finite set. 
Then $\{A_i\}_{i\in\Lambda}$ is the basis 
of Khovanov chain complex for $L$.  
We call each $A_i$  {\it Khovanov basis element}  
($i\in\Lambda$), and $\{A_i\}_{i\in\Lambda}$ {\it Khovanov basis}.
\\

Let $Kh^{q,i}(L)$ (respectively, $C^{q,i}(L)$)
denote Khovanov homology (respectively, Khovanov chain complex)
with the $\Z$ coefficient 
of quantum grading $q$ and homological grading $i$ 
for a virtual link diagram $L$. 
We sometimes omit the words, $\Z$ coefficient, when it is clear from the context.\\

By using  the Khovanov chain complex  $C^{q,i}(L;\Z)$, 
we can define 
the Khovanov chain complex $C^{q,i}(L;\Z_2)$, 
and  
the Khovanov homology $Kh^{q,i}(L;\Z_2)$ 
with  $\Z_2$ coefficients, quantum grading $q$ and homological grading $i$ 
for a virtual link diagram $L$. 
\end{definition}


See \cite{DKK, Man, Igor}. 
For a fixed virtual link diagram and its associated Khovanov chain complex, 
Khovanov homology is independent of the placement of the starting star 
in each immersed circle of each labeled resolution configuration. 
Khovanov homology for a given virtual link diagram $L$ 
does not change by Reidemeister moves on $L$.  
Thus the following is well-defined. 

\begin{definition}\label{cohomology} 
Let $\mathcal L$ be a virtual link. 
Let $L$ be a virtual link diagram which represents $\mathcal L$. 
Define Khovanov homology $Kh^{q,i}(\mathcal L)$  to be $Kh^{q,i}(L)$ for $L$. 
We can define $Kh^{q,i}(\mathcal L;\Z_2)$, as well.  
\end{definition}

The Jones polynomial of any virtual link $L$ is a graded Euler characteristic of 
the $\Z$-coefficient Khovanov homology of $L$.
See \cite{DKK, Man}.

\begin{definition}\label{homology}
Take $L$ and $\mathcal L$ above.  
Use Hom$(C^{q,i}(L;\Z), \Z)$, and  
$<\partial a, \alpha>=<a, \delta\alpha>$ for a dual chain $a$ and 
a Khovanov chain $\alpha$, 
where $< , >$ is the Kronecker product as usual. 
Thus we can define the {\it dual Khovanov chain complex}  $C_{q,i}(L)$, and  
the {\it homology of the dual Khovanov chain complex},  $Kh_{q,i}(L)$.  
We can define  $Kh_{q,i}(\mathcal L;\Z_2)$, as well.  
It does not matter if we write 
$K_{q,i}(L)$ (respectively, $C_{q,i}(L)$) 
as 
$K^q_i(L)$ (respectively, $C^q_i(L)$).  \\

Let $A_i$  be each Khovanov basis element 
for a virtual link diagram $L$ 
($i\in\Lambda$, where $\Lambda$ is a finite set).  
 Let $a_i$ ($i\in\Lambda$) be 
each basis element 
of the dual Khovanov chain complex for $L$ 
such that 
$<A_l,a_k>=\delta_{l,k}$ for 
two arbitrary elements, $l$ and $k$, in $\Lambda$. 
 We call each $a_i$ the {\it dual Khovanov basis element} 
of the dual Khovanov chain complex for $L$, 
and $\{a_i\}$  {\it dual Khovanov basis}.  
We call the dual Khovanov basis element 
the Khovanov basis element 
when it is clear from the context what is meant.\\  

Define the homological grading gr$_ha_i$ to be gr$_hA_i$ for any $i\in\Lambda$. \\

Let 
$\{a_*\}_{*\in\Lambda}$ be the dual Khovanov basis. 
Define a partial order $\prec$
on the set  $\{a_*\}_{*\in\Lambda}$  
 as follows: 
Let $k,l\in\Lambda$. 
$a_k\prec a_l$ if and only if $A_l\prec A_k.$ \\

Define $c[a_l;a_k]$ to be the coefficient in 
$\partial a_k=\displaystyle\sum_{l\in\Lambda}c[a_l;a_k]a_l$.
Note  that $c[a_l;a_k]=c[A_k;A_l]$. 
\end{definition}


\begin{thm}\label{ittoku} 
Let $\mathcal L$ be a classical link. 
Note that  $\mathcal L$ is a virtual link since any classical link is also a virtual link.
Then the virtual Khovanov homology for the link $\mathcal L$ is 
the original  Khovanov homology for the link $\mathcal L$ as a classical link. 
\end{thm}

\noindent{\bf Proof of Theorem \ref{ittoku}.}
Recall $\mathcal A(\alpha^\text{no}, \beta^\text{no})$  in Remark \ref{kazoe}.

Let $K$ be a classical link diagram. 
Note that there is no cut loci in any Kauffman state made from $K$ 
if we use the rule in 
Figure \ref{fig:canonical_source_sink_orientation}.  
Hence we do not need $\mathcal P(\quad,\quad)$ when we define the coefficient.
We have the following.

\begin{cla}\label{seyade}
Let $K$ be a classical link diagram. 
The Khovanov chain complex of $K$ defined 
by the coefficients 
$\mathcal A(\alpha^\text{no}, \beta^\text{no})$
is chain isomorphic to that by 
$\mathcal O(\alpha^\text{no}, \beta^\text{no})$. 
\end{cla}

\h{\bf Proof of Claim \ref{seyade}.}   
If the number of classical crossings of a given virtual link diagram is zero, 
Claim \ref{seyade} is true. 
\\


Let $K$ be a virtual link diagram with $m$  classical crossings.  

Make all Kauffman states and 
all enhanced Kauffman states. 

Let $P$ be the poset of all (non-enhanced) Kauffman states 
made 
by the natural partial order on the set of the vectors. 

Let $x,y\in P$
satisfy the following condition:  
Only one component is different, comparing the components of of $x$ with those of $y$. 
Then we give $+1$ or $-1$ 
to the pair $x,y\in P$ 
in two ways defined 
by $\mathcal O(\quad,\quad)$ and  $\mathcal A(\quad,\quad)$.

Let $i,j\in\{1,...,n+1\}$. Let \\
$v_{00}=(...,v_{i-1},0,v_{i+1},...,v_{j-1},0,v_{j+1},...)$\\
$v_{01}=(...,v_{i-1},0,v_{i+1},...,v_{j-1},1,v_{j+1},...)$\\
$v_{10}=(...,v_{i-1},1,v_{i+1},...,v_{j-1},0,v_{j+1},...)$\\
$v_{11}=(...,v_{i-1},1,v_{i+1},...,v_{j-1},1,v_{j+1},...)$.
\\

By the definition of 
$\mathcal A(\quad,\quad)$, we have 
\begin{equation}\label{eq:A}
\mathcal A(v_{00},v_{01})\mathcal A(v_{01},v_{11})=
-\mathcal A(v_{00},v_{10})\mathcal A(v_{10},v_{11}).
\end{equation}

By the definition of 
$\mathcal O(\quad,\quad)$, we have 
\begin{equation}\label{eq:O}
\mathcal O(v_{00},v_{01})\mathcal O(v_{01},v_{11})=
-\mathcal O(v_{00},v_{10})\mathcal O(v_{10},v_{11}).
\end{equation}

Thus we obtain a co-chain complex 
$C_\mathcal O$ (respectively, $C_\mathcal A$)
whose basis is $P$ 
which is 
made by 
$\mathcal O(\quad,\quad)$ (respectively, $\mathcal A(\quad,\quad)$). 

We assume the following. 
\\

\h{\bf Cliam.}
{\it Let $n\in\N$. 
There is a chain isomprphism from 
$C_\mathcal O$ to $C_\mathcal A$ if $m=n$. }
\\

Therefore  
Claim \ref{seyade} is valid 
if $m=n$.  

We prove that Claim is true if $m=n+1$. 
Then 
Claim \ref{seyade} is valid for all $m\in\{0\}\cup\N$.   
\\


Let $P'$ (respectively, $P''$)
be a subposet of $P$, made of 
all (non-enhanced) Kauffman states with the vector 
$v_{n+1}=0$ (respectively, $v_{n+1}=1$).

Make cochain complexes  
$C_\mathcal O$,  
$C_\mathcal A$,  
$C'_\mathcal O$,  
$C'_\mathcal A$,  
$C''_\mathcal O$, and  
$C''_\mathcal A$.  
Note that 
$C'_\mathcal O$ and 
$C''_\mathcal O$ 
are sub-cochain complex of 
$C_\mathcal O$,  
and that 
$C'_\mathcal A$ and 
$C''_\mathcal A$ 
are sub-cochain complex of 
$C_\mathcal A$.

By the assumption of the induction, 
we have chain isomorphisms 
$f':C'_\mathcal O\to C'_\mathcal A$ and 
$f'':C''_\mathcal O\to C''_\mathcal A$. 

We construct 
$f:C_\mathcal O\to C_\mathcal A$ as below. 

If
 $\mathcal O((0,...,0), (0,...,1))= \mathcal A((0,...,0), (0,...,1))$,  
$f$ is defined by $f'$ and $f''$. 

If
 $\mathcal O((0,...,0), (0,...,1))= -\mathcal A((0,...,0), (0,...,1))$,  
$f$ is defined by $f'$ and $-f''$.  

Then $f$ is a chain isomorphism 
by the equations (\ref{eq:A}) and (\ref{eq:O}). 
 \qed\\

Let $K$ be a virtual link diagram. 
Make $\mathcal P(\alpha, \beta)$ by using the rule 
in Figure \ref{fig:canonical_source_sink_orientation}.  
We ask a question. 
Can we define a cochain complex  by using 
$\mathcal A(\alpha^\text{no}, \beta^\text{no})\mathcal P(\alpha, \beta)$? 
If so, is
the cochain complex chain isomorphic (or chain homotopy equivalent) to 
that defined 
by the coefficients 
$\mathcal O(\alpha^\text{no}, \beta^\text{no})\mathcal P(\alpha, \beta)$? 
We know that we can use 
$\mathcal A(\alpha^\text{no}, \beta^\text{no})$   
in the classical case, but in the virtual case we only know at present how to make a consistent definition using 
$\mathcal O(\alpha^\text{no}, \beta^\text{no})\mathcal P(\alpha, \beta)$.

\bigbreak
\section{An example of Khovanov basis  
of 
the Khovanov chain complex for a virtual link diagram}
\label{eK}


\noindent
We show an example of Khovanov basis 
of a Khovanov chain complex for a virtual link diagram in Figures \ref{zero1}and \ref{zero2}, and  
another one in Figures \ref{eK1}-\ref{eK4}.
We draw posets of decorated resolution configurations 
in Figures 
\ref{eK3} and \ref{eK4}.

\begin{figure}
\centering\includegraphics[width=80mm]{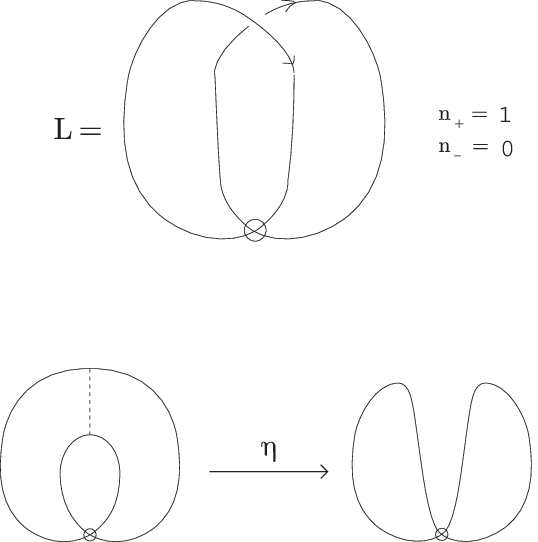}
\caption{{\bf A virtual link $L$, and 
the relation,  which is made by surgeries, 
among all (non-labeled) resolution configurations made from $L$.  
See Definition \ref{2.15gr} 
for the definition of $n_+$ and $n_-$.
$\eta$ denotes a single cycle surgery.}\label{zero1}}   
\end{figure}

\begin{figure}
\centering\includegraphics[width=133mm]{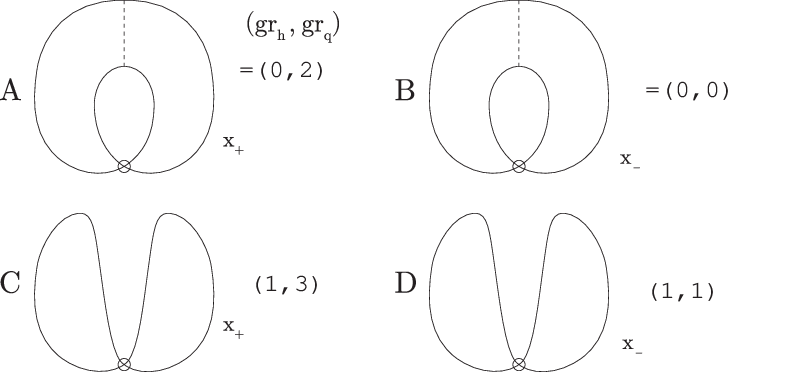}
\caption{{\bf The Khovanov basis 
of 
Khovanov chain complex for the virtual link diagram $L$ in Figure \ref{zero1}: 
They are all labeled resolution configurations made from $L$. 
$gr_h$ denotes the homological degree, and $gr_q$ the quantum one. 
We have 
$\delta(A)=\delta(B)=\delta(C)=\delta(D)=0$.
Recall Fact \ref{qde}, which is an important comment on a single cycle surgery. 
}\label{zero2}}   
\end{figure}

\begin{figure}
\centering\includegraphics[width=80mm]{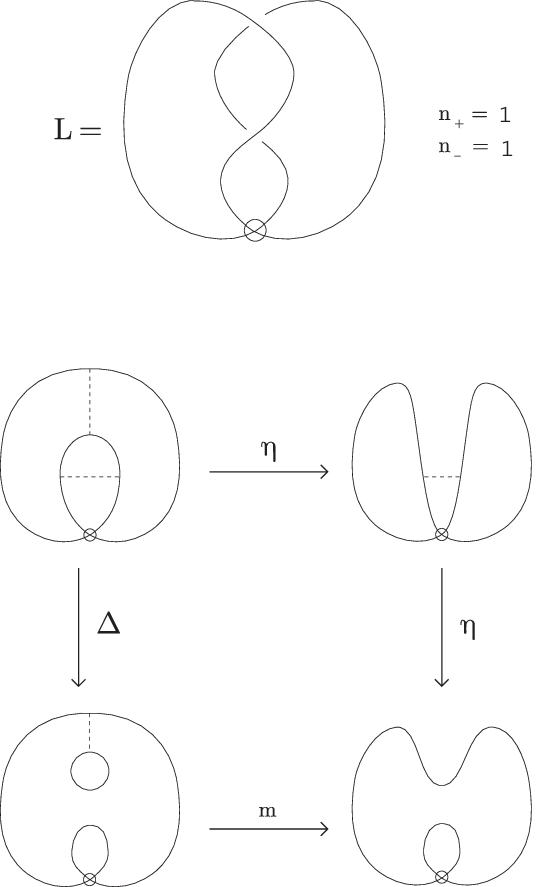}
\caption{{\bf A virtual link $L$, and 
the relation,  which is made by surgeries, 
among all (non-labeled) resolution configurations made from $L$.  
See Definition \ref{2.15gr} 
for the definition of $n_+$ and $n_-$.
$\eta$ denotes a single cycle surgery.}\label{eK1}}   
\end{figure}

\begin{figure}
\includegraphics[width=133mm]{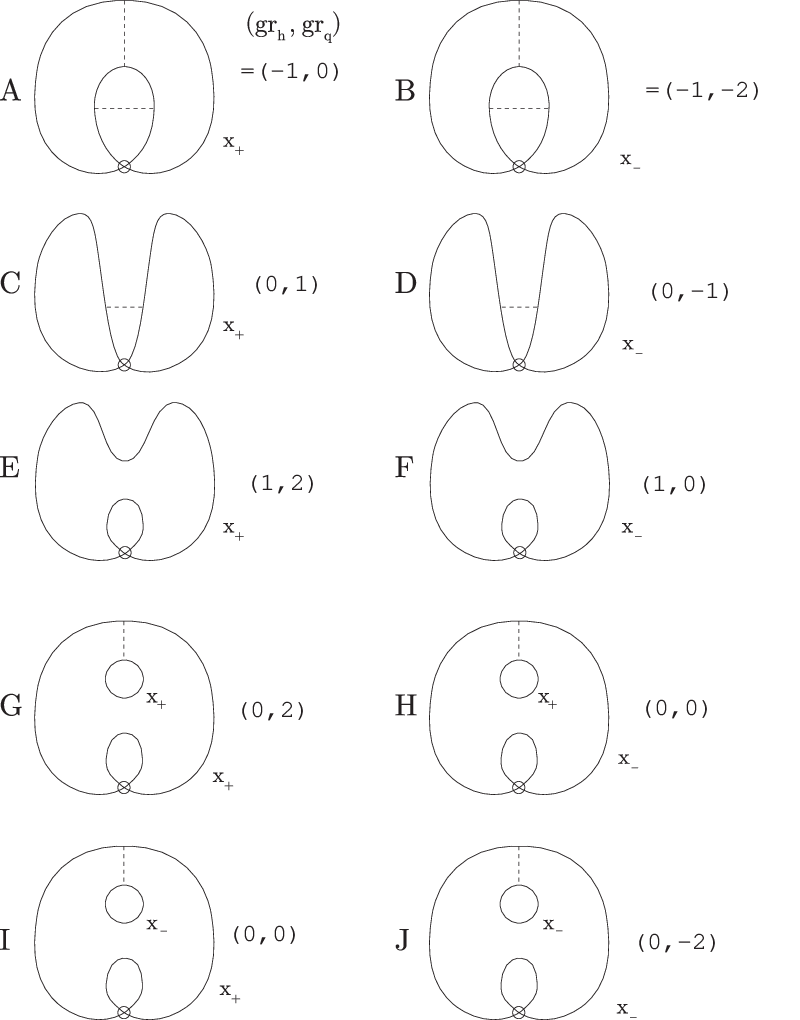}
\caption{{\bf The Khovanov basis 
of Khovanov chain complex for the virtual link diagram $L$ in Figure \ref{eK1}: 
They are all labeled resolution configurations made from $L$. 
$gr_h$ denotes the homological degree, and $gr_q$ the quantum one. 
}\label{eK2}}   
\end{figure}

\begin{figure}
\includegraphics[width=111mm]{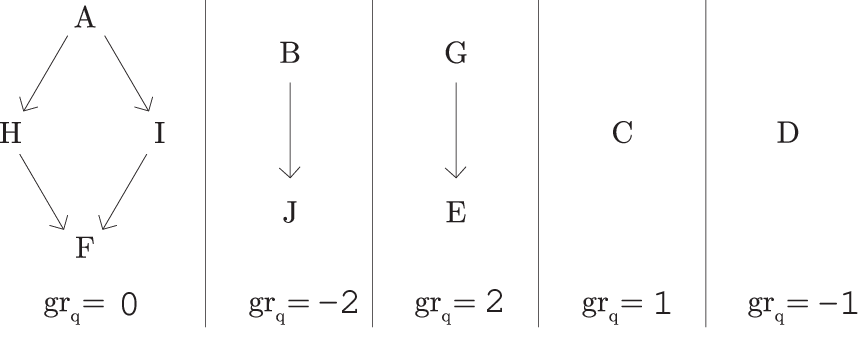}
\caption{{\bf 
The relation, which is made by surgeries, 
among all Khovanov basis elements 
in Figure \ref{eK1}. 
Let $P$ and $Q$ be labeled resolution configurations. 
If $c[P;Q]\neq0$, we connect $P$ and $Q$ by an arrow from $P $ to $Q$. 
There are five chunks. All Khovanov basis elements 
in each chunk have the same quantum grading.
Note that  if $c[P;Q]\neq0$, then $c[P;Q]=\pm1$. 
Recall that whether $+1$ or $-1$ is determined by using starting stars and  cut loci. 
Recall Fact \ref{qde}, which is an important comment on a single cycle surgery. 
}\label{eK3}}   
\end{figure}

\begin{figure}
\centering\includegraphics[width=111mm]{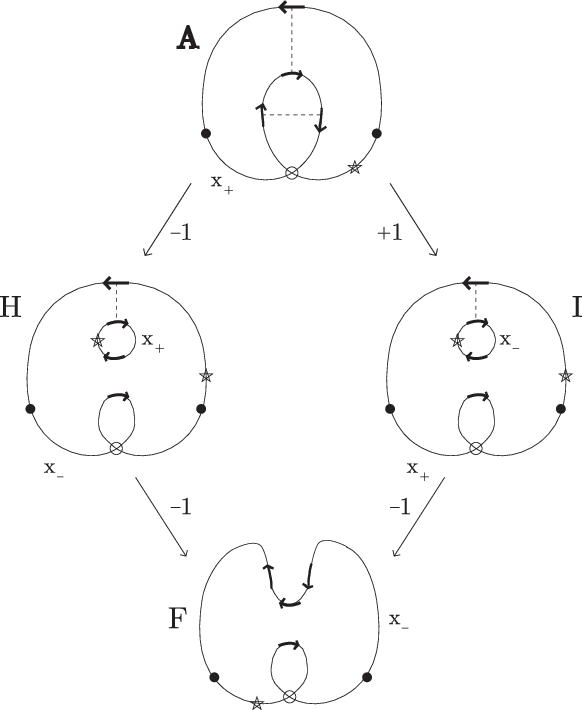} 
\caption{{\bf 
Examples of starting stars, cut loci, and $\zeta(\quad)$: 
By the above calculation, we know that $\delta\cdot\delta(A)=0$.  
}\label{eK4}}   
\end{figure}

\bigbreak
\part{Khovanov-Lipshitz-Sarkar CW complexes and 
the second Steenrod square for virtual links}\label{partsht}

\bigbreak
\section{Framed flow category: framings, mdulis, and CW complexes}\label{secfrcat}

As we announced in \S\ref{secffc}, 
we will make a CW complex from a Khovanov chain complex for a  given virtual link diagram.
We assign to a given chain complex, which is not  necessarily a Khovanov chain complex, 
a framed flow category (framings and modulis). 
For our purpose, we consider a stable homotopy type of CW complexes 
whose chain complex is a given one.

\begin{example}\label{exam}
Both 
$\Sigma^k$(the one point union of $S^2\cup S^4$) 
and 
$\Sigma^k(\C P^2)$,  
 where $\Sigma^k$ denotes the $k$-times suspension and $k$ is large, 
have a natural CW decomposition\\
 (the base point)$\cup e^{2+k}\cup e^{4+k}$. 
Consider a set of moduli spaces associated with \\
$\Sigma^k($the one point union of $S^2\cup S^4)$  
(respectively, $\Sigma^k(\C P^2)$).  
In $\partial e^{2+k}$, there is no moduli space. 
In $\partial e^{4+k}$, take an embedded circle. 
It is a moduli space. 
Take the normal bundle of the circle in  $\partial e^{4+k}$, 
and take the trivial (respectively, nontrivial) framing.  
It is a framing on the moduli space. 
 \end{example}

\begin{example}\label{nofr}
Regard $D^2$ as a union of the base point, a 1-cell $e^1$, and 2-cell $e^2$. 
We can regard $\Sigma^k(D^2)$ as 
a union of the base point, a $(k+1)$-cell $e^{k+1}$, and $(k+2)$-cell $e^{k+2}$. 
In $\partial e^{k+1}$, there is no moduli. 
In $\partial e^{k+2}$, take one point. It is a moduli space. 
The framing on the normal bundle of the point in $\partial e^{k+2}$ is unique 
after we give an orientation.
\end{example}

In short, for any given chain complex, 
once we have moduli spaces in each cell, and framings on them, we have a CW complex 
by the generalized Pontrjagin-Thom construction. 
Note that the chain complex may not be a Khovanov chain complex.
There is more than one way to associate a framed flow category 
(See 
Definition 3.12 of \cite{LSk} in page \pageref{pageda} of this paper, 
\cite[\S3.3 and \S4]{LSk} and \cite[\S3.3]{LSs}.)
to a chain complex in general. 
The cube moduli in \cite[\S4]{LSk}
for Khovanov chain complexes of classical links
lets Lipshitz and Sarkar  choose a specific way to make the framed flow category.



\bigbreak
\section{The ladybug configuration for classical link diagrams}\label{lady}

We review the ladybug configuration for classical link diagrams,  
which is introduced in \cite[section 5.4]{LSk}.
Lipshitz and Sarkar introduced it to define a CW complex for any classical link digram.   
We cite the definition of it, 
that of the right pair, and that of the left pair associated with it 
from \cite[section 5.4.2]{LSk}. \\

\begin{definition}\label{teten}  {\bf (\cite[Definition 5.6]{LSk}).}
An index 2 basic resolution configuration $D$ is said to be a ladybug configuration 
if the following conditions are satisfied (See Figure \ref{tento}.). 

$\bullet$ $Z(D)$ consists of a single circle, which we will abbreviate as $Z$;

$\bullet$ The endpoints of the two arcs in $A(D)$, say $A_1$ and $A_2$, 
alternate around $Z$ 

\hskip3mm (that is, $\partial A_1$ and $\partial A_2$ are linked in $Z$).\\

\end{definition}

\begin{definition}\label{rl}  {\bf (\cite[section 5.4.2]{LSk}).}
Let $Z$ denote the unique circle in $Z(D)$. 
The surgery $s_{A_1}(D)$ (respectively,  $s_{A_2}(D)$) consists of two circles; 
denote these $Z_{1,1}$ and $Z_{1,2}$ (respectively, $Z_{2,1}$ and $Z_{2,2}$); 
that is, $Z(s_{A_i}(D)) = \{Z_{i,1}, Z_{i,2}\}$.  
Our main goal is to find a bijection between 
$\{Z_{1,1}, Z_{1,2}\}$  and $\{Z_{2,1}, Z_{2,2}\}$; 
this bijection will then tell us which points in 
$\partial_{\rm exp}\mathcal M(x, y)$ to identify.

As an intermediate step, we distinguish two of the four arcs in 
$Z - (\partial A_1\cup \partial A_2)$. 
Assume that the point $\infty\in S^2$ is not in $D$, 
and view $D$ as lying in the plane $S^2-\{\infty\}\cong\R^2$. 
Then one of $A_1$ or $A_2$ lies outside $Z$ (in the plane) 
while the other lies inside $Z$. 
Let $A_i$ be the inside arc and $A_o$ the outside arc. 
The circle $Z$ inherits an orientation from the disk it bounds in $\R^2$. 
With respect to this orientation, each component of 
$Z - (\partial A_1\cup\partial A_2)$ 
either runs from the outside arc $A_o$ to an inside arc $A_i$ or vice-versa. 
The {\it right pair} is the pair of components of $Z-(\partial A_1\cup\partial A_2)$ 
which run from the outside arc $A_o$ to the inside arc $A_i$. 
The other pair of components is the {\it left pair}. See \cite[Figure 5.1]{LSk}.
\end{definition}

We explain why the ladybug configuration is important, below.

\begin{proposition}\label{4} 
Let ${\bf x}$ $($respectively, ${\bf y})$ be a labelled resolution configuration of 
homological grading $n$ $($respectively, $n+2).$
 Then the cardinality of the set  

\hskip3mm $\{p| p$ is a labelled resolution configuration.  
${\bf x}\prec p, p\prec{\bf y},$ $p\neq{\bf x}$, $p\neq{\bf y}\}$  

\noindent
is 0, 2, or 4, where $\prec$ is defined in Definition \ref{2.10}.
  \end{proposition} 

Let $D$ be the ladybug configuration.  
Give $D$ (respectively, $s(D)$) a labeling $x_+$ (respectively, $x_-$). 
We call the resultant labeled resolution configuration $(D,x)$ (respectively, $(s(D),y)$).  
The resultant decorated resolution configuration $(D,y,x)$ is called the 
the {\it decorated resolution configuration associated with the 
ladybug configuration $D$}. 
We draw the poset of  $(D,y,x)$ in Figure \ref{uen}. 

\begin{fact}\label{tenten}
The case of 4 in Proposition \ref{4} occurs when we have 
the decorated resolution configuration associated with the ladybug configuration. 
\end{fact}

Fact \ref{tenten}  is also explained in \cite[section 5.4]{LSk}.  

\begin{figure}
\centering\includegraphics[width=40mm]{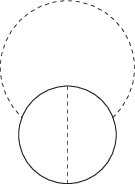}
\caption{{\bf The ladybug configuration 
}\label{tento}}   
\end{figure}

\begin{figure}
\includegraphics[width=155mm]{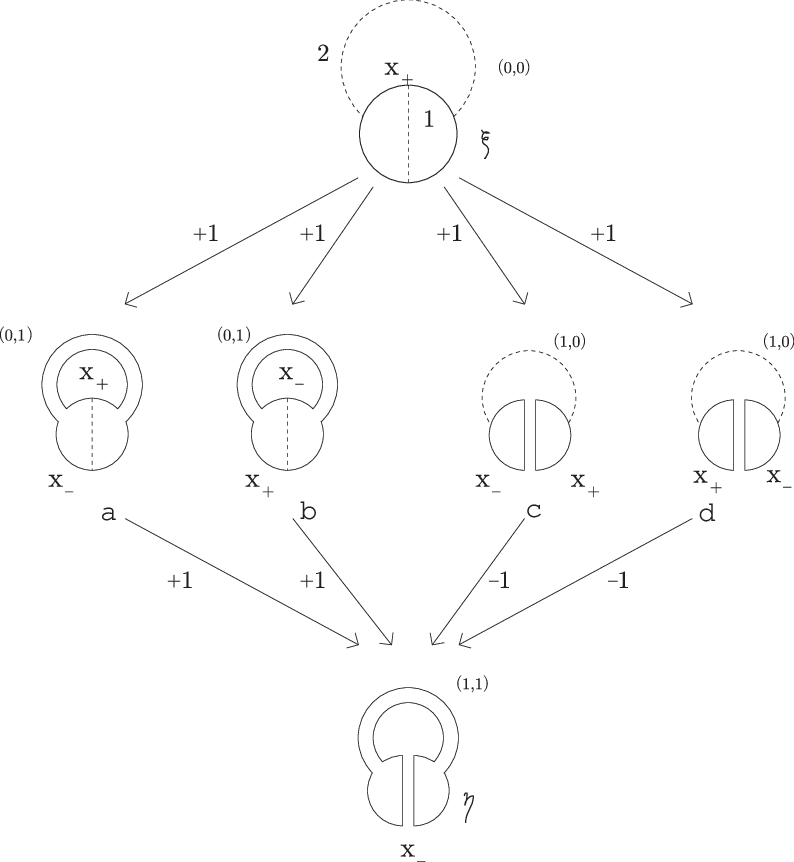}
\caption{
                 {\bf 
                        The boundary operator acting on the ladybug configuration:
                        The notation $a,b,c,$ and $d$ are defined in \cite[section 5.4.2]{LSk}. 
                        In \cite[Definition 2.2 and 2.15]{LSk}, 
                        $(*,\#)$ is defined associated with arcs.   
                        The numbers $+1$ and $-1$ denote the coefficient 
                        $[$the labeled resolution configuration at the arrowtail; that at the arrowhead$]$.
                  }\label{uen}
             }   
\end{figure}

\begin{proposition}\label{konnya}{\bf\rm(\cite[\S6]{LSk}.)}
The stable homotopy type of  the Khovanov-Lipshitz-Sarkar construction for classical link diagrams 
does not depend on 
whether we use the right pair or the left pair of each ladybug configuration 
in the classical link case when we construct Khovanov-Lipshitz-Sarkar CW complexes. 
\end{proposition}

\begin{proposition}\label{kya}  {\bf 
(This follows from results in \cite{LSk}. 
See the comments below.)}
The stable homotopy type of  the Khovanov-Lipshitz-Sarkar construction 
for classical link diagrams 
does not depend on 
the choice of  framing on modulis.
\end{proposition}

Of course Proposition \ref{kya} does not hold in the general case of construction of CW complexes.
 See Example \ref{exam}: 
 Framings change the stable homotopy types of CW complexes. 
 
However, Proposition \ref{kya}  is true in this case. 
Example \ref{nofr} is an example of a sub CW complex of a Khovanov CW complex. 
In Example \ref{nofr}, we use a framing when we construct a CW complex, but 
framings do not play such an important role, comparing Example \ref{exam}. 
 
Proposition \ref{kya} is the same as 
\cite[(4) in the first part of section six]{LSk}, which is proved in the proof of \cite[Proposition 6.1]{LSk}: 
In three lines above \cite[Definition 3.4]{LSs}, 
it is written, ``all such framings lead to the same Khovanov homotopy type
\cite[Proposition 6.1]{LSk}''.  
See also \cite[Lemma 4.13, which is cited in the proof of Proposition 6.1]{LSk}.

\bb

It is an outstanding property of 
the Khovanov chain complex 
and Khovanov stable homotopy type for classical links  
 that, 
if $\mathcal M_{\mathcal C_K(L)}(\bold{x},\bold{y})\neq\phi$,  
each connected component of 
$\mathcal M_{\mathcal C_K(L)}(\bold{x},\bold{y})$ 
is determined only by 
gr$_h\bold{x}-$gr$_h\bold{y}$.  
Chain complexes in other cases do not have this property in general.

Recall Remark \ref{remiron}. We mix the right pair and the left one. 

\section{Ladybug configurations and quasi-ladybug configurations 
 for 
 virtual link diagrams}\label{qlady}

\begin{figure}
\centering\includegraphics[width=50mm]{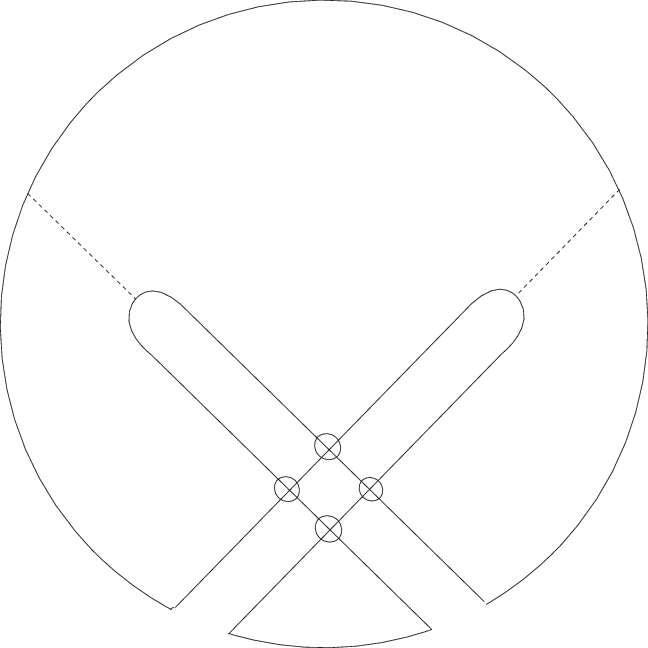} 
\caption{{\bf 
A quasi-ladybug configuration 
}\label{quasi}}   \end{figure}

\h
When we introduce  the Steenrod square for virtual links, 
 we also have the ladybug configuration. 
Furthermore, in the case of virtual links, 
we have quasi-ladybug configurations that is different from  the ladybug configuration. 
\\

Take a resolution configuration 
which is made of one immersed circle and two 
m-arcs.  \\

Stand at a point in the immersed circle 
where you see 
an arc to your right. 
Go ahead along the immersed circle. Go around one time.
Assume that you encounter the following pattern: 
 In the order of travel you next touch the other 
 arc.  
Then you touch the  first 
arc. 
Then you touch the other 
arc again.  
Finally, you came back the point at the beginning. \\

Since both 
arcs are m-arcs,
both satisfy the following property: 
At both endpoints of each arc, 
you see the 
arc in the same side -- either on  
the right hand side and on the left hand side.\\

If you see 
arcs 
both in the right hand side and in the left hand side 
(respectively, only in the right 
hand side)
while you go around one time, 
the resolution configuration is 
called a {\it ladybug configuration} 
(respectively, {\it quasi-ladybug configuration}). 

If a given resolution configuration has no virtual crossing, 
this definition of ladybug configurations is the same as 
that 
in \S\ref{lady}. 

No quasi-ladybug configuration appears 
if a given resolution configuration has no virtual crossing.  
However, a quasi-ladybug configuration may exist 
if a given resolution configuration has a virtual crossing.  
An example is drawn in Figure \ref{quasi}.

\bb
Let $D$ be a ladybug configuration. 
Let $C$ be the only one immersed circle in $Z(D)$.  
Cut $C$ at the four points where the 
arcs meet the endpoints. 
The immersed circle is then divided into four pieces. 
Recall that, at the beginning point of your trip, you see 
an arc on the right 
hand side. 
The first 
and third 
pieces of the four,  
which you are in while your trip, 
is the {\it right pair},  
and call the other two the {\it left pair.}
Note that the orientation of your trip 
and the place where you stand at the beginning of your trip
do not change the right and the left pair.
Note also that, 
if a given  resolution configuration does not have a virtual crossing, 
this definition is the same as that in \S\ref{lady}. 

It is important that 
we cannot determine the right and left pair 
in the case of quasi-ladybug configurations 
by this method.
(We pause to ask a question: Can one find a method to define the right and the left pair 
for quasi-ladybug configurations, 
to be compatible with the construction of Khovanov-Lipshitz-Sarkar stable homotopy type?)
\bb

Let $D$ be a ladybug configuration (respectively, quasi-ladybug configuration). 
Make $s(D)$.  
Give $D$ (respectively, $s(D)$) a labeling $x_+$ (respectively, $x_-$).
Call the resultant labeled resolution configuration 
$(D,x)$ (respectively, $(s(D),y)$). 
The decorated resolution configuration $(D,y,x)$
is called  the {\it decorated resolution configuration
associated with the ladybug configuration 
$($respectively, quasi-ladybug configuration$)$ 
$D$}.  
We draw the poset of an example in Figure \ref{uen} 
(respectively, Figure \ref{quasideco}).
$P(D,y,x)$ includes four labeled resolution configurations other than $(D,x)$ and $(s(D),y)$. 
We have two ways to make a 2-dimensional CW complex associated with $P(D,y,x)$.

In the case of ladybug configurations, 
by using the right and left pairs introduced above, 
we determine the right and left pair of the labeled resolution configurations 
in the middle raw of 
the poset of the decorated resolution configuration associated with 
a given ladybug configuration 
as in \cite[Figure 5.1 and its explanation in \S5.4.2]{LSk}. 

Recall Remark \ref{remiron}. 
In  Lipshitz and Sarkar's paper \cite{LSk} 
it is important 
how to assign  to 
 the ladybug configuration a moduli. 
They give all ladybug configurations one of the right and the left pair. 
That is,  all ladybug configurations have the same pair.  
In our paper, we consider the following situation: 
Each ladybug configuration may have different pairs. 
We consider all cases. Therefore we may consider more than one set of modulis for one Khovanov chain complex. 
%
We discuss 
both 
the condition that Lipshitz and Sarkar\cite{LSk} discussed 
and 
the condition that Lipshitz and Sarkar\cite{LSk} did not discuss.
See also Proposition \ref{konnya} and Remark \ref{tsun}.

%
%


However, 
in the case of quasi-ladybug configurations, 
we cannot distinguish two cases.  
Therefore we make two CW complexes 
if there is a quasi-ladybug configuration. 
This means that we consider a multiplicity of Khovanov homotopy types and take the collection of them as an invariant,  
as we explain below.

\begin{figure}
\includegraphics[width=170mm]{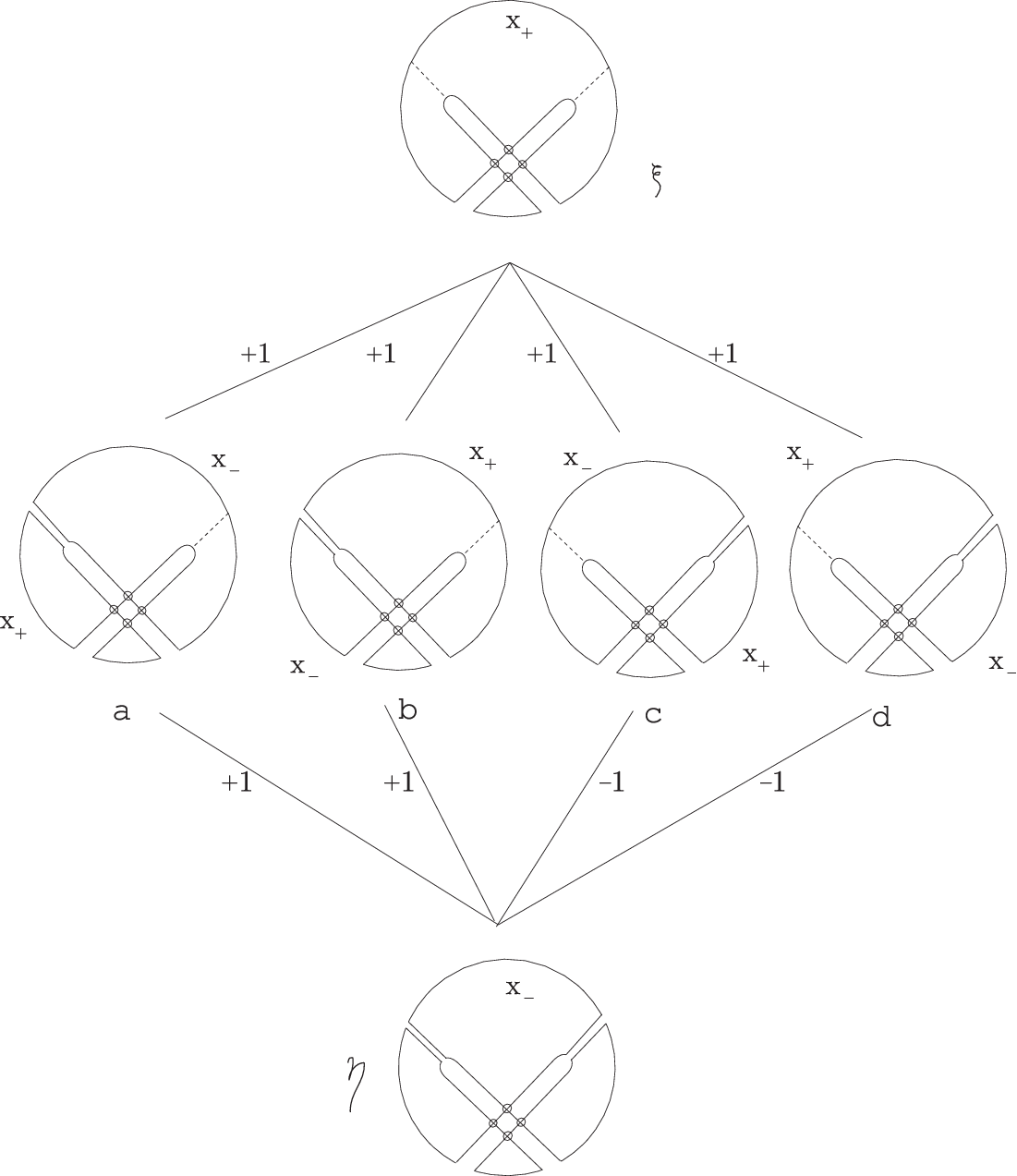} 
\caption{{\bf 
The decorated resolution configuration associated with 
a quasi-ladybug configuration 
}\label{quasideco}}   \end{figure}


Recall Remark \ref{remiron}. We mix the right pair and the left one.

\begin{remark}\label{kurage} 
In both Figures \ref{uen} and \ref{quasideco},
we have the following identities of coefficients in the differentials.

\h
$[\xi;a][a;\eta]=[\xi;b][b;\eta]=-[\xi;c][c;\eta]=-[\xi;d][d;\eta]$ and 
 $[\xi;\ast],[\ast;\eta],[\xi;\ast],[\ast;\eta]\in\{+1,-1\}$, 
where $\ast=a,b,c,d$.

\h{\it Reason}. 
The vector (in Definition \ref{2.2}) of $a$ and that of $b$ are the same. That of $c$ and that of $d$ are the same. That of $a$ (respectively, $b$) is different from that of $c$ (respectively $d$). \\

Therefore 
the moduli of the decorated resolution configurations in Figures \ref{uen} and \ref{quasideco} 
is the disjoint union of two segments, and  
its boundary is $a\amalg b\amalg c\amalg d$. 
We have just two cases. 

(i) One segment connects $a$ and $c$, and    
the other $b$ and $d$. 

(ii) One segment connects $a$ and $d$, and    
the other $b$ and $c$. 

\h 
We never have the following case: 
One segment connects $a$ and $b$, and    
the other $c$ and $d$. 

\h{\it Reason.} 
If we choose ``$a\amalg b$, and  $c\amalg d$'', then $\delta\circ\delta\neq0$ by the definition of $\delta$ 
and therefore $\partial\circ\partial\neq0$, 
associated with each moduli.  
\end{remark}
\bb

In the case of ladybug configurations, 
we can distinguish the cases (i) and (ii) above, by using the right and left pairs. 
However, in the case of quasi-ladybug configurations, 
we cannot distinguish them.  
As we wrote above, this means that, in general,  
we may associate more than one CW complex to a single virtual link diagram. 
See the following sections. 


\section{Why is it more difficult to define Khovanov-Lipshitz-Sarkar stable homotopy type 
for virtual links than for classical links?} \label{secs1}\hskip10mm

In the following sections, 
we will define a second Steenrod square for virtual links 
that is stronger than Khovanov homology for virtual links. 
We use a CW complex associated 
with a dual Khovanov chain complex. 
\\

The presence of the single cycle zero map in the virtual Khovanov chain complex 
and $\mathcal P(\alpha, \beta)$ 
in the coefficient of the virtual Khovanov chain complex  
take us out of the cube complex method for defining a framed flow category. For this reason we use a truncated homotopy type for this paper.

In the case of classical link diagrams,
for Khovanov basis elements 
$\bold{x}=(D_L(u),x)$ and $\bold{y}=(D_L(v),y)$), 
of the Khovanov chain complex of a classical link diagram $L$,  
we can assign to the moduli space 

$$\mathcal M_{\mathscr C_K(L)}(\bold{x},\bold{y})=\mathcal M(D_L(v)-D_L(u),x|, y|),$$

\h a disjoint union of the gr$_h \bold x-$gr$_h \bold y$ dimensional cube moduli.  
Here,  $\bold y\prec \bold x$. 
On the other hand, 
in the case of virtual link diagrams,   
we cannot assign to 
the moduli space \\
$\mathcal M_{\mathscr C_K(L)}(\bold{x},\bold{y})=\mathcal M(D_L(v)-D_L(u),x|, y|)$
 an $m$-dimensional cube, in general. \\
 
%
\begin{figure}
\centering\includegraphics[width=70mm]{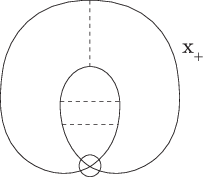}
\caption{{\bf 
A labeled resolution configuration
}\label{rei4}}   
\end{figure}

An example is shown in Figures \ref{rei4}-\ref{rei2}.
 A labeled resolution configuration $A$ is drawn in Figure \ref{rei4}.  
The sequence of the labeled resolution configurations starting from $A$  
is drawn in Figure \ref{rei1}.
Of course this sequence is made by surgery along arcs. 
A decorated resolution configuration 
associated with $A$ is drawn 
in Figure \ref{rei3}. 
The poset of the decorated resolution configuration 
is drawn in Figure \ref{rei2}.
Note that Figures \ref{rei1} and  \ref{rei2} are different.

\begin{figure}
\centering\includegraphics[width=150mm]{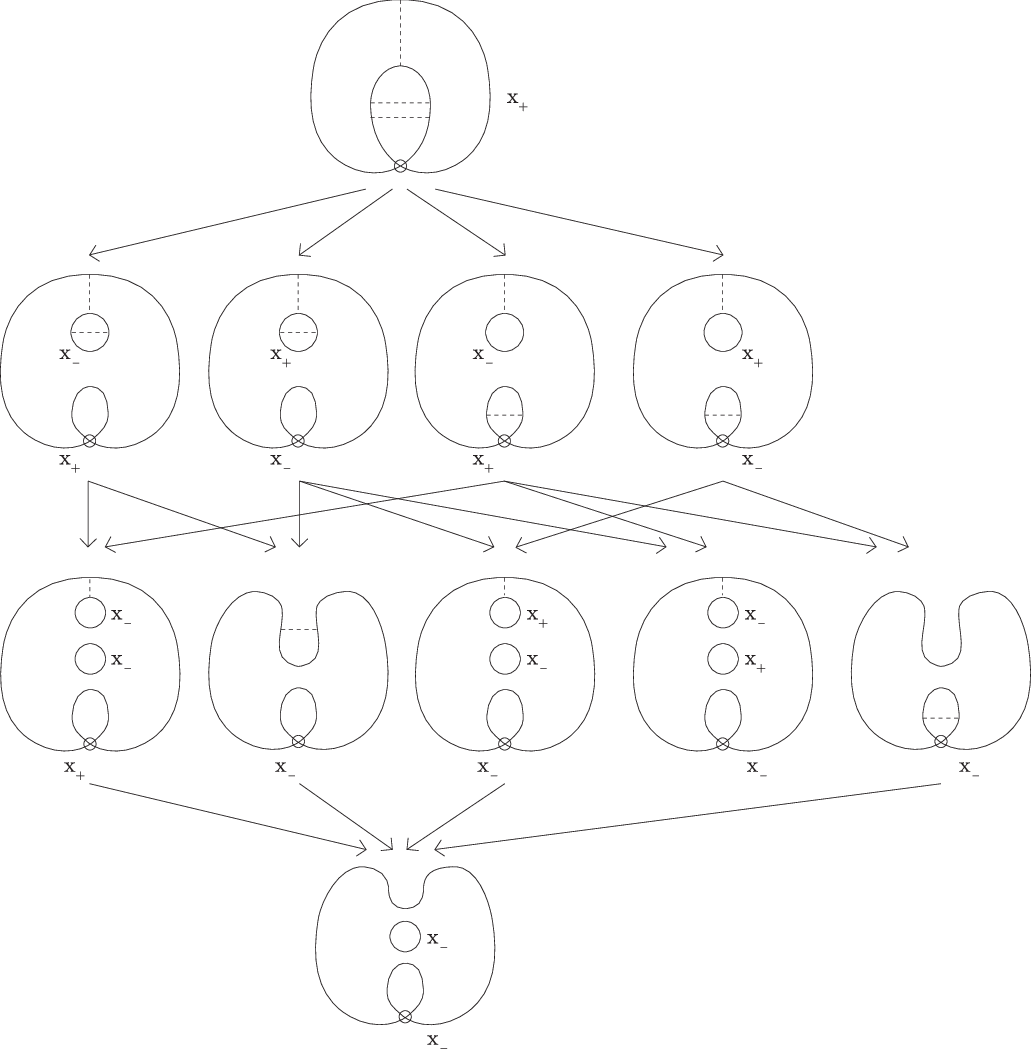}
\caption{{\bf 
The sequence of 
the labeled resolution configurations made from a labeled resolution configuration
in Figure \ref{rei4}}\label{rei1}}   
\end{figure}

\begin{figure}
\includegraphics[width=140mm]{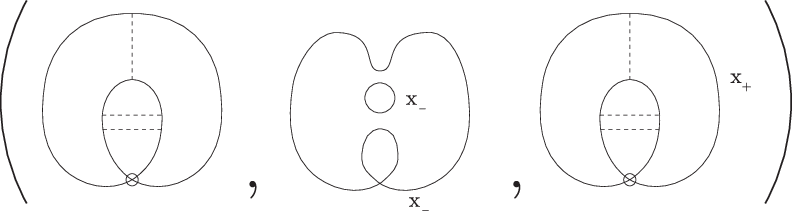}
\caption{{\bf 
A decorated resolution configuration
}\label{rei3}}   
\end{figure}

\begin{figure}
\includegraphics[width=140mm]{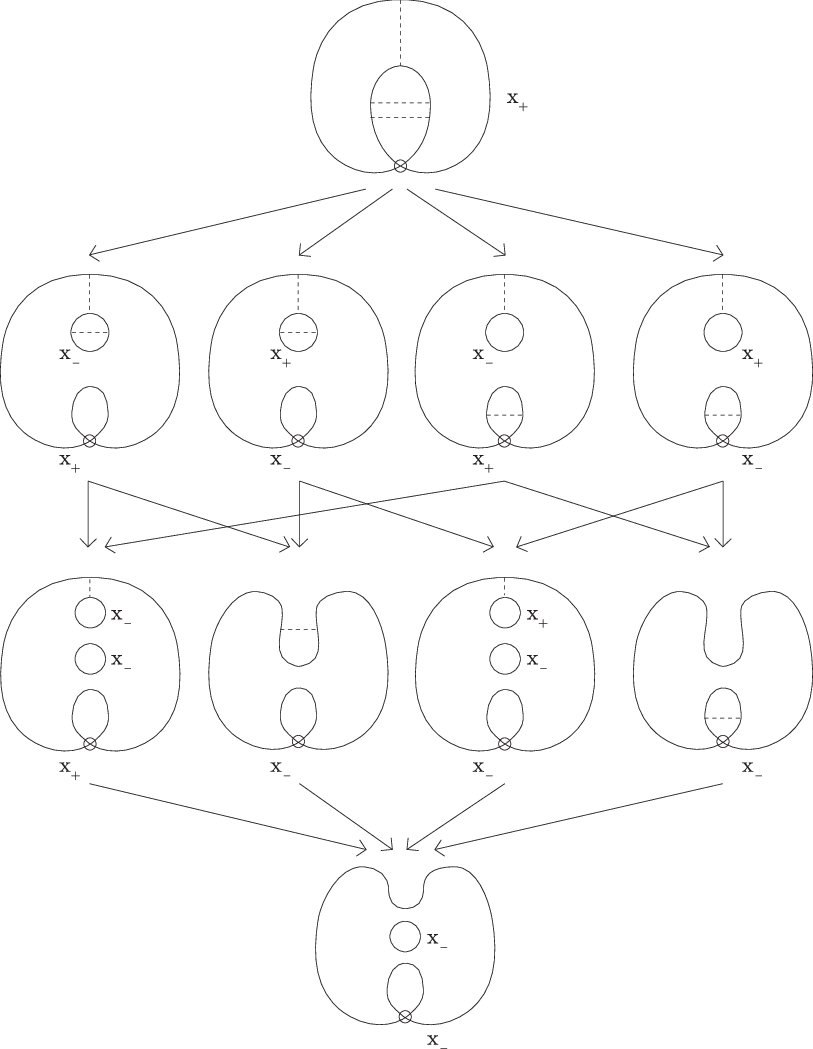}
\caption{{\bf 
The poset 
of the decorated resolution configuration in Figure \ref{rei3}
}\label{rei2}}   
\end{figure}

\begin{figure}
\vskip-10mm
\includegraphics[width=140mm]{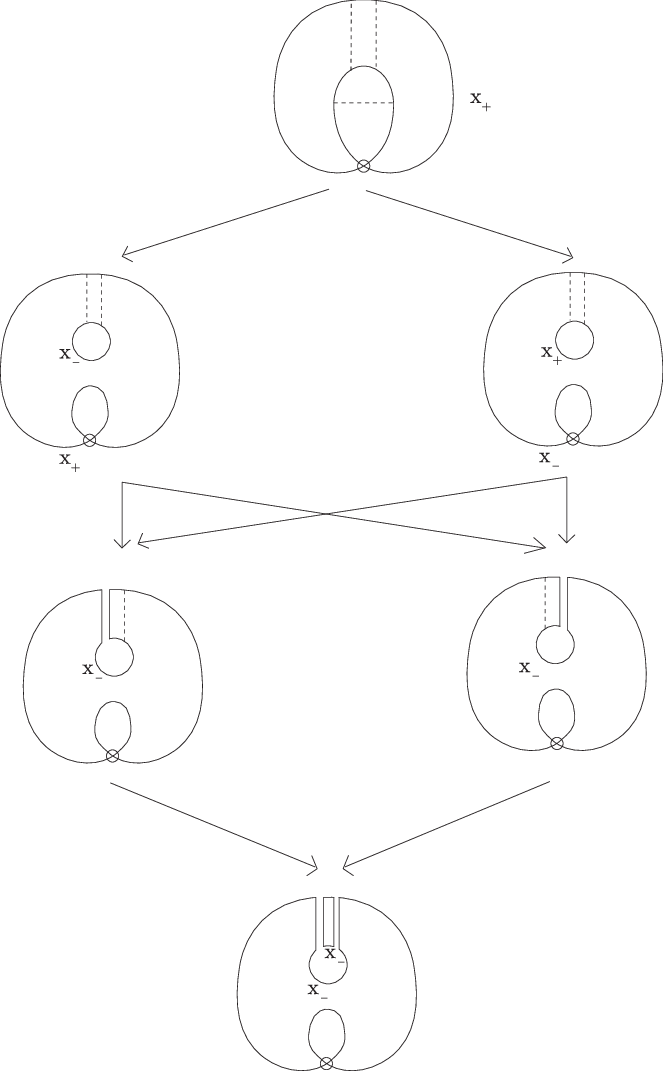}
\caption{{\bf 
The poset 
of the decorated resolution configuration of an index 3 labeled resolution configuration
}\label{mouikko}}   
\end{figure}

In \cite[section 5.1]{LSk} it is proved  that,  
if Ind$(D,x,y)$=3,  $\mathcal M(D,x,y)$ is a disjoint union of the 3-dimensional cube moduli 
in the case of classical link diagrams. 
However, Figure \ref{rei2} indicates the following:  
 Ind$(D,x,y)$=3 does not imply that  $\mathcal M(D,x,y)$ is 
 a disjoint union of the 3-dimensional cube moduli in general. 
Under this condition, 
if  Ind$(D,x,y)>3$, the situation would be more difficult. 
Therefore it is more difficult to define moduli spaces for Khovanov chain complex 
in the case of virtual links 
than in  that of classical links. 

We draw another poset in Figure \ref{mouikko}: 
There is not a natural onto map from this poset to $\mathcal C(3)$. 
We can assciate to this poset a moduli homeomorphic to the 2-disc, 
but it is not the 3-cube moduli $\mathcal M_{\mathcal C(3)}(1,0)$.
\\

We have two other reasons: 
Remark \ref{kinoko} and the comment below Remark \ref{kinoko}. 

\begin{remark}\label{kinoko}
See Figure \ref{eK4}.
In the case of virtual link diagrams with virtual crossing points, 
we have the following.   
We have $[A;H][H;F]=-[A;I][I;F]=1$, and 
we have the vector of $H$ is the same as $I$. 


 In Figure \ref{eK4}, the moduli space $\mathcal M(F, A)$ is one segment.  
 Consider a map from this moduli to the 2-dimensional cube moduli which is written in \cite[\S5.1]{LSk}, and say $\pi$ . 
 The image of the boundary, just two points,  of $\mathcal M(F, A)$, by this map  is one point because the vector of $H$ is the same as $I$. 
We do not have this phenomenon in the case of classical link diagrams. 

Let $x,p,q$ and $y$ be labelled resolution configurations for a virtual link diagram $D$ of a virtual link $L$. Suppose that $[x;p][p;y]=-[x;q][q;,y]=1$. If $D$ does not have a virtual crossing, the vector (in Definition \ref{2.2}) of $p$ is different from that of $q$. On the other hand, if $D$ has a virtual crossing, the vector of $p$ is the same as that of $q$ in some cases as written right above, 
and  the former is different from the latter in the other cases. 
\end{remark}

As we stated in \S\ref{qlady}, 
we may have a quasi-ladybug configuration in general, 
in the case of non-classical link diagrams.  Then 
we may associate more than one CW complex to a single virtual link diagram in general.

\bb
\section{Our strategy of the construction of the second Steenrod square 
for virtual links
}\label{s2}\hskip9mm

We explain the detail of 
our strategy that we announced 
in \S\ref{secstr}

It is easy to prove that 
a dual Khovanov chain complex always associates no less than one CW complex 
(see for example \cite[Theorem in Exercise 4,  section 39, page 231]{Munkres}). 
However, we do not know whether we can define the moduli space consistently 
by only the information of the relations of 
Khovanov basis elements 
in general, in the virtual link case. 

See \cite[Definition 5.5]{LSk}: 
 Lipshitz and Sarkar use the $(n-1)$ cube moduli, which we call the `generating moduli' 
 for a classical link $n$ crossings 
 and define all moduli spaces  and framings consistently for 
Khovanov chain complex of all classical link diagrams, 
and define the CW complex. 
They proved that the stable homotopy type of it is invariant under all Reidemeister moves. \\

Take a set of all CW complexes which are associated with 
the dual Khovanov chain complex for an arbitrary virtual link diagram. 
In fact, this set may be 
an invariant for virtual links if we can prove the invariance under Reidemeister moves. 
However, even if so, we do not know whether this invariant is stronger than 
Khovanov homology, nor whether we can calculate this invariant.\\

In the virtual case, it is very complicated to construct a `generating moduli' 
because of the property of coefficients in the virtual case (Definition \ref{korekoso}) 
is different from the classical case. 
In this paper  we show an explicit way to assign to the moduli space 
$\mathcal M_{\mathscr C_K(L)}(\bold{x},\bold{y})\\=\mathcal M(D_L(v)-D_L(u),x|, y|)$, 
a compact topological space which admits a CW structure,  
when  {\rm gr}$_h \bold x-${\rm gr}$_h \bold y$ is 
one, two, and 
three.  
We construct a CW complex which consists of 
only 
$(m-1)$-cells, 
 $m$-cells, 
$(m+1)$-cells,  
$(m+2)$-cells, and 
$(m+3)$-cells, 
where $m$ is any integer,    
for the dual Khovanov chain complex in this case by using these moduli spaces. 
We prove that the second Steenrod square of the CW complex is invariant under any Reidemeister move 
although we do not prove 
whether the stable homotopy type of the CW complex is invarint under any Reidemeister move.  

We have not constructed modulis 
when  
{\rm gr}$_h \bold x-${\rm gr}$_h \bold y\geqq5$, 
to be compatible with those in the case of  
{\rm gr}$_h \bold x-${\rm gr}$_h \bold y\leqq4$. 

\section{0-dimensional modulis  and framings }\label{secmo0}
Let $\mathcal L$ be a virtual link. 
Let $L$ be a virtual link diagram which represents $\mathcal L$. 
Let $\bold{x} =(D_L(u),x)$ and $\bold{y}=(D_L(v),y)$ be  
Khovanov basis elements. 
We will make moduli spaces 
$\mathcal M_{\mathscr C_K(L)}(\bold{x},\bold{y})=\mathcal M(D_L(v)-D_L(u),x|, y|)$ 
when 
gr$_h\bold x-$gr$_h\bold y=1,2,3,4$.  
In this section we suppose that the modulis are not the empty set. 
Of course, this discussion includes the case of classical link diagrams because 
any classical link diagram is a virtual link diagram.

We have the following proposition.

\begin{proposition}\label{ten}
If {\rm gr}$_h \bold x-${\rm gr}$_h \bold y=1$,  we can assign to 
the moduli space $\mathcal M_{\mathscr C_K(L)}(\bold{x},\bold{y})$ 
a single point.  
\end{proposition}

\noindent{\bf Proof of Proposition \ref{ten}.}
The proposition follows because $c[x:y]$ is $+1, 0,$ or $-1.$ 
\qed

\begin{remark}\label{n1}
In the case of classical links, Proposition \ref{ten} is the same as \cite[section 5.3]{LSk}.  
\end{remark}

We give a framing on $\mathcal M_{\mathscr C_K(L)}(\bold{x},\bold{y})$ 
so that it satisfies the conditions  
in \cite[Definitions 3.18 and 3.20]{LSk}.


\section{1-dimensional modulis  and framings }\label{secmo1}

We have the following.

\begin{lem}\label{2or4}
Let $\bold a, \bold b$ be two dual Khovanov basis elements 
such that the difference of the homological gradings are two. 
Let $\#=\#\{\bold x|\bold a\prec \bold x, \bold x\prec \bold b, 
x\neq a$, $x\neq b\}$. 
If $\#\neq0$, it is two or four. 
\end{lem}

\noindent{\bf Proof of Lemma \ref{2or4}.}
If $\bold a$ and $\bold b$ are (respectively, are not) associated with 
a ladybug configuration or a quasi-ladybug configuration, 
$\#$ is four (respectively, two). 
\qed\\

It is easy to prove the following proposition.   

\begin{proposition}\label{hen}
Let {\rm gr}$_h \bold x-${\rm gr}$_h \bold y=2$. 
Under the conditions of Proposition \ref{ten}, we have just three cases.

\bs\h$(1)$
We can assign to    
the moduli space 
$\mathcal M_{\mathscr C_K(L)}(\bold{x},\bold{y})$ 
one closed segment 
if $\#$ in Lemma $\ref{2or4}$ is two.


\bs\h$(2)$
We can assign to    
the moduli space $\mathcal M_{\mathscr C_K(L)}(\bold{x},\bold{y})$ 
two kinds of  a disjoint union of two closed segments 
if $\#$ in Lemma $\ref{2or4}$ is four. 
\end{proposition}

\h{\bf Remark.} 
In the case of virtual links, we may have both ladybug configurations and quasi-ladybug configurations.  
Note that quasi-ladybug configurations are different from ladybug configurations. 
Recall \S\ref{qlady}.
We must note that, in the case of classical links, 
there is not a quasi-ladybug configuration.   
\\

We use Lemma \ref{2or4}, and prove Proposition \ref{hen}   
as in the case of classical links in \cite[section 5.4]{LSk}. 
Proposition \ref{hen} in the case of classical links includes \cite[section 5.4]{LSk}.  
They are not the same. See the following Remark. 

\begin{remark}\label{tsun}
As we stated in \S\ref{qlady}, 
if we have a quasi-ladybug configuration, we may assign to a single virtual link diagram 
more than one CW complex.  
As we stated in  Remark \ref{remiron}, 
if we have a ladybug configuration, we may also assign to a single virtual link diagram  
more than one CW complex.

\end{remark}

Compare \cite[The property associated with (2) in the proof of Proposition 6.5]{LSk} with 
the following.

\begin{fact}\label{arere} 
{\it Fix a link diagram $L$, and let $L'$ be the result of reflecting 
$L$ across the $y$-axis, say, and reversing all of the classical crossings. 
Then $L$ and $L'$ do not always represent the same virtual link. }
\end{fact} 
\bb

The boundary of one segment in $\mathcal M_{\mathscr C_K(L)}(\bold{x},\bold{y})$ 
is two points.  The framings on two points have different orientations. 
Therefore we can extend the framing on the boundary to 
$\mathcal M_{\mathscr C_K(L)}(\bold{x},\bold{y})$. 
Take a framing on $\mathcal M_{\mathscr C_K(L)}(\bold{x},\bold{y})$. 
See \cite[\S4]{LSk}.

\section{2-dimensional modulis  and framings }\label{secmo2}

\begin{rev}\label{hexa}
The 3-cube moduli $\mathcal M_{\mathcal C(3)}(\bar1,\bar0)$  is a hexagon, 
and is homeomormohic to the 2-disc. See \cite[\S4.1]{LSk}. 

In Figure \ref{roku}, we draw 
the partial ordered set $\mathcal C(3)$, the set of six vertices of the 3-dimensional cube, 
and cells associated with the six vertices: Here we show an example 
the cells are 3, 2, 1, and 0-dimensional case for convenience.
This moduli space  $\mathcal M_{\mathcal C(3)}(\bar1,\bar0)$  is in 
the 2-sphere, $\R^2\cup\{\infty\}$, which is the boundary of the 3-cell $x$. 
We draw 
$\mathcal M_{\mathcal C(3)}(\bar1,\bar0)$ 
and 
$\partial\mathcal M_{\mathcal C(3)}(\bar1,\bar0)$ 
in this $\R^2$ in 
Figures \ref{rokkaku} and \ref{rokkaku2}. 








\begin{figure}
\includegraphics[width=100mm]{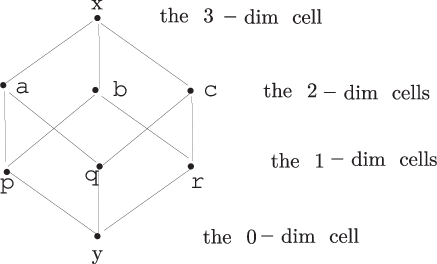}
\caption{{\bf Cells associated with 
the partial ordered set $\mathcal C(3)$, 
the set of six vertices of the 3-dimensional cube    
}\label{roku}}   
\end{figure}

\begin{figure}
\centering\includegraphics[width=150mm]{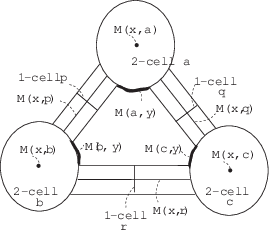}
\caption{{\bf 
The 3-dimensional cube moduli. 
Take only the hexagon  $\mathcal M_{\mathcal C(3)}(1,0)$ from Figure \ref{rokkaku}.  
}\label{rokkaku2}}   
\end{figure}

\begin{figure}
\centering\includegraphics[width=150mm]{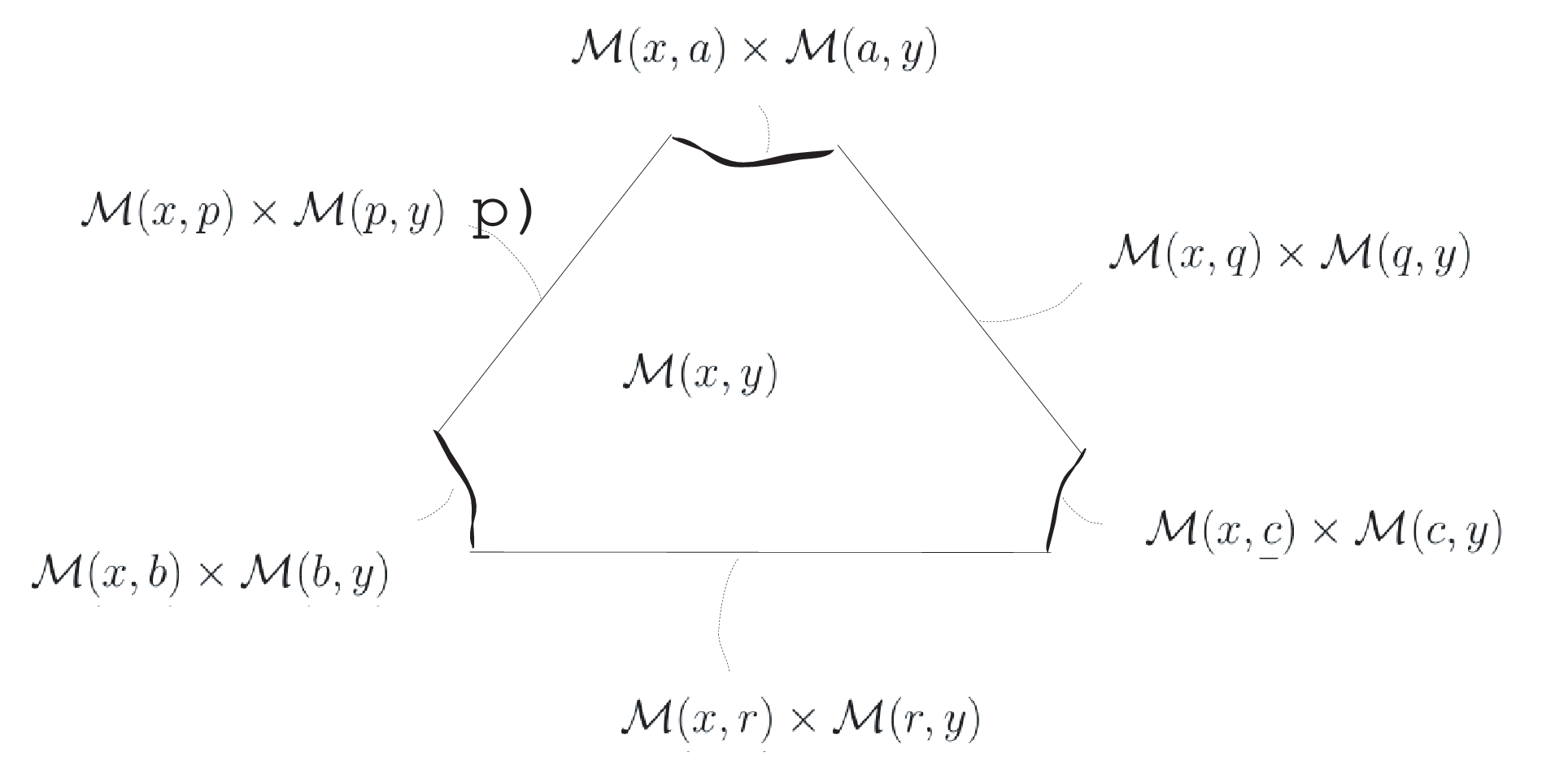}
\caption{{\bf 
The 3-dimensional cube moduli. See also Figure \ref{rokkaku2}
}\label{rokkaku}}   
\end{figure}
\end{rev}

We have the following proposition when  {\rm gr}$_h \bold x-${\rm gr}$_h \bold y=3$.

\begin{proposition}\label{rippo}
Let {\rm gr}$_h \bold x-${\rm gr}$_h \bold y=3$. 
Under the conditions of Propositions $\ref{ten}$ and $\ref{hen}$, 
we can assign to 
each moduli space $\mathcal M_{\mathscr C_K(L)}(\bold{x},\bold{y})$,  
a space homeomorphic to a finite number of  2-discs. 
\end{proposition}

\begin{remark}\label{n3} 
In Proposition \ref{rippo}, 
in the case of virtual links 
$\mathcal M(D,x,y)$  is not a cube flow category 
in general 
as is shown  
by the example of Figures \ref{rei4}-\ref{rei2} in \S\ref{secs1}.  
Consider 
 $\mathcal M(D,x,y)$ associated with the example of Figures \ref{rei4}-\ref{rei2}. 
  We can check directly that 
$\partial\mathcal M(D,x,y)$ is homeomorphic to a disjoint union of circles, 
and can suppose that 
  $\mathcal M(D,x,y)$ is homeomorphic to a disjoint union of 2-discs.

In the case of classical links,  Proposition \ref{rippo} is the same as \cite[section 5.5]{LSk}. 
\end{remark}

\begin{figure}
\centering\includegraphics[width=100mm]{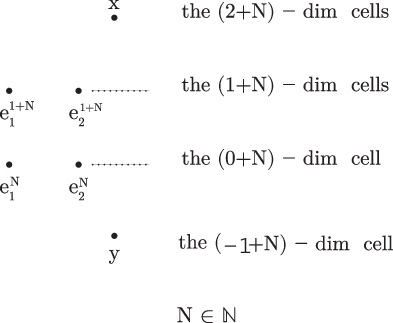}
\caption{{\bf  Khovanov basis elements 
in Proposition \ref{rippo} 
}\label{8}}   
\end{figure}


\noindent{\bf Proof of Proposition \ref{rippo}.}  
Make a basic index 3 decorated resolution configuration 
 $(D,x,y)$ associated with  
$\bf x$ and   $\bf y$  
by the same method of \cite[Definition 5.3]{LSk}. 
$P(D,x,y)$ 
has the maximal element $(D,y)$ 
and the minimal element $(s(D), x)$. 

We assign a cell to each element of   $P(D,x,y)$ as drawn in Figure \ref{8}: 
$\bf x$ (respectively, $\bf y$) corresponds to $(s(D), x)$ (respectively, $(D,y)$). 
Assume that 
(the homological degree of $\bf x$) 
$-$
that of $\bf y$ is $3$. 
There are $N$-cells $e^{N}_i$ and 
$(N+1)$-cells $e^{N+1}_j$ between $\bf x$ and $\bf y$. 
Here,  $N$ is a large integer, $i=1,...,\nu_i$ and $j=1,...,\nu_j$. 

By the definition, ${\bf y}=(D,y)$ has just three arcs.

If all arcs in each labelled resolution configuration of $P(D,x,y)$ are mc arcs, 
we have the same results as ones in \cite[section 5.5]{LSk}. 
In this case Proposition \ref{rippo} holds. 
\\

If all arcs in $(D,y)$ are scs arcs,    
$\delta (D,y)=0$. However $(D,y)$ is an element of $P(D,x,y)$ 
which is not the empty set.  
We arrived at a contradiction.   
Hence this case does not occur. 
\\

We prove the other cases.

Note that, 
even if 
$(D,y)$ has only three mc arcs, 
a labeled resolution configuration $e^{1+N}_*$ may have a scs arc. 
See Figure \ref{mcscs} for an example.

\begin{figure}
\includegraphics[width=140mm]{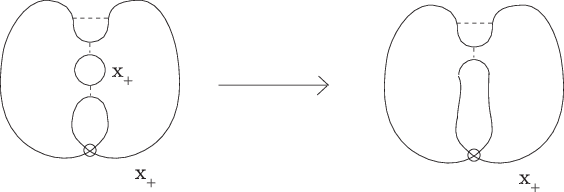}
\caption{{\bf 
An example of surgery which changes 
a labelled resolution configuration with only mc arcs 
into 
a labelled resolution configuration with a scs arc. 
}\label{mcscs}}   
\end{figure}

If we can construct $\mathcal M({\bf x,y})$,  
by \cite[Definition 3.12.(M-2)]{LSk},   
$\partial\mathcal M({\bf x,y})$ is 
a union of 

$$\displaystyle\coprod_{i=1,...,\nu_i} 
\mathcal M(e^{N}_i,{\bf y})\x\mathcal M({\bf x},e^{N}_i)$$
$$(=\partial_1\mathcal M({\bf x,y}))$$

\h and 

$$\displaystyle\coprod_{j=1,...,\nu_j} 
\mathcal M(e^{N+1}_j,{\bf y})\x\mathcal M({\bf x},e^{N+1}_j) $$
$$(=\partial_2\mathcal M({\bf x,y})).$$

Each of 
$\mathcal M(e^{N}_i,{\bf y})\x\mathcal M({\bf x},e^{N}_i)$
and  
$\mathcal M(e^{N+1}_j,{\bf y})\x\mathcal M({\bf x},e^{N+1}_j)$
is a  segment. 

We have that 
the set of the boundary of each segment above is 
the set of points, 
$\mathcal M(e^{N+1}_i,e^{N}_j)$ ($i=1,...,\nu_i$ and $j=1,...,\nu_j$). 


Furthermore, we have  the following: 
Fix $i$ and $j$.  
There is an only one segment in 
$\partial_1\mathcal M({\bf x,y})$ 
(respectively, $\partial_2\mathcal M({\bf x,y})$)
which touches a point $\mathcal M(e^{N+1}_i,e^{N}_j)$. 
These two segments touch each other at the point.

Therefore $\partial\mathcal M({\bf x,y})$ is a disjoint union of circles.

Therefore we can define as follows: $\mathcal M(x, y)$  
is a finite disjoint union of CW complexes 
which are homeomorphic to the 2-ball. 
This  $\mathcal M(x, y)$ is 
what we call $\mathcal M_{\mathscr C_K(L)}(\bold{x},\bold{y})$ here. 
\\

In all cases of Proposition \ref{rippo}, 
we have defined 
$\mathcal M_{\mathscr C_K(L)}(\bold{x},\bold{y})$ to be 
a finite disjoint union of CW complexes 
which are homeomorphic to the 2-ball 
(Each of these CW complexs is a $<2>$-manifold defined in \cite[Definition 3.1]{LSk}.)
This completes the proof of Proposition \ref{rippo}.
\qed
\\

\h{\bf Remark.}  
Since ${\bf y}=(D,y)$ has just three arcs, 
the number of $e_*^{1+N}$ is $\leqq6$. 
{\it Reason.} One surgery on a labeled resolution configuration 
makes no greater than two labeled resolution configurations. 
Recall Definitions \ref{2.10}. 
Furthermore the number of $e_*^{2+N}$ is $\leqq6$ 
by \cite[
Lemma 2.13]{LSk}.
\\

Suppose that 
$y_1,...y_{\nu_{N-1}}$ are 
all $(-1+N)$-dimensional cells to which 
the $(2+N)$-dimensional cell $x$ 
is attached. 
Here, we also use $x$ and $y_\ast$ for cells. 
For each $i$,  
$\mathcal M(x,y_i)$ 
is embedded in $\partial x$. 
If $i\neq j$, we have 
$\mathcal M(x,y_i)\cap \mathcal M(x,y_i)=\phi$.  

Segments 
$\mathcal M(x,e^{N+1}_\natural)\x\mathcal M(e^{N+1}_\natural,y_\ast)$ 
have been framed before attaching $x$. 
Note that 
$\mathcal M(x,e^{N+1}_\natural)$ is 0-dimensional 
and that 
$\mathcal M(e^{N+1}_\natural,y_\ast)$ is 1-dimensional. 
We give framings to segments $\mathcal M(x,e^N_\sharp)\x\mathcal M(e^N_\sharp,y_\ast)$ 
 so that we can extend a framing on $\partial\mathcal M(x,y_\ast)$ to  $\mathcal M(x,y_\ast)$.  
(A part of 
\cite[Proof Proposition 4.12]{LSk}  
explains how we extends a sign assignment. We use the same way written there.)

\section{3-dimensional modulis  and framings }\label{secmo3}

A {\it handle-body} is a compact oriented 3-manifolds 
with a CW decomposition of one 0-handle and $g$ 1-handles, 
where $g\in\N\cup\{0\}$: $g$ is called the {\it genus}.

\begin{proposition}\label{yon}
Let {\rm gr}$_h \bold x-${\rm gr}$_h \bold y=4$. 
Under the conditions of Propositions $\ref{ten}$, $\ref{hen}$, and $\ref{rippo}$, 
we can assign to 
each moduli space $\mathcal M_{\mathscr C_K(L)}(\bold{x},\bold{y})$,  
a space homeomorphic to a finite number of  handle-bodies. 
\end{proposition}



\h{\bf Proof of Proposition \ref{yon}.}
Assume that 
$\bold x$ 
(respectively, $\bold y$) 
is an $(m+3+N)$-cell   
(respectively, $(m-1+N)$-cell) 
which corresponds to a Khovanov basis  
$(s(D), x)$(respectively, $(D,y)$)  
whose homological grading is 
$m+3$ 
(respectively, $m-1$).
Let $N$ be a fixed large integer. 
We define cells to be closed balls, not open balls.

Take $P(D,x,y)$. 
%
%
Let 
 $g^{m}_i$ (respectively, $g^{m+1}_j$, $g^{m+2}_k$) 
be 
all dual Khovanov basis elements other than $\bold x$ and $\bold y$ 
in  $P(D,x,y)$, 
whose homological gradings are  
$m$ (respectively,  $m+1$, $m+2$).     
We assign to 
 $g^{m}_i$ 
 (respectively,
$g^{m+1}_j$, 
 $g^{m+2}_k$) 
a cell  
$e^{m+N}_i$   
(respectively, 
$e^{m+1+N}_j$,  
$e^{m+2+N}_k$),  
where the right upper suffix of the notation of cells denote the degree. 
Note that $i, j$ and $k$ run over a fixed finite set of natural numbers. 

We assume the following. 

\begin{assumption}\label{proast} 
 We can attach $\bold x$ to $\bold y$ 
to be compatible with the modulis and the framings in 
 Propositions $\ref{ten}$, $\ref{hen}$, and $\ref{rippo}$. 
\end{assumption}

If Assumption \ref{proast} is true, 
  $\partial\mathcal M(\bold x, \bold y)$ 
 satisfies 
\cite[Definition 3.12.(M-2)]{LSk}, which is quoted in \S\ref{secmo2}. 
$\partial\mathcal M(\bold x, \bold y)$ 
is made of 
 $\partial\mathcal M(p,q)$, 
 where  
$p$ (respectively, $q$) is one of 
$\bold x, \bold y$, 
$e^{m+N}_i$,    
$e^{m+1+N}_j$,  and 
$e^{m+2+N}_k$,  
and $(p,q)\neq(\bold x, \bold y)$. 
\\

Note the following: 
Even if Assumption \ref{proast} does not hold, 
we can define a surface $F$ 
which is $\partial\mathcal M(\bold x, \bold y)$ if Assumption \ref{proast} holds. 
$F$ is embedded in the $(m+2+N)$-sphere $\partial \bold x$. 

\begin{cla}\label{claF}
 The above $F$ is orientable.
 \end{cla}
 
 \h{\bf Proof of Claim \ref{claF}.}
Each 
$\mathcal M(\bold x, e^{m+N+1}_\ast)\x\mathcal M(e^{m+N+1}_\ast, \bold y)$ in $F$ 
is a disjoint collection of squares. 
$F$ 
is orientable if and only if \\
$F-$
(all $\mathcal M(\bold x, e^{m+N+1}_\ast)\x\mathcal M(e^{m+N+1}_\ast, \bold y)$)  
is orientable. 

Let 
$A=$($\mathcal M(\bold x, e^{m+2+N}_\sharp)\x\mathcal M(e^{m+2+N}_\sharp, \bold y)$)   \\
and $B=$($\mathcal M(\bold x, e^{m+N}_\ast)\x\mathcal M(e^{m+N}_\ast, \bold y)$).    \\
Then $A\cap B$ 
is 
$C=\mathcal M(\bold x, e^{m+2+N}_\sharp)\x
\mathcal M(e^{m+2+N}_\sharp,  e^{m+N}_\ast)\x\mathcal M(e^{m+N}_\ast, \bold y)$. 
By the rule of \cite[Definition 3.18]{LSk}, we have the following fact. 

\begin{fact}\label{factaiso} 
The framing on $A$ and $C$ are compatible at $C$. 
\end{fact}

Therefore $F$ is orientable. 
This completes the proof of Claim \ref{claF}.  \qed\\

This completes the proof of Proposition \ref{yon}.\qed\\

We prove the following.

\begin{proposition}\label{clanobi}
 Assumption $\ref{proast}$ is true.
   \end{proposition}
\h{\bf Remark.}   Since $\Omega_2^\text{fr}=\Z_2$, 
only 
Proposition \ref{yon}
and Fact \ref{factaiso} never imply Claim \ref{clanobi}. \bb

\h{\bf Proof of Proposition \ref{clanobi}.}    
We fix $l$ and $h$ and consider 
$\mathcal M(e^{m+N+3}_l, e^{m+N-1}_h)$. 
Note that $\bold x=e^{m+N+3}_l$ and $\bold y=e^{m+N-1}_h$. 

A framing on 
$\mathcal M(e^{m+N+3}_l, e^{m+N+2}_k)$
is determined by the differential operator.

A framing on 
$\mathcal M(e^{m+N+2}_l, e^{m+N-1}_j)$ 
is determined when we attach $e^{m+N+2}_l$ to the $(m+N+1)$-skeleton.

Hence 
a framing on 
$\mathcal M(e^{m+N+3}_l, e^{m+N+2}_k)\x \mathcal M(e^{m+N+2}_k, e^{m+N-1}_h)$ 
 in the $(m+2+N)$-skeleton 
has been determined. 
This is a dodecagon or a hexagon.

Here, we consider all possibilities of $h,i,j$, and $k$.
Framings on 0-dimensional modulis \\
$$
\mathcal M(e^{m+N+3}_l, e^{m+N+2}_k)\x
\mathcal M(e^{m+N+2}_k, e^{m+N+1}_j)\x
\mathcal M(e^{m+N+1}_j, e^{m+N}_i)\x
\mathcal M(e^{m+N}_i, e^{m+N-1}_h)
$$
have been given by the differential operator. 

When we attach $e^{m+N+3}_l$, 
 $\mathcal M(e^{m+N+3}_l, e^{m+N+2}_k)$ has been framed, and 
neither 
$\mathcal M(e^{m+N+3}_l, e^{m+N+1}_j)$, \\ 
$\mathcal M(e^{m+N+3}_l, e^{m+N}_i)$, nor 
$\mathcal M(e^{m+N}_i, e^{m+N-1}_h)$ has been framed.

We can change framings on 1-dimensional modulis\\
$$
\mathcal M(e^{m+3+N}_l,  e^{m+1+N}_j)\x
\mathcal M(e^{m+1+N}_j, e^{m+N}_i)\x
\mathcal M(e^{m+N}_i, e^{m-1+N}_h)
$$
and 
$$
\mathcal M(e^{m+3+N}_l, e^{m+2+N}_k)\x
\mathcal M(e^{m+2+N}_k, e^{m+N}_i)\x
\mathcal M(e^{m+N}_i, e^{m-1+N}_h). 
$$

\bb

Recall $F=\partial\mathcal M(e^{m+N+3}_l,e^{m+N-1}_h)$. 
Suppose that a moduli $\mathcal M$ is embedded in $F$. 
\\
Let $X=\amalg_{\text{all $l$}}(\mathcal M(e^{m+N+3}_l, e^{m+N+2}_k)\x \mathcal M(e^{m+N+2}_k, e^{m+N-1}_h))$. 
If 
Int $\mathcal M$ 
is included in 
Int$\displaystyle F-X$, 
then  
Int $\mathcal M$ 
 has not been framed. 
 We will give framing to it.

Give a framing $fr$ on $F$ that extends to $\mathcal M(e^{m+N+3}_l,e^{m+N-1}_h)$.
Each pair of \\
$\mathcal M(e^{m+N+3}_l, e^{m+N+2}_k)\x \mathcal M(e^{m+N+2}_k, e^{m+N-1}_h)$
are disjoint.  
Each $\mathcal M(e^{m+N+3}_l, e^{m+N+2}_k)\x \mathcal M(e^{m+N+2}_k, e^{m+N-1}_h)$ 
is contractible.

All 0-dimensional modulis in $F$ is in $X$. 

Using an isotopy of framings, 
change $fr$ on $F$ so that $fr|X$ coincides with the framing on $X$ that we have given already 
in the lower dimensional skeleton. {\it Reason}: Each component of $X$ is contractible.  
 Restrict this new $fr$ to $F-X$. 

Take any square moduli. 
Only two disjoint 1-dimensional modulis, or edges, are in $X$. 
The other part is in $F-X$. 

For any square moduli $S$, we have the following: 
There is a connected component $E$, a segment,  of 
$\mathcal M(e^{m+N+3}_\ast, e^{m+N+1}_j)$ 
and $F$ of $\mathcal M(e^{m+N+1}_\ast, e^{m+N-1}_h)$ for a suffix $\ast$. 
We have $S=E\x F$. 

Let $\partial E$ be two points $P\amalg P'$. 
Let $\partial F$ be two points $Q\amalg Q'$. 
Then $F\x P$ and $F\x P'$ are in $X$. 
We have Int $E\x Q$ and Int $E\x Q'$ are in $F-X$. 

We must change the framing $fr$ on $F$ so that $fr|S$ 
is the product of a framing on \\
$\mathcal M(e^{m+N+3}_\ast, e^{m+N+1}_j)$ 
 and that on \\ $\mathcal M(e^{m+N+1}_\ast, e^{m+N-1}_h)$ 
 (See \cite[Definition 3.18]{LSk}).

We can do it by using an isotopy of framings because 
$S$ is contractible. 
 
 \bb
Restrict $fr$ to $F-X-\text{all } S$. 
%
%
Therefore there is a framing on the 3-dimensional moduli 
$\mathcal M(e^{m+3+N}_l, e^{m-1+N}_h)$ 
which is compatible with the framing which has been fixed in the  $(m+2+N)$-skeleton.
This completes the proof of Proposition \ref{clanobi}. 
\qed\\

Now we have constructed 
0,1,2, and 3-dimensional modulis 
and framings on them. 
We will use these modulis and frmiangs and 
construct partial Khovanov CW complexes.

\bigbreak
\section{Review of the first Steenrod square operator $Sq^1$ }\label{sq1q}
\hskip10mm

\h
In \cite{Steenrod,SE}
the Steenrod square $Sq^* (*\in\Z)$ is defined.
Let $X$ and $X'$ be compact CW complexes.
Let $\{C_i\}_{i\in\Z}$  be a chain complex. Assume that $\{C_i\}_{i\in\Z}$ is associated with both a CW decomposition  on $X$ and a CW decomposition  on $X'$.
It is well-known that \\$Sq^1(X)=Sq^1(X')$ 
(see e.g. \cite[Introduction]{LSs}) 
and that $Sq^2(X)$ and $Sq^2(X')$ are different in general 
(see e.g. \cite{Seed}). 

Therefore $Sq^1$ is not useful as classical (respectively, virtual) link invariants.
We consider $Sq^2$ below.

\bigbreak
\section{Review of the second Steenrod square operator $Sq^2$ }\label{sq1q}
\hskip10mm

\noindent
We review the definition of the second Steenrod square. 

\begin{definition}\label{secondsq}
{\rm\bf(\cite{Steenrod, SE}.)} 
Let $K$ be the Eilenberg–MacLane space $K(\Z_2,m)$ for any natural number $m>1$. 
We have that $K$ is connected.  
We have that $\pi_i(K)\cong\Z_2$ (respecctively, $0$) 
if $i=m$ (respectively, $i\neq m$ and $ i>1$). 
It is well-known that the homotopy type of $K$ is unique. 
It is known that $H^{m+2}(K;\Z_2)\cong\Z_2$. Let $\xi$ be the generator of 
$H^{m+2}(K;\Z_2)\cong\Z_2$. 

Let $X$ be a CW complex. 
Let $[X,K]$ be the set of all homotopy classes of continuous maps $X\to K$. 
It is well-known that $[X,K]=H^m(X;\Z_2)$. 

For an arbitrary element $x\in H^m(X;\Z_2)$, take $f_x$ by this bijection.
Define the second Steenrod square $Sq^2(x)$ to be $f^\ast_x(\xi)$. 
\end{definition}

This definition is reviewed and explained very well  in \cite[section 3.1]{LSs}.  
See also \cite[The item (4) between Proposition 4L.1 and Theorem 4L.2]{Ha}.
We review an important property of the second Steenrod square operator $Sq^2$.

\begin{proposition}\label{sq2}
{\rm\bf(See \cite[section 12]{Steenrod} and \cite[\S4.L]{Ha}.)} 
Let $Y$ be any compact CW complex. 
Let $Y^{(*)}$ be the $*$-skeleton of $Y (*\in\Z)$. 
Then the second Steenrod square 
$Sq^2(Y):H^{m}(Y;\Z_2)\to H^{m+2}(Y;\Z_2)$ 
is determined by the stable homotopy type of $Y^{(m+2)}/Y^{(m-1)}$. 
\end{proposition}

This proposition is reviewed and explained very well  in \cite[section 3.1]{LSs}.  


We explain how we know $Sq^2(Y)$ from  $Sq^2(Y^{(m+2)}/Y^{(m-1)})$. 
There are maps 


$$
  \begin{CD}
Y @<\text{inclusion}<< Y^{(m+2)} @>\text{quotient}>> Y^{(m+2)}/Y^{(m-1)}
  \end{CD}
$$

By the CW structure of $Y$, that of $Y^{(m+2)}$ and that of $Y^{(m+2)}/Y^{(m-1)}$, 
we have the following commutative diagram (\ref{inyo}). 
Call the homomorphisms as written there. 


{\normalsize

\begin{equation}\label{inyo}
\begin{CD}
      H^{m+2}(Y;\Z_2)  @>\text{injective, $g'$}>> H^{m+2}(Y^{(m+2)};\Z_2)  
@<\text{isomorphism, $f'$}<< H^{m+2}(Y^{(m+2)}/Y^{(m-1)};\Z_2) \\
   @A\text{$Sq^2$}AA   @A\text{$Sq^2$}AA        @A\text{$Sq^2$}AA \\
     H^{m}(Y;\Z_2)   @>\text{isomorphism, $g$}>>  H^{m}(Y^{(m+2)};\Z_2) 
 @<\text{surjective, $f$}<< H^m(Y^{(m+2)}/Y^{(m-1)};\Z_2) \\
\end{CD}
\end{equation}
}
First we explain how we know $Sq^2$ on  $H^{m}(Y^{(m+2)};\Z_2)$ by using 
$Sq^2$ on  \\  
$H^m(Y^{(m+2)}/Y^{(m-1)};\Z_2)$. 
Take any element $x\in H^{m}(Y^{(m+2)};\Z_2)$. 
Take
 $y\\
\in H^m(Y^{(m+2)}/Y^{(m-1)};\Z_2)$ such that $f(y)=x$. 
 By the naturality of $Sq^2$, we have 
$f'(Sq^2(y))=Sq^2(f(y))$. 
So we can know that $Sq^2(x)$ is  $f'(Sq^2(y))$. 
Here, note the following fact: 
Suppose that there is an element \\ $\check{y}\in H^m(Y^{(m+2)}/Y^{(m-1)};\Z_2)$
such that $f(\check{y})=x$ and such that $y\neq\check{y}$. 
Then  $f'(Sq^2(y))=f'(Sq^2(\check{y}))$ although $y\neq\check{y}$. \\

Second we explain how we know $Sq^2$ on  $H^{m}(Y;\Z_2)$ by using 
$Sq^2$ on  $H^m(Y^{(m+2)};\Z_2)$. 
Take any element $z\in H^{m}(Y;\Z_2)$. 
 By the naturality of $Sq^2$, we have 
$g'(Sq^2(z))\\
=Sq^2(g(z))$. Note that $g'$ is injective. 
We can know that $Sq^2(z)$ is ${g'}^{-1}(Sq^2(g(z)))$. \\

Last we explain how we know $Sq^2$ on  $H^{m}(Y;\Z_2)$ by using 
$Sq^2$ on  $H^m(Y^{(m+2)}/Y^{(m-1)};\Z_2)$. 
Take any element $z\in H^{m}(Y;\Z_2)$. 
Take $y\in H^m(Y^{(m+2)}/Y^{(m-1)};\Z_2)$ such that $f(y)=g(z)$. 
So we can know that $Sq^2(z)$ is  ${g'}^{-1}(f'(Sq^2(y)))$. 
Here, note the following fact. \\

\begin{fact}\label{piki}  
Take any element $z\in H^{m}(Y;\Z_2)$. 
Take $y$ in the previous paragraph.   
Suppose that there is an element $\check{y}  
\in H^m(Y^{(m+2)}/Y^{(m-1)};\Z_2)$
such that 
$f(\check{y})=g(z)$ and such that $y\neq\check{y}$. 
Then  ${g'}^{-1}(f'(Sq^2(y))) 
={g'}^{-1}(f'(Sq^2(\check{y})))$ although $y\neq\check{y}$. 

We also have that $Sq^2(\{f^{-1}(g(z))\})$ has only one element. 
\end{fact}

\section{The second Steenrod square operator for virtual links
}\label{sq}

\begin{definition}\label{Zm+2}  
Let $\mathcal L$ be a virtual link. 
Let $L$ be a virtual link diagram which represents $\mathcal L$.
We define a CW complex 
$Z(L)$ 
for $L$ below. 

Let $\{a_p\}_{p\in\Lambda}$ be 
the dual Khovanov basis for $L$.  Note that $\Lambda$ is a finite set.

Let $\{a_p\}_{p\in\Lambda}$ be the dual Khovanov basis for $L$.  
Note that $\Lambda$ is a finite set. 
Take an element $a$ in $\{a_p\}_{p\in\Lambda}$.  
If gr$_ha$ is less than $m$, we assign a base point. 
If gr$_ha$ is greater  than 
$m+4$, 
we assign nothing.

Let 
$g^{m-1}_h$ (respectively, $g^{m}_i$, $g^{m+1}_j$, $g^{m+2}_k$, $g^{m+3}_l$) 
be 
all dual Khovanov basis elements  
whose homological gradings are  
$m-1$ (respectively, $m$, $m+1$, $m+2$, $m+3$)  
in $\{a_p\}_{p\in\Lambda}$.  
We assign to 
$g^{m-1}_h$  
 (respectively,
 $g^{m}_i$, 
$g^{m+1}_j$, 
 $g^{m+2}_k$, $g^{m+3}_l$) 
a cell  
$e^{m-1+N}_h$  
(respectively, 
$e^{m+N}_i$,  
$e^{m+1+N}_j$,  
$e^{m+2+N}_k$, $e^{m+3+N}_l$), 
where the right upper suffix of the notation of cells denote the degree 
and $N$ is a large integer.  Fix $N$.

We attach the cells,  
$e^{m-1+N}_h$, 
$e^{m+N}_i$, 
$e^{m+1+N}_j$, 
$e^{m+2+N}_k$,  and 
$e^{m+3+N}_l$:    
We use the moduli spaces  and the framings 
defined in 
\S\ref{secmo0}-\S\ref{secmo3}.  
The result is $Z(L)$. 
\end{definition} 
\bigbreak

\h{\bf Remark.} 
We consider all possibilities of modulis and framings.  
We may assign to a single virtual link diagram, 
more than one CW complexes if there is a ladybug or a quasi-ladybug configuration. 
See Remark \ref{tsun}. 
\\

Construct $Z(L)$ for each set of modulis and framings. 
We have the following commutative diagram (\ref{sososo}).
The left column is Khovanov chain complex for $L$. 
The middle column is made as follows: 
Remove $C^*(L;\Z_2)  \hskip2mm (*>m+4)$ 
from the left one. 
Put $0$ instead. Let 
$C^{m+4}(L;\Z_2)\to 0$. 

The right column is a {cochain} complex associated with each CW complex $Z(L)$.

\begin{equation}\label{sososo}
  \begin{CD}
      \cdot \\
      \cdot \\
      \cdot \\
  @A\text{$\delta$}AA      \\
     C^{m+4}(L;\Z_2) @>>> 0 @<<< 0\\
@A\text{$\delta$}AA  @A\text{$\delta$}AA @A\text{$\delta$}AA      \\
   C^{m+3}(L;\Z_2)  @>\text{isomorphism}>> C^{m+3}(L;\Z_2) 
@<\text{isomorphism}<< C^{m+3+N}(Z(L);\Z_2) \\
   @A\text{$\delta$}AA   @A\text{$\delta$}AA        @A\text{$\delta$}AA \\
    C^{m+2}(L;\Z_2)  @>\text{isomorphism}>> C^{m+2}(L;\Z_2)  
@<\text{isomorphism}<< C^{m+2+N}(Z(L);\Z_2) \\
   @A\text{$\delta$}AA   @A\text{$\delta$}AA        @A\text{$\delta$}AA \\
    C^{m+1}(L;\Z_2)  @>\text{isomorphism}>>   C^{m+1}(L;\Z_2)  
@<\text{isomorphism}<< C^{m+1+N}(Z(L);\Z_2) \\
 @A\text{$\delta$}AA   @A\text{$\delta$}AA        @A\text{$\delta$}AA \\
   C^{m}(L;\Z_2)   @>\text{isomorphism}>>  C^{m}(L;\Z_2) 
 @<\text{isomorphism}<< C^{m+N}(Z(L);\Z_2) \\
@A\text{$\delta$}AA  @A\text{$\delta$}AA        @A\text{$\delta$}AA \\
     C^{m-1}(L;\Z_2)  @>\text{isomorphism}>> C^{m-1}(L;\Z_2) @<<< 
     C^{m-1+N}(Z(L);\Z_2)\\
@A\text{$\delta$}AA  @A\text{$\delta$}AA        @A\text{$\delta$}AA \\
     C^{m-2}(L;\Z_2)  @>\text{isomorphism}>>    C^{m-2}(L;\Z_2)  @<<< 0\\
  @A\text{$\delta$}AA   @A\text{$\delta$}AA        \\
      \cdot @>>>  \cdot\\ 
      \cdot @>>>  \cdot \\
      \cdot @>>>  \cdot \\
  \end{CD}
\end{equation}
\bigbreak
\bigbreak

We define the second Steenrod square for virtual link diagrams.

In the following definition, recall that there may be many choices of 
a set of modulis. 


\begin{definition}\label{onsen}
Let $L$ be a virtual link diagram. 
Fix a quantum grading. 
Let $x$ be any element in $Kh^m(L;\Z_2)$. 
Define the second Steenrod square $Sq^2(x)\in  Kh^{m+2}(L;\Z_2)$ 
to be 
${\kappa}^{-1}(Sq^2(\{{\pi}^{-1}(x)\}))$ for each $Z(L)$, where 
$Z(L)$ is defined above and 
$\pi$ and $\kappa$ are defined below.

By the commutative diagram (\ref{sososo}) 
we have 
two natural isomorphisms 

  $$\pi:H^{m+N}(Z(L);\Z_2)\to  Kh^m(L;\Z_2)$$ 

\h and 
  
  $$\kappa: Kh^{m+2}(L;\Z_2) \to H^{m+2+N}(Z(L);\Z_2).$$ 

\h
for each $Z(L)$. 
\end{definition} 
\bigbreak
 
 \h{\bf Remark.} 
 For each $m$, $Z(L)$ is a CW complex. 
 Hence $Z(L)$ has the second Steenrod square $Sq^2$.  
By using $Sq^2$ for each $Z(L)$, 
and these two homomorphisms,  
we define the second Steenrod square  $Sq^2$ for $L$ as above.


\section{Our second Steenrod square does not depend on modulis or framings}\label{secdnt}

We have the following.

\begin{thm}\label{daijida}
  $Sq^2$ for $L$ in Definition \ref{onsen} is independent of the choice of modulis and framings. 
  \end{thm}

\h{\bf Proof of Theorem \ref{daijida}.} 
The second Steenrod square of $Z(L)$ is determined by $Z^{(m+2+N)}(L)/Z^{(m-1+N)}(L)$ 
(See \cite[\S3.1]{LSs}).

$Z^{(m+2+N)}(L)/Z^{(m-1+N)}(L)$ are made of 
the cells  
$e^{m+N}_i$,  
$e^{m+1+N}_j$, 
and  
$e^{m+2+N}_k$ 
( Recall 
Definition \ref{Zm+2}.)


\begin{figure}
\includegraphics[width=80mm]{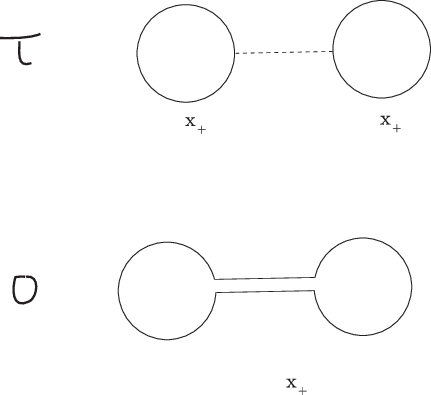}
\caption{{\bf 
Two (resepectively, one) component enhanced Kauffman state
$\tau$ (respectively, $o$)
}\label{figLp}}   
\end{figure}

See the part of 
\cite[Proof of Proposition 6.1]{LSk} corresponding to the following discussion. 
Consider the link diagram $L'$
 obtained by
taking the disjoint union of $L$ with a 1-crossing unknot $U$, 
drawn so that the 0-resolution of $U$ consists of two circles 
(and hence the 1-resolution of $U$ consists of a single circle).
Make  
two (resepectively, one) component enhanced Kauffman state
$\tau$ (respectively, $o$) from $U$ 
as in Figure \ref{figLp}.

Let $P$ be a poset of 
all 
$e^{m+N}_i$,  
$e^{m+1+N}_j$, 
and 
$e^{m+2+N}_k$ of $Z^{(m+2+N)}(L)/Z^{(m-1+N)}(L)$. 

Let $\xi$ be a cell of $P$. Hence $\xi$ corresponds to an enhanced Kauffman state. 
Make posets 
$T=\{\xi\amalg\tau\}$, 
$O=\{\xi\amalg o\}$, and 
$Q=T\cup O$. 

Note the following. 
$[\delta(\tau); o]=1$. \\
$[\delta(\xi\amalg\tau); \xi\amalg o]=1$.  \\
$gr_\text{h}(\xi\amalg\tau)=gr_\text{h}(\xi\amalg o)+1$.  \\  
$gr_\text{q}(\xi\amalg\tau)=gr_\text{q}(\xi\amalg o)$.  

There is an order preserving one-to-one map from $P$ to $O$ (respectively, $T$). 

Consder 
a set $\mathcal S$ (respectively, $\hat{\mathcal S}$)  
of pairs $(\mathcal M(\xi, \xi'), fr)$ (respectively, $(\hat{\mathcal M}(\xi, \xi'), \hat{fr})$) 
of modulis and framings 
for $P$. Constrct a CW complex $W$ (respectively, $\hat W$).

For 
$\xi\amalg\tau\in T$, let 
$(\mathcal M(\xi\amalg\tau, \xi'\amalg\tau), fr)$
$=(\mathcal M(\xi, \xi'), fr)$. 
Make a CW complex $W_o$ from $O$. 

For 
$\hat\xi\amalg o\in O$, let 
$(\hat{\mathcal M}(\hat\xi\amalg o, \hat\xi'\amalg o), \hat{fr})$
$=(\hat{\mathcal M}(\xi, \xi'), \hat{fr})$.

Let 
$\xi\amalg\tau\in T$ and 
$\hat\xi\amalg o\in O$.  
Note 
$\xi\amalg\tau\in T$,  
$\hat\xi\amalg o\in O$ $\in T\cup O=Q$. 
Make a moduli 
$\mathcal M(\xi\amalg\tau, \hat\xi\amalg o)$ 
as in 
\S\ref{secmo0}-\S\ref{secmo3}. 
We can extend framings on 
this moduli $\mathcal M(\xi\amalg\tau, \hat\xi\amalg o)$ 
 to be compatible with those two ones 
 $(\mathcal M(\xi\amalg\tau, \xi'\amalg\tau), fr)$
and 
$(\hat{\mathcal M}(\hat\xi\amalg o, \hat\xi'\amalg o), \hat{fr})$. 
{\it Reason}: 
Let $e\in T$ and $f\in O$. Let $g\in P(e,f)$. 
Let $\text{gr}_h(e)^\text{gr}_h(f)=3$ and  
 $\text{gr}_h(e)^\text{gr}_h(g)=2$. 
Then $\mathcal M(e,g)\subset \partial \mathcal M(e,f)$ has not been framed when we attach $\mathcal M(e,f)$.  We give  a framing on each $\mathcal M(e,g)$ so that 
the framings on $\partial \mathcal M(e,f)$ extends to  $\mathcal M(e,f)$.  

Make a CW complex $W_Q$ from $Q$.

%
%

By the same discussion of the part of \cite[Proof of Proposition 6.1]{LSk}, we have the following: 
We have $W_Q/W_O=\Sigma \hat{W}$ by the construction of these CW complexes. 
By Puppe theorem which is cited in \cite[Lemma3.32.(3)]{LSk}, 
 $W_Q/W_O$ is stable homotopy type equivalent to $\Sigma W$.

Therefore 
$Z^{(m+2+N)}(L)/Z^{(m-1+N)}(L)$ 
does not depend on which moduli and which framings we use.  
\qed



\begin{remark}\label{taiyo}
If we give a framing to modulis of $Z^{(m+2)}(L)/Z^{(m-1)}(L)$ 
 different from the framing  defined in Definition \ref{Zm+2},  
$Sq^2$ does not change (Theorem 
\ref{daijida}).

Although the second Steenrod square of $Z(L)$ is determined by $Z^{(m+2+N)}(L)/Z^{(m-1+N)}(L)$ 
(Proposition \ref{sq2}), 
we do not have  
 $H^{m+N}(Z^{(m+2+N)}(L)/Z^{(m-1+N)}(L);\Z_2)\cong Kh^m(L;\Z_2)$ \\
 and 
$Kh^{m+2}(L;\Z_2) \cong H^{m+2+N}(Z^{(m+2+N)}(L)/Z^{(m-1+N)}(L);\Z_2)$ 
in general (see (\ref{sososo})).
Therefore 
we need to let $\rho\geqq4$ to define $Sq^2$,  
where $\rho$ is 
the difference between 
the largest degree of 
chain groups and the smallest degree of them  
in the right column in 
(\ref{sososo}). 


For $Z(L)$, we make a commutative diagram (\ref{inyo}). 
We can check maps in  (\ref{inyo}) explicitly by using the commutative diagram (\ref{sososo}). Hence we can apply Proposition \ref{sq2} well. 
 We can also apply Fact \ref{piki}.

Furthermore, therefore, we have the following. 
Suppose that $Z(L)$ and $Z^{\bullet}(L)$ are not homotopy type equivalent. 
Assume that $Z^{(m+2)}(L)/Z^{(m-1)}(L)$ and 
$(Z^{\bullet}(L))^{(m+2)}/(Z^{\bullet}(L))^{(m-1)}$ are homotopy type equivalent. 
Then $Z(L)$ and $Z^{\bullet}(L)$  determine the same $Sq^2$ for $L$. 
\end{remark}

 \h{\bf Remark.}
There may be many choices of $Z(L)$. 
Let $Z(L)$ and $Z^{\bullet}(L)$ be such ones. 
Although $Z(L)$ and $Z^{\bullet}(L)$ may not be homotopy type equivalent, 
there is a natural chain isomorphism 
between  $C^*(Z(L);\Z_2)$ and  $C^*(Z^{\bullet}(L);\Z_2)$
by the commutative diagram  (\ref{sososo}). 
 However there is no continuous map between $Z(L)$ and $Z^{\bullet}(L)$ 
 that induces this natural chain  isomorphism, in general.

\bs

Take 
$f$ (respectively, $f', g, g'$) 
for $Z(L)$ in the commutative diagram (\ref{inyo}),  
and call one for   $Z^{\bullet}(L)$, 
$f^{\bullet}$ (respectively, ${f^{\bullet}}', g^{\bullet}, {g^{\bullet}}'$).  
There are two commutative diagrams below. 

{\normalsize
\begin{equation}\label{tomoni}
\begin{CD}
     H^{m}(Z(L);\Z_2)   @>\text{isomorphism, $g$}>>  H^{m}((Z(L))^{(m+2)};\Z_2) 
 @<\text{surjective, $f$}<< H^m((Z(L))^{(m+2)}/(Z(L))^{(m-1)};\Z_2) \\
   @A\text{isomorphism}AA   @A\text{isomorphism}AA        @A\text{isomorphism}AA \\
     H^{m}(Z^{\bullet}(L);\Z_2)   @>\text{isomorphism, $g^{\bullet}$}>>  
     H^{m}((Z^{\bullet}(L))^{(m+2)};\Z_2)  @<\text{surjective, $f^{\bullet}$}<< 
     H^m((Z^{\bullet}(L))^{(m+2)}/(Z^{\bullet}(L))^{(m-1)};\Z_2) \\
\end{CD}
\end{equation}

\begin{equation}\label{kazeto}
\hskip-10mm
\begin{CD}
      H^{m+2}(Z(L);\Z_2)  @>\text{injective, $g'$}>> H^{m+2}(Z(L)^{(m+2)};\Z_2)  
@<\text{isomorphism, $f'$}<< H^{m+2}(Z(L)^{(m+2)}/Z(L)^{(m-1)};\Z_2) \\
   @A\text{isomorphism}AA   @A\text{isomorphism}AA        @A\text{isomorphism}AA \\
      H^{m+2}(Z^{\bullet}(L);\Z_2)  @>\text{injective, ${g^{\bullet}}'$}
      >> H^{m+2}((Z^{\bullet}(L))^{(m+2)};\Z_2)  @<\text{isomorphism, ${f^{\bullet}}'$}<< 
      H^{m+2}((Z^{\bullet}(L))^{(m+2)}/(Z^{\bullet}(L))^{(m-1)};\Z_2) \\
\end{CD}
\end{equation}
}

There is no continuous map between $Z(L)$ and $Z^{\bullet}(L)$ 
 that induces the vertical isomorphisms in 
  (\ref{tomoni})  and (\ref{kazeto}), in general.

Make $Sq^2$ by using each of $Z(L)$ and $Z^{\bullet}(L)$
 according to Definition \ref{onsen}. 
By the commutative diagrams (\ref{tomoni})  and (\ref{kazeto}), 
the difference of each $Sq^2$ occurs by 
the  the most right-handed vertical homomorphism in 
the commutative diagrams (\ref{inyo}) for $Z(L)$ and $Z^{\bullet}(L)$:

{\normalsize

$\begin{CD}
H^{m+2}((Z(L))^{(m+2)}/(Z(L))^{(m-1)};\Z_2) \\
   @A\text{$Sq^2$}AA  \\
   H^m((Z(L))^{(m+2)}/(Z(L))^{(m-1)};\Z_2) \\
\end{CD}$ 
\hskip5mm and \hskip5mm 
$\begin{CD}
H^{m+2}((Z(L)^{\bullet})^{(m+2)}/(Z(L)^{\bullet})^{(m-1)};\Z_2) \\
   @A\text{$Sq^2$}AA  \\
   H^m((Z(L)^{\bullet})^{(m+2)}/(Z(L)^{\bullet})^{(m-1)};\Z_2). \\
\end{CD}$

 }
 \bb

\section{Reidemeister moves do not change our second Steenrod square}\label{secourSq}

\bb

Although  the homotopy type of each $Z(L)$ may change 
by Reidemeister moves on $L$, 
we have the following.  It is a main result of this paper.

\begin{thm}\label{main} 
Let $L$ and $L'$ be virtual link diagrams which represent the same virtual link. 
Note that  $Kh^m(L;\Z_2)\cong Kh^m(L';\Z_2)$. 
Then the second Steenrod square  $Sq^2(L)$ 
is the same as 
$Sq^2(L')$. 
\end{thm} 


\noindent{\bf Proof of Theorem \ref{main}.}
We have the following.
\begin{fact}\label{brIIInashi}
It is enough to prove the case where $L$ is changed into $L'$ by 
a single Reidemeister move. 
The Reidemeister moves are shown in Figure \ref{all-1}.    
\end{fact}

Fix a quantum grading.
Let $\{C^*(L)\}$ (respectively, $\{C^*(L')\}$) be Khovanov chain complex 
for the virtual link diagram $L$  (respectively, $L'$), and   
 $\{C_\#(L)\}$ (respectively, $\{C_\#(L')\}$) the dual Khovanov chain complex 
for $L$  (respectively, $L'$).   
\\

%
%
\noindent{\bf The case of  non-classical Reidemeister moves.}  
Assuming that the Reidemeister move is non-classical, 
then there is a bijective map from 
Khovanov basis 
of $L$ to that of $L'$ 
such that the homological and quantum gradings are kept 
and such that the partial order of the dual Khovanov basis 
are kept.  (See \S\ref{seclo}.).
This identity map induces a chain identity map $C_\#(L)\to  C_\#(L')$.  \\

By using the dual Khovanov basis 
of $L$ (respectively, $L'$) and the moduli spaces, 
construct each $Z(L)$ (respectively, $Z(L')$) 
as in Definition \ref{Zm+2}. 
By using the above bijection from the dual Khovanov basis 
of $L$ to that of $L'$, 
we can make a homeomorphism  
from a disjoinnt union of all $Z(L)$ 
to that of all $Z(L')$, where we give appropriate orders to the components of the disjoint unions.

By this homeomorphism map, 
Theorem \ref{main} holds in this case.  
\\

\noindent{\bf The case of  classical Reidemeister moves.}  
%
If the Reidemeister move is the classical Reidemeister move $I$ or $II$, 
we can prove,  
by the same way as the one in \cite[section 6]{LSk}, 
that  
there is an injective map from 
Khovanov basis 
(respectively, 
the dual Khovanov basis) 
of one of $L$ and $L'$ to that of the other,    
such that the homological and quantum gradings are kept 
and such that the partial order of  
Khovanov basis 
(respectively, the dual Khovanov basis) 
are kept. 
Without loss of generality,  we can assume that this injective map  is from that of $L$ to that of $L'$. 
Furthermore, this injective map induces injective chain homotopy equivalence maps,  $C^*(L')\to  C^*(L)$ and $C_\#(L)\to  C_\#(L')$. \\

Suppose that the Reidemeister move is the classical Reidemeister move $III$.  
Note that both  \cite[section 3.5.5]{B} and \cite[section 6]{LSk} proved that 
there is a chain homotopy equivalence map 
$C^*(L)\to  C^*(L')$ or $C^*(L')\to  C^*(L)$ 
if $L$ and $L'$ are  classical links 
and 
if $L$ is changed into $L'$ by one classical Reidemeister $III$ move. 
The result from \cite[section 6]{LSk} improves on the result in \cite[section 3.5.5]{B}. \\

By the same method of  \cite[section 6]{LSk}, we can prove the following. 
Let $\N$ be the set of natural numbers.  
For $L$ and $L'$ above, 
there are  
virtual link diagrams, $M_1$,..., $M_\mu$ ($\mu\in\N-\{1\}$), 
with the following properties: 
$M_1$ is one of $L$ and $L'$,   
and $M_\mu$ the other. 
If $1\leqq i\leqq\mu-1$, 
$M_i$ is made into $M_{i+1}$ by one classical Reidemeister move $II$ 
or one classical braid-like Reidemeister move $III$ in 
\cite[Figure 6.1]{LSk} and in \cite[section 7.3]{Baldwin}. 
The classical braid-like Reidemeister move $III$ is shown in 
Figure \ref{braid}. \\

\begin{figure}
\includegraphics[width=70mm]{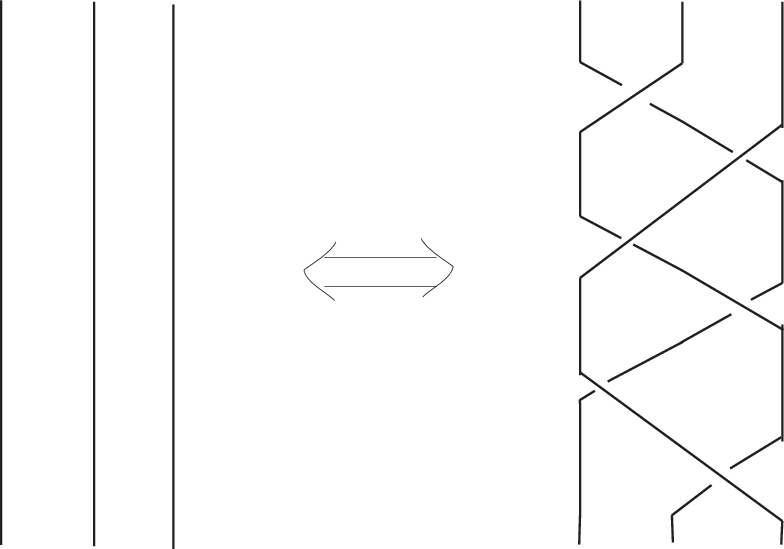}
\caption{{\bf The classical braid-like Reidemeister move $III$
}\label{braid}}   
\end{figure}

Let $M$ be obtained from $M'$ by one classical braid-like Reidemeister move $III$.   
Then there is an injective map from
Khovanov basis 
(respectively, 
the dual Khovanov basis) 
of one of $M$ and $M'$
to that of the other     
such that the homological and quantum gradings are kept 
and such that the partial order of 
Khovanov basis 
(respectively, the dual Khovanov  basis) 
are kept. \\

So we interpret Fact \ref{brIIInashi} as follows. 

\begin{fact}\label{brIII}
It is enough to prove Theorem \ref{main} 
in the case where $L$ is changed into $L'$ by 
a single move. 
The single move is 
each  of 
the classical braid-like Reidemeister move $III$ shown in Figure \ref{braid}  
and 
all other Reidemeister moves  
except for the classical Reidemeister move $III$ shown  in Figure \ref{all-1}.    
\end{fact}

From here, we use the notations $\{C^*(L)\}$,  $\{C^*(L')\}$, $\{C_\#(L)\}$, 
and $\{C_\#(L')\}$ for $L$ and $L'$ in Fact \ref{brIII}.   \\


We proved the case of  non-classical Reidemeister moves 
in the first part of this proof. \\

We prove the case of  classical Reidemeister moves below. \\

In all classical cases
of Fact \ref{brIII}, 
we can assume,  
by the discussion right above,  
that 
we have injective maps from  
the dual Khovanov set 
(repectively, Khovanov set) 
of $L$ to that of $L'$, 
and  
injective chain homotopy equivalence maps \\
$\alpha: C_\#(L)\to  C_\#(L')$ (respectively,  $\beta: C^*(L)\to  C^*(L')$). 
\\

We can prove the following by the same method of \cite{LSk}.

\begin{cla}\label{kaze}
Let 
${\bf x}=(D_L(u),x)$  $($respectively, ${\bf y}=(D_L(v),y))$ 
correspond to 
${\bf x'}=(D_{L'}(u'),x')$  $($respectively, ${\bf y'}=(D_{L'}(v'),y'))$ by this map. 
Then we have the following facts.

  Let $B$ be a 2-dics where the Reidemeister move is carried out.  
After we remove 
arcs in $B$ from ${\bf x'}$ $($respectively, ${\bf y'})$, 
${\bf x}$ and ${\bf x'}$ $($respectively, ${\bf y}$ and ${\bf y'})$ are the same. 
See \cite[Figure 6.2-6.4]{LSk}: 
In these figures, dotted arcs $($their red arcs$)$ are omitted.

$(D_L(u)-D_L(v),x|)$  $($respectively, $(D_L(v)-D_L(u),y|)$$)$  
is the same as \\
$(D_{L'}(u')-D_{L'}(v'),x'|)$  $($respectively, $(D_{L'}(v')-D_{L'}(u'),y'|)$$)$.   

There is a map from 
$P(D_L(v)-D_L(u),x|, y|)$  
to 
$P(D_L(v')-D_L(u'), x'|, y'|)$ 
preserving the homological degree and the quantum one.

The set of all poossible
$\mathcal M({\bf x,y})$   
is the same as 
that of all possible 
$\mathcal M({\bf x',y'})$.   
\end{cla}

\h{\bf Remark.} 
For a fixed pair of $\bf x$ and $\bf y$, 
there are no more than one $\mathcal M({\bf x,y})$
in general.  
When we make modulis in Proposition \ref{rippo}, 
there may be many choices by the existence of quasi-ladybug configurations. 
When we make modulis in Proposition \ref{yon}, 
there may be many choices of homeomorphisms of the boudary of handle bodies.
\\

We use the above injective maps between 
the dual Khovanov basis 
(not Khovanov basis), 
and  
$\alpha: C_\#(L)\to  C_\#(L')$ (not $\beta: C^*(L)\to  C^*(L')$). 
Note that 
this injective chain homotopy equivalence map 
 $\alpha: C_\#(L)\to  C_\#(L')$ 
induces a chain homotopy equivalence map
 $\wp:C^*(L')\to  C^*(L)$ by using the usual ``Hom'' duality.  
Note  that 
 $\wp:C^*(L')\to  C^*(L)$ is not 
$\beta: C^*(L)\to  C^*(L')$,  
and that $\wp$ is not necessarily injective, in general. \\

\section{Sub CW complexes}\label{secsubC}

We are now in a position to construct the needed CW complexes, as we explain below. 
We define each $Z(L)$ (respectively, $Z(L')$) for $L$ (respectively $L'$) 
by using the moduli spaces above. 
There may be more than one $Z(L)$ (respectively, $Z(L')$). 
By the definition of $\alpha$ and 
Claim \ref{kaze}, 
we have the following. 

\begin{cla}\label{daijobu}
Each $Z(L)$ is a sub CW complex of one $Z(L')$ and 
the inclusion map is obtained by using $\alpha$. 
For Each $Z(L')$, there is one $Z(L)$ 
which is a sub CW complex of the  $Z(L')$, and  
the inclusion map is obtained by using $\alpha$. 
\end{cla}

This inclusion map $Z(L)\to Z(L')$ is important for our proof.

We want to prove that 
$Sq^2$ for $Z(L)$ 
and $Sq^2$ for $Z(L')$ 
are the same.

Let $f^\ast_\sharp$ be a cell in $Z(L')$. 
Then  we have only two cases:  
$\text{Int} f^\ast_\sharp$ is in $Z(L)$.  
$\text{Int} f^\ast_\sharp$ is out $Z(L)$.

When we attach $f^\ast_\sharp$, 
framings in the lower dimensional skeleton and in $Z(L)$ are determined. 
Note that this situation is the same as the following situation: 
When we attach $f^\ast_\sharp$, 
framings in the lower dimensional skeleton are determined. 

We proved that we can attach $f^\ast_\sharp$ in this situation. 
Therefore we can define framings on modulis so that we can construct $Z(L')$.

Since 
 $\wp:C^*(L')\to  C^*(L)$ is a chain homotopy equivalence map, 
we have two  isomorphisms 
$$Kh^*(L')\to  Kh^*(L)$$ 
and 
 $$\rho: Kh^*(L';\Z_2)\to  Kh^*(L;\Z_2)$$ 
for all $*$.\\

By using the dual Khovanov basis elements 
of $L'$ and the moduli spaces,  
construct $Z(L')$ as in Definition \ref{Zm+2}. 
By using the above injective map 
from the dual Khovanov basis 
of $L$  
to that of $L'$,    
we can let one $Z(L)$ be a sub-CW complex of one $Z(L')$ 
as written above.  
Hence there is an inclusion map $Z(L)\to Z(L')$. \\

By the diagram (\ref{sososo}) 
we have the following two commutative diagrams.

\np
\begin{equation}\label{kakaka}
  \begin{CD}
C^{m+4}(L;\Z_2)  
 @<<< 0\\
  @A\text{$\delta$}AA @A\text{$\delta$}AA      \\
C^{m+3}(L;\Z_2)  
@<\text{isomorphism}<< C^{m+3+N}(Z(L);\Z_2) \\
  @A\text{$\delta$}AA        @A\text{$\delta$}AA \\
C^{m+2}(L;\Z_2)  
@<\text{isomorphism}<< C^{m+2+N}(Z(L);\Z_2) \\
    @A\text{$\delta$}AA        @A\text{$\delta$}AA \\
C^{m+1}(L;\Z_2)  
@<\text{isomorphism}<< C^{m+1+N}(Z(L);\Z_2) \\
   @A\text{$\delta$}AA        @A\text{$\delta$}AA \\
 C^{m}(L;\Z_2) 
 @<\text{isomorphism}<< C^{m+N}(Z(L);\Z_2) \\
  @A\text{$\delta$}AA @A\text{$\delta$}AA      \\
C^{m-1}(L;\Z_2)  
@<\text{isomorphism}<< 
C^{m-1+N}(Z(L);\Z_2) \\
  @A\text{$\delta$}AA        @A\text{$\delta$}AA \\
  C^{m-2}(L;\Z_2)  @<<< 0\\
  \end{CD}
\end{equation}

\bigbreak

\np
\begin{equation}\label{kakaka}
  \begin{CD}
C^{m+4}(L';\Z_2)  
 @<<< 0\\
  @A\text{$\delta$}AA @A\text{$\delta$}AA      \\
C^{m+3}(L';\Z_2)  
@<\text{isomorphism}<< C^{m+3+N}(Z(L');\Z_2) \\
  @A\text{$\delta$}AA        @A\text{$\delta$}AA \\
C^{m+2}(L';\Z_2)  
@<\text{isomorphism}<< C^{m+2+N}(Z(L');\Z_2) \\
    @A\text{$\delta$}AA        @A\text{$\delta$}AA \\
C^{m+1}(L';\Z_2)  
@<\text{isomorphism}<< C^{m+1+N}(Z(L');\Z_2) \\
   @A\text{$\delta$}AA        @A\text{$\delta$}AA \\
 C^{m}(L';\Z_2) 
 @<\text{isomorphism}<< C^{m+N}(Z(L');\Z_2) \\
  @A\text{$\delta$}AA @A\text{$\delta$}AA      \\
C^{m-1}(L';\Z_2)  
@<\text{isomorphism}<< 
C^{m-1+N}(Z(L');\Z_2) \\
  @A\text{$\delta$}AA        @A\text{$\delta$}AA \\
  C^{m-2}(L';\Z_2)  @<<< 0\\
  \end{CD}
\end{equation}

\bigbreak

\np
The following is a part of the chain homotopy equivalence map $\wp:C^*(L')\to  C^*(L)$. 

\begin{equation}\label{kiki}
  \begin{CD}
C^{m+3}(L';\Z_2)  @>>>  C^{m+3}(L;\Z_2) \\
  @A\text{$\delta$}AA        @A\text{$\delta$}AA \\
C^{m+2}(L';\Z_2)  @>>>  C^{m+2}(L;\Z_2) \\
    @A\text{$\delta$}AA        @A\text{$\delta$}AA \\
C^{m+1}(L';\Z_2)  @>>>  C^{m+1}(L;\Z_2) \\
   @A\text{$\delta$}AA        @A\text{$\delta$}AA \\
C^{m}(L';\Z_2)  @>>>  C^{m}(L;\Z_2) \\
  @A\text{$\delta$}AA        @A\text{$\delta$}AA \\
C^{m-1}(L';\Z_2)  @>>>  C^{m-1}(L;\Z_2) \\
  @A\text{$\delta$}AA        @A\text{$\delta$}AA \\
C^{m-2}(L';\Z_2)  @>>>  C^{m-2}(L;\Z_2) \\
\end{CD}
\end{equation}
\bigbreak


By using the commutative diagrams 
(\ref{kakaka})-(\ref{kiki}), 
we make the following

\begin{equation}\label{kuku}
  \begin{CD}
 0 @>>> 0\\
  @A\text{$\delta$}AA @A\text{$\delta$}AA      \\
C^{m+3+N}(Z(L');\Z_2)  @>>> C^{m+3+N}(Z(L);\Z_2) \\
  @A\text{$\delta$}AA        @A\text{$\delta$}AA \\
C^{m+2+N}(Z(L');\Z_2)  @>>> C^{m+2+N}(Z(L);\Z_2) \\
    @A\text{$\delta$}AA        @A\text{$\delta$}AA \\
C^{m+1+N}(Z(L');\Z_2)  @>>> C^{m+1+N}(Z(L);\Z_2) \\
   @A\text{$\delta$}AA        @A\text{$\delta$}AA \\
C^{m+N}(Z(L');\Z_2)  @>>> C^{m+N}(Z(L);\Z_2) \\
  @A\text{$\delta$}AA        @A\text{$\delta$}AA \\
C^{m-1+N}(Z(L');\Z_2)  @>>> C^{m-1+N}(Z(L);\Z_2) \\
  @A\text{$\delta$}AA @A\text{$\delta$}AA      \\
0  @>>> 0\\
  \end{CD}
\end{equation}

\noindent
that  induces
%
%
%
%
two isomorphisms 
  $$\zeta: H^{m+N}(Z(L');\Z_2)\to H^{m+N}(Z(L);\Z_2)$$ 
and 
$$\theta: H^{m+2+N}(Z(L');\Z_2)\to H^{m+2+N}(Z(L);\Z_2).$$

\bb

We have the following two commutative diagrams by 
the construction of $Z(L)$, 
that of $Z(L')$,  
the definition of $\wp$, 
and the commutative diagrams (\ref{sososo}) and  (\ref{kakaka})-(\ref{kuku}).

The homomorphism 
$   H^{m+N}(Z(L');\Z_2)  \to  H^{m+N}(Z(L);\Z_2)$
is an epimorohism by the construction of $Z(L)$ and that of $Z(L')$.

$$
  \begin{CD}
      Kh^m(L';\Z_2) @>\text{isomorphism, $\rho$}>> Kh^m(L;\Z_2) \\
  @A\text{isomorphism, $\pi$}AA          @A\text{isomorphism, $\pi$}AA \\
     H^{m+N}(Z(L');\Z_2)   @>\text{
isomorohism,      $\zeta$}>>  H^{m+N}(Z(L);\Z_2)
  \end{CD}
$$
\vskip10mm

$$\begin{CD}
H^{m+2+N}(Z(L');\Z_2)   @>\text{isomorphism, $\theta$}>>  H^{m+2+N}(Z(L);\Z_2)\\
  @A\text{isomorphism, $\kappa$}AA      @A\text{isomorphism, $\kappa$}AA \\
    Kh^{m+2}(L';\Z_2) @>\text{isomorphism, $\rho$}>> Kh^{m+2}(L;\Z_2) \\
  \end{CD}
$$

\vskip10mm
The naturality of Steenrod squares makes the following commutative diagram.

$$
  \begin{CD}
     H^{m+N}(Z(L');\Z_2)   @>\text{
isomorphism,     $\zeta$}>>  H^{m+N}(Z(L);\Z_2)\\
  @V Sq^2 VV    @V Sq^2 VV    \\
     H^{m+2+N}(Z(L');\Z_2)   @>\text{isomorphism, $\theta$}>>  H^{m+2+N}(Z(L);\Z_2)
  \end{CD}
$$
\vskip10mm

Combine these three commutative diagrams to obtain a commutative diagram below. 

$$
  \begin{CD}
      Kh^m(L';\Z_2) @>\text{isomorphism, $\rho$}>> Kh^m(L;\Z_2) \\
  @A\text{isomorphism, $\pi$}AA           @A\text{isomorphism, $\pi$}AA        \\
   H^{m+N}(Z(L');\Z_2)   @>\text{
isomorphism,   $\zeta$}>>  H^{m+N}(Z(L);\Z_2)\\
  @V Sq^2 VV    @V Sq^2 VV    \\
   H^{m+2+N}(Z(L');\Z_2)   @>\text{isomorphism, $\theta$}>>  H^{m+2+N}(Z(L);\Z_2)\\
@A\text{isomorphism, $\kappa$}AA      @A\text{isomorphism, $\kappa$}AA \\
    Kh^{m+2}(L';\Z_2) @>\text{isomorphism, $\rho$}>> Kh^{m+2}(L;\Z_2) \\
  \end{CD}
$$
\vskip10mm

By this diagram, Theorem \ref{main} holds. \qed\\

By Theorem \ref{main}, the following is well-defined.

\begin{definition}\label{ofuro}
Let $\mathcal L$ be a virtual link. 
Let $L$ be a virtual link diagram which represents $\mathcal L$. 
We define the second Steenrod square 
$Sq^2$ for $\mathcal L$ 
to be  that on $Sq^2$ for $L$. 
\end{definition}


The following is a main result of this paper.

\begin{thm}\label{YY}
If $\mathcal L$ is a classical link, 
the second Steenrod square $Sq^2$ 
specifies 
the second Steenrod square 
in the case of classical links which is defined in \cite{LSs}.
\end{thm}

\noindent{\bf Proof of Theorem \ref{YY}.} 
By 
Theorem \ref{daijida} and 
the definition of  the second Steenrod square of classical links in \cite{LSk}. 
\qed\\

Note that 
Theorem \ref{main} includes  Theorem \ref{YY}. 
Note Remark \ref{remiron}. 

In the case of classical links,  
in \cite{LSs} Lipshitz and Sarkar showed a way to calculate $Sq^2$ 
by using classical link diagrams. 
We must consider whether the method can be extended to the case of virtual links.
\\

Recall Theorem \ref{LSSe}: 
There are classical links $K_1$ and $K_2$ 
such that 
Khovanov homology of $K_1$ is the same as that of $K_2$, 
but that the second Steenrod square  of $K_1$ is different from that of $K_2$. \\

Note that all classical links are virtual links by the definition. 
%
Therefore, by Theorem  \ref{LSSe}, 
there are virtual links $K_1$ and $K_2$ 
such that 
Khovanov homology of $K_1$ is the same as that of $K_2$, 
but that the second Steenrod square  of $K_1$ is different from that of $K_2$. \\
\\

It is very natural to ask whether 
there is a pair of non-classical, virtual links 
which have different second Steenrod squares 
and which have the same Khovanov homology. 
We answer this question below. 
This is a main result of this paper.\\

\begin{thm}\label{oreikhalkos}
In the case of non-classical, virtual links, our second Steenrod square $Sq^2$ 
is stronger than Khovanov homology. 
That is, there is a pair of non-classical, virtual links 
which have different second Steenrod squares  
and which have the same Khovanov homology. 
\end{thm}

\noindent{\bf Proof of Theorem \ref{oreikhalkos}.}  
Let us consider the example in Figure \ref{reiX}.
Let  $K$ be any classical link diagram. 
Then this represents a non-classical, virtual  link. 
{\it Reason.} The virtual link diagram in Figure \ref{reiX} has only one virtual crossing point. 
No Reidemeister move changes the parity of virtual crossing points 
which are made by two different components of a virtual link diagram. \\

In this case, we do not have a quasi-ladybug configuration. 
The right pair and the left one of ladybug situations 
give the same Steenrod second square 
by explicit calculus which uses that about the classical link diagram $K$.  
(Note that the Steenrod square is only one element in this case.) \\ 

Take $K_1$ and $K_2$ in Theorem \ref{LSSe}, cited above.   
Let $K$ be $K_1$ (respectively, $K_2$).   
These two cases have different Steenrod squares and the same Khovanov homology. 
Recall that Seed used Lipshitz and Sarkar's work to produce two classical knots $K_1$ and $K_2$
 with same Khovanov homology that are distinguished by using the Steenrod square. 
 We observe that, by using  $K_1$ and $K_2$ in our virtual examples, 
 we can achieve the same phenomenon.
\\

\h{\bf An alternative proof:}  
The operation 
flanking (See \cite{DKK}) 
of each of Seed's pair of classical knots \cite{Seed} 
makes 
infinitely many pair that satisfy the condition of Theorem \ref{oreikhalkos}.   
(Two operations, flanking and virtualization, are different although they are very related. See \cite{DKK, Ru}.)
\qed\\

The above examples are just the beginning of many possible applications of the result in this paper. 
Further applications require deeper computations of the virtual Khovanov homology 
and will be the subject of a subsequent paper.\\

\begin{figure}
\includegraphics[width=60mm]{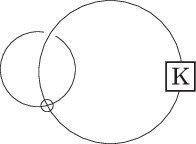}
\caption{{\bf 
If $K$ is any classical link diagram, 
this virtual link diagram represents a non-classical, virtual link.  
}\label{reiX}}   
\end{figure}

\bigbreak

In \cite{Ru}, Rushworth gives another way to  construct a Khovanov homology for virtual links, called the doubled Khovanov Homology. It is natural to ask whether we can make a CW complex and Steenrod squares for the doubled Khovanov chain complex. We hope to answer this question in later work.

\bigbreak
\section{Open problems 
} 
\label{tuke}

\noindent
One reason why virtual knots and links are introduced is as follows: 
The definition of the Alexander polynomial and those of many other invariants 
for links in $S^3$ are extended into the case of links in any 3-manifold other than $S^3$
consistently. It is very natural to ask the following question.

\begin{que}\label{Jones}
Can we generalize the definition of the Jones polynomial for classical links in $S^3$ 
to the case of any 3-manifold?  
\end{que}

A significant partial answer to this question  
is given by using virtual knot theory, for the case of manifolds of the form of a thickened surface.
The Jones polynomial for virtual knots is defined in
 \cite{Kauffman1,Kauffman, Kauffmani}.   
This fact means that , by using virtual 1-knots, we can define the Jones polynomial for links in 
(the closed oriented surface)$\x[-1,1]$.  
The case of knots in the 3-ball $B^3$ and that in the 3-space $\R^3$ 
are trivially the same as that in $S^3$. We omit comments about 
these cases and such similar other trivial cases ($\R^3-$(an open 3-balls), etc.). 
Furthermore see Remark \ref{!} below (see also \cite{KOS}). \\

Note also that for knots and links in thickened surfaces, taken up to handle stabilization, virtual knot theory has a fully diagrammatic formulation, and it can be studied by using
the minimal embedding genus for the virtual knot or link. This gives the theory a flexibility that has led to the discovery of many new invariants of virtual links and relationships with classical knot theory. One finds that by using the generalization of classical knot theory to virtual knot theory, there are infinitely many non-trivial virtual knots with unit Jones polynomial, all occurring in higher genus surfaces so far (see \cite{DKK}). 
Virtual knot theory is a context for studying the conjecture that the Jones polynomial detects the classical unknot.  
(Recall the following facts. 
Let $K$ be a classical knot diagram for  a classical knot. 
In \cite{H}, Haken introduced an algorithm which 
detects whether $K$ represents the classical unknot or not.
After that, in \cite{KM}, 
Kronheimer and Mrowka proved   
the Khovanov homology for classical knots can 
detect that. 
After that, in \cite{OSg}, 
Ozsv\'ath and Szab\'o proved that 
their introduced knot Floer homology can detect that.) 
\\

\begin{figure}
\includegraphics[width=40mm]{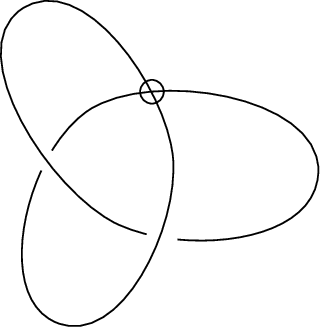}    
\caption{{\bf The Jones polynomial of this virtual knot is not that of any classical knot. }\label{vtre}}   
\end{figure}

There is a virtual 1-knot whose Jones polynomial is not that of any classical knot. 
An example is shown in Figure \ref{vtre}. 
Therefore there is a knot in a thickened surface 
whose Jones polynomial is not that of any  knot in $S^3$. 
We can say that the Jones polynomial of knots in a thickened surface 
defined by using virtual knots is 
not that of knots made by the following way (i) nor that by (ii). 

\bs\h(i) 
Let $K$ be a knot in a thickened surface $F\x[-1,1]$. 
Make $K$ into $K'$: 
Embed $F\x[-1,1]$ in $S^3$, for example, in the standard position. 
We obtain a new knot $K'$ in $S^3$. 

\bs\h(ii)  
Let $K$ be a knot in a thickened surface $F\x[-1,1]$. 
Make $K$ into $K'$: 
Take the universal covering space, which is $\R^2\x[-1,1]$,  of 
$F\x[-1,1]$. Lift $K$ to the universal covering space. 
Note that the lift of $K$ has many components in general. 
Take one of them if there are many.
Thus we obtain $S^1$ or $\R$ embedded in $\R^2\x[-1,1]$. 
We let it be made into  a new knot $K'$ in $S^3$.

\begin{remark}\label{!}
In \cite[Theorem 3.3.3, page 560]{RT} there are defined invariants for links in a closed oriented 
3-manifold $M$. These invariants depend on the use of the colored Jones polynomials at roots of unity and the use of the Kirby calculus (see the last paragraph of this note).
We are interested in more direct constructions for invariants of links in three manifolds, and we believe that the formulation of virtual knot theory is
a step in this direction. Similarly in \cite {W} Witten formulated his functional integral for any link in $S^3$  at physics level. 
This gives a heuristic three dimensional
definition of specializations of the Jones polynomial. 
However, Question \ref{Jones} is open even in physics-level.  
Witten's path integral has not been calculated explicitly in the case of links in 
all 3-manifolds.   
Another aim for combinatorial topology is to find rigorous combinatorial definitions for invariants such as the Jones polynomial.
Here we do not claim that new insight is gained from virtual knot theory. But now we can further ask for fully three dimensional definitions of the extension of the Jones polynomial that we have defined in the virtual theory.

Let $L$ be a link in $S^3$. 
Let $RT(L)$ denote  all invariants defined in \cite[Theorem 3.3.3, page 560]{RT}. 
Note that we consider the case $M=S^3$.   
It is an open question that which is stronger, $RT(L)$ or the Jones polynomial of $L$. 
Furthermore, even if $RT(L)$ is no weaker than 
the Jones polynomial of $L$,  
then do we know the Jones polynomial explicitly 
by a finite times of algorithm which uses a piece of information of $RT(L)$?

 \end{remark}
\vskip3mm

We also have very natural outstanding open questions below.

\begin{que}\label{Khov}
\noindent$(1)$
Can we generalize the definition of the Khovanov homology for links in $S^3$
to the case of any 3-manifold?  

\noindent$(2)$
 Can we generalize the definition of 
the second Steenrod square operator on the Khovanov homology for links in $S^3$ 
to the case of any 3-manifold?

\noindent$(3)$
Can we generalize the definition of 
the Khovanov stable homotopy type for links in $S^3$ to the case of any 3-manifold?  
\end{que}




In \cite{APS},  
 Asaeda,  Przytycki, and Sikora gave 
a partial answer  to Question \ref{Khov}.(1):  
It is given in the case of links in thickened surfaces. 

In \cite{Man}, Manturov introduced the Khovanov homology for virtual links. 
It is a partial answer to Question \ref{Khov}.(1):  
It is the case of links in thickened surfaces by using virtual knot theory.

In \cite{Ru}, 
Rushworth introduced the Khovanov homology for virtual links 
in a different way from that in \cite{Man}. 
It is a partial answer to Question \ref{Khov}.(1):  
It is the case of links in thickened surfaces by using virtual knot theory.

See also Tubbenhauer \cite{Tub} and  Viro \cite{Viro}. 

In \cite{MN}, Manturov and Nikonov made 
an alternative definition of that in \cite{APS}, 
and obtained a new result by using it. 

In \cite{DKK}, 
Dye, Kaestner, and Kauffman gave an alternative definition of that of \cite{Man}, 
and obtained a new result by using it.  

In \cite{Igor}, Nikonov described an alternative definition of integral virtual Khovanov homology  of \cite{Man}. 
Nikonov's definition is written elegantly and explicitly. 
\\

In this paper we define 
the second Steenrod square operator on 
the Khovanov homology for virtual links. 
Therefore  
we can define  the second Steenrod square operator on 
the Khovanov homology 
for links in (the closed oriented surface)$\x[-1,1]$,   
and give a partial answer to Question \ref{Khov}.(2).  
This is only one consistent partial answer to Question \ref{Khov}.(2),  for now.  
This is the main result of this paper (Main theorem \ref{ogon}). \\

We do not give an answer to Question \ref{Khov}.(3) in this paper. 
We give the partial solution to Question \ref{Khov}.(2)   toward 
answering Question \ref{Khov}.(3) in the future. 

\vs

There are found many relations between 
  the Khovanov homology for knots in $S^3$  
  and  knot Floer homology for knots in $S^3$.  
See  \cite{MOST, OSS, OSkf, Ras}
for knot Floer homology. 
See \cite{OS3, OSap} for 
the Heegaard Floer homology, 
from which knot Floer homology is made. 
See  \cite{Baldwin, BSS, MO} and  \cite[section 1.4]{OSS} etc. for their relations.

Since we extended 
the definition of the Jones polynomial for knots in $S^3$,   
that of the Khovanov homology for them,  and  
that of the Steenrod square acting on the Khovanov homology for them  
to the case of thickened surfaces, 
it is very natural to ask questions below.

\begin{que}\label{FQ}
(1) 
Can we extend knot Floer homology for knots in $S^3$ to 
the case of thickened surfaces?  
Is there an invariant for knots in thickened surfaces which is made 
from knot Floer homology or which is much related to  knot Floer homology? 

\bs\h(2)
If we can make an invariant in the question  right above, 
is there a relation between the new invariant and 
the Khovanov homology for knots in thickened surfaces? 

\bs\h(3) 
If we can make an invariant in the above question (1), 
is it virtual knot invariant? 
\end{que}

Note that, by using virtual knot theory, 
we can define the Alexander-Conway polynomial for 
 knots in thickened surfaces 
even if knots are non-vanishing cycles in thickened surfaces 
 (\cite{Kauffman1,Kauffman, Kauffmani}).
\\

\begin{remark}\label{remFloer}
An observation to Question \ref{FQ}. 
The following  is an invariant for knots in thickened surafces, 
which is made by using knot Floer homology. 
Let $F$ be a closed oriented surface.
Let $i=1,2$. Let $K_i$ be a null-homologous 1-knot in $F\x[-1,1]$. 
Regard $F\x S^1$ as the double of $F\x[-1,1]$. 
Note that $K_i\subset F\x[-1,1]\subset F\x S^1$, and   
call this knot in $F\x S^1$, $K'_i$. 
If $K_1\subset F\x[-1,1]$ is obtained 
from 
$K_2\subset F\x[-1,1]$
by \\
(a diffeomorphism of $F$)$\x$(the identity map of $[-1,1]$), 
 $K'_1\subset F\x S^1$ is obtained 
from 
$K'_2\subset F\x S^1$ 
by 
(a diffeomorphism of $F$)$\x$(the identity map of $S^1$). 
  If $K_1\subset F\x[-1,1]$ is obtained 
from 
$K_2\subset F\x[-1,1]$
by a classical move, 
 $K'_1\subset F\x S^1$ is obtained 
from 
$K'_2\subset F\x S^1$
by a classical move. 

It is important that knot Floer homology of null-homologous knots 
in $F\x S^1$ is defined in 
\cite{OSkf}. 
Knot Floer homology of $K'_i\subset F\x S^1$ 
gives an invariant of  $K_i\subset F\x[-1,1]$. 

We could prove that, by using \cite[Theorem 1.1]{OSg}, 
not all invariant of this kind is a virtual knot invariant. 

For 3-manifolds with non-vacuous boundary, 
some Floer homologies are defined:  
 bordered Heegaard Floer homology 
in \cite{borderpaper, borderbook} and 
sutured Floer homology in \cite{J}.  
Can we obtain an invariant of links in thickened surfaces by 
using them?
\end{remark}

\h\noindent{\bf Remark.} 
It is trivial that there is a meaningful map from 
 the set $KI$ of knots in  $F\x[-1,1]$
 to that $KS$ of ones in  $F\x S^1$, 
 but that 
 it is very difficult to make 
 a `meaningful' map from $KS$ to $KI$. 
So 
we can 
obtain an invariant of knots in $F\x[-1,1]$ 
by using 
knot Floer homology for null-homologous knots in $F\x S^1$. 
We know 
the Jones polynomial (resp. the Khovanov homology, the Steenrod square) 
for knots in $F\x[-1,1]$, 
but have not known that of  knots in $F\x S^1$ 
 (Recall Questions \ref{Jones} and \ref{Khov}.). 
\bb

\begin{remark}\label{remato}
 (1) After submitting this paper to arXiv, 
and before publishing this paper in this journal,  
Kauffman, Nikonov, and Ogasa \cite{KauffmanNikonovOgasa, KauffmanNikonovOgasaT2}
 constructed the Khovanov homotpy type for links in thickened surfaces. 
 It is the first partial consistent answer to  Question \ref{Khov}.(3) and the second  partial consistent answer to Question \ref{Khov}.(2).  
\bb
\h(2)
Virtual links are not only generalizations of a classical links, 
but also 
Virtual links introduce 
new topological quantum invariants of links in the 3-sphere, 
that is, classical links.  
See Kauffman and Ogasa's paper \cite{KauffmanOgasaquantum}. 
In the paper, they also defined new topological quantum invariants of 3-manifolds with non-vacuous boundary. 
Dye, Kauffman, and Ogasa's paper 
\cite{DKO}  
is a sequel of 
\cite{KauffmanOgasaquantum}. 
\bb
\h(3)
After submitting this paper to arXiv, 
and before publishing this paper in this journal,  
Juhász, Kauffman, and Ogasa \cite{JKO} 
generalized an idea in \ref{remFloer} and 
introduced a Floer homology for all knots (not only null-homologous knots but  also 
knots that are non-vanishing cycles) in thickened surfaces and all virtual knots. 
\end{remark}

\noindent
{\bf Acknowledgment.}  
The authors would like to thank 
Igor Mikhailovich Nikonov for the valuable discussion.

Kauffman's work was supported by the 
Laboratory of Topology and Dynamics, 
Novosibirsk State University 
(contract no. 14.Y26.31.0025 
with the Ministry of Education and Science 
of the Russian Federation.)

\np
\vskip8mm
\noindent
Louis H. Kauffman

\noindent
Department of Mathematics, Statistics and Computer Science

\noindent
University of Illinois at Chicago

\noindent
851 South Morgan Street

\noindent
Chicago, Illinois 60607-7045

\noindent
USA

\noindent
and

\noindent
Department of Mechanics and Mathematics

\noindent
Novosibirsk State University

\noindent
Novosibirsk

\noindent
Russia

\noindent
kauffman@uic.edu
\\

\noindent
Eiji Ogasa

\noindent
Meijigakuin University, Computer Science 
 
\noindent
Yokohama, Kanagawa, 244-8539 

\noindent
Japan 

\noindent
pqr100pqr100@yahoo.co.jp  

\noindent
ogasa@mail1.meijigkakuin.ac.jp 

\end{document}